\documentclass[11pt]{article}

\usepackage{setspace}
\usepackage[utf8]{inputenc}

\usepackage{bbold}
\usepackage{dsfont}
\usepackage{mathtools}
\usepackage{stmaryrd}
\usepackage{verbatim}
\usepackage{amsthm}
\usepackage{amssymb}
\usepackage{titlesec}
\usepackage{mathrsfs}

\usepackage{float}
\usepackage{amsfonts}
\usepackage{comment}
\usepackage{soul}

\usepackage{lipsum}
\usepackage{mhchem}
\usepackage{graphicx}
\graphicspath{ {images/} }
\usepackage{array}
\newcolumntype{P}[1]{>{\centering\arraybackslash}p{#1}}
\usepackage{longtable}
\usepackage{mathrsfs}
\usepackage{caption}
\DeclareCaptionLabelFormat{adja-page}{\hrulefill\\#1 #2 \emph{(previous page)}}
\usepackage{subcaption}
\usepackage{xcolor}
\allowdisplaybreaks

\usepackage[noline]{algorithm2e}
\RestyleAlgo{tworuled}

\usepackage{soul}

\usepackage{cite}
\usepackage{tikz}

\usepackage{tabulary}
\usepackage{booktabs}
\usepackage{array} 
\usepackage[flushleft]{threeparttable}
\usepackage{blkarray}

\usepackage{hyperref}
\hypersetup{
	colorlinks=true,
	linkcolor=blue,
	filecolor=blue,      
	urlcolor=blue,
	citecolor=blue,
	pdftitle={Overleaf Example},
	pdfpagemode=FullScreen,
}

\usepackage[toc,page,titletoc]{appendix}
\usepackage[capitalise,nameinlink]{cleveref}

\crefname{supp}{Supplement}{Supplements}

\usepackage{amsthm}
\newtheorem{theorem}{Theorem}[section]

\newtheorem{lemma}{Lemma}[section]

\newtheorem{assumption}{Assumption}[section]
\numberwithin{equation}{section}

\theoremstyle{definition}
\newtheorem{definition}{Definition}[section]

\theoremstyle{definition}
\newtheorem{remark}{Remark}[section]

\theoremstyle{definition}
\newtheorem{example}{Example}[section]

\usepackage{authblk}

\title{Comparison Theorems for Stochastic Chemical Reaction Networks}
\author[1,*]{Felipe A. Campos}
\author[2,*]{Simone Bruno}
\author[1]{Yi Fu}
\author[2]{Domitilla Del Vecchio}
\author[1]{Ruth J. Williams}
\affil[1]{Department of Mathematics, University of California, San Diego, 9500 Gilman Drive, La Jolla CA 92093-0112. Email: {\tt\small (fcamposv,yif064,rjwilliams)@ucsd.edu}}
\affil[2]{Department of Mechanical Engineering, Massachusetts Institute of Technology, 77 Massachusetts Avenue, Cambridge, MA 02139. Emails: {\tt\small (sbruno,ddv)@mit.edu}}
\affil[*]{These authors contributed equally: F. A. Campos and S. Bruno}
\date{}                     
\renewcommand{\l}{\ell}

\newcommand{\eps}{\varepsilon}
\newcommand{\intervalx}{(K_A + x) \cap \X}

\newcommand{\initialx}{x^\circ}
\newcommand{\initialxbreve}{\Breve{x}^\circ}

\newcommand{\Ss}{\mathscr{S}}
\newcommand{\Rs}{\mathscr{R}}

\newcommand{\PP}{\mathds{P}}
\newcommand{\R}{\mathds{R}}

\newcommand{\Z}{\mathds{Z}}
\newcommand{\E}{\mathds{E}} 

\newcommand{\D}{\mathcal{D}}
\newcommand{\F}{\mathcal{F}}

\renewcommand{\L}{\mathcal{L}}

\newcommand{\prp}{\Lambda}  
\newcommand{\rate}{\Upsilon} 

\newcommand{\X}{\mathcal{X}}

\newcommand{\V}{\mathcal{V}}

\newcommand{\one}{\mathds{1}}

\newcommand{\Dtot}{\text{D}_{\text{tot}}}

\newcommand{\PPhi}{\Psi}
\newcommand{\Kb}{\textbf{K}}

\DeclareMathOperator{\rank}{rank}
\DeclareMathOperator{\sgn}{sgn}

\DeclareMathOperator{\vspan}{span}
\DeclareMathOperator{\interior}{int}

\DeclarePairedDelimiterX{\norm}[1]{\lVert}{\rVert}{#1}
\DeclarePairedDelimiterX{\floor}[1]{\lfloor}{\rfloor}{#1}
\DeclarePairedDelimiterX{\qvar}[1]{\langle}{\rangle}{#1}
\DeclarePairedDelimiterX{\inn}[1]{\langle}{\rangle}{#1}

\begin{document}
	\maketitle

	\noindent	
	{\bf Abstract. } 
	Continuous-time Markov chains are frequently used as stochastic models for chemical reaction networks, especially in the growing field of systems biology. A fundamental problem for these Stochastic Chemical Reaction Networks (SCRNs)  is to understand the dependence of the stochastic behavior of these systems on the chemical reaction rate parameters. Towards solving this problem, in this paper we develop theoretical tools called comparison theorems that provide stochastic ordering results for SCRNs. These theorems give sufficient conditions for monotonic dependence on parameters in these network models, which allow us to obtain, under suitable conditions, information about transient and steady state behavior. These theorems exploit structural properties of SCRNs, beyond those of general continuous-time Markov chains. Furthermore, we derive two theorems to compare stationary distributions and mean first passage times for SCRNs with different parameter values, or with the same parameters and different initial conditions. These tools are developed for SCRNs taking values in a generic (finite or countably infinite) state space and can also be applied for non-mass-action kinetics models. When propensity functions are bounded, our method of proof gives an explicit method for coupling two comparable SCRNs, which can be used to simultaneously simulate their sample paths in a comparable manner. We illustrate our results with applications to models of enzymatic kinetics and epigenetic regulation by chromatin modifications.
	
	
		
\section{Introduction}
\label{sec:Introduction}

\subsection{Overview}

Stochastic Chemical Reaction Networks (SCRNs) are a class of continuous-time Markov chain models used to describe the stochastic dynamics of a chemical system undergoing a series of reactions which change the numbers of molecules of a finite set of species over time. These models provide a framework for the theoretical study of biochemical systems in areas such as intracellular viral kinetics (see Srivastava et al. \cite{Srivastava} and Haseltine \& Rawlings \cite{HaseltineRawlings}), enzymatic kinetics (see Kang et al. \cite{KangKhudaBukhshKoepplRempala} for example) and epigenetic regulation by chromatin modifications (see Bruno et al. \cite{BrunoDelVecchio} for a recently developed model of chromatin regulation).

One of the most interesting questions for biochemical system models is: ``What effect does changing reaction rate parameters have on system dynamics?" Indeed, different rate parameters for chemical processes can lead to different stochastic behaviors. One possible approach to evaluate the effect of parameter variations on system dynamics is through comparison theorems for stochastic processes. More precisely, this type of theorem provides inequalities between stochastic processes (see Muller \& Stoyan \cite{MullerStoyan} for a general reference on this topic).

In this paper, we employ uniformization and coupling methods (see Grassmann \cite{uniformization} and Keilson \cite{KeilsonRarityExponentiality}) to derive comparison theorems for SCRNs under verifiable sufficient conditions. These theoretical results enable us to develop two novel theorems yielding a direct comparison of mean first passage times and stationary distributions between SCRNs with different rate parameters or initial conditions. We apply these theorems to several examples to illustrate how they can be used to understand how key biological parameters affect stochastic behavior. While a major motivator for our work has been the study of SCRNs, we state our theorems in the context of continuous-time Markov chains, for which the state space is a subset of $\Z_+^d$ (the set of $d$-dimensional vectors with non-negative integer entries), and the set of all possible transition vectors is a finite set. This thereby allows for other applications that have similar characteristics to SCRNs. In addition, for the case of bounded transition intensities satisfying our conditions, we give an explicit concrete coupling of two comparable Markov chains, which can be used to simultaneously simulate them in such a way that their sample paths are monotonically related.

The paper is structured as follows: in Section \ref{Sec2:prelim} we introduce some background on stochastic chemical reaction networks needed for this article. We present the main results in Section \ref{sec:MainResults}, with proofs provided in Section \ref{sec:proofmainresults}. In Section \ref{sec:examples} we apply our theoretical tools to several examples, such as epigenetic regulation by chromatin modifications and enzymatic kinetics. Concluding remarks are presented in Section \ref{sec:conclusion}. The Supplementary information (SI) file contains some further details and extensions of the main results and examples in the paper.

\subsection{Related work}

Due to the growing field of systems biology, the mathematical study of chemical reaction networks has seen a wealth of activity lately. Concerning comparison results, considerable work has been conducted on monotonicity properties for deterministic models of chemical reaction networks, i.e., systems of ordinary differential equations describing the dynamics of species concentrations. For example, Angeli et al. \cite{2006Sontag} proposed a graphical method, based on the monotonicity properties of the reaction rates with respect to species concentrations, to determine global stability properties for the models. More recently, Gori et al. \cite{2019Nasti} introduced sufficient conditions to verify the existence of a monotonicity property for the concentrations of species for any positive time with respect to their initial concentrations. However, these works do not address how changing parameters affects the behavior of stochastic models.

To the best of our knowledge, no systematic study of stochastic ordering has been conducted for stochastic chemical reaction networks. On a more general level, theorems have been established for stochastic processes and have been specialized for particular classes such as for queueing systems and point processes (see Muller \& Stoyan \cite{MullerStoyan} for an introduction to the topic). For Markov chains, the work of Massey \cite{MasseyStochasticOrdering} is of special interest, since he establishes criteria for comparison of continuous-time Markov chains in terms of their infinitesimal generators. For relevant work prior to Massey, there is a nice summary in \cite{MasseyStochasticOrdering}. In particular, Kamae et al. \cite{KamaeKrengelOBrien} showed that for Markov processes, a comparison between transition probability functions, at all fixed times and for all partially ordered starting points, can be realized in a pathwise stochastic comparison between versions of the Markov processes. In relation to Massey's work, our results provide simplified conditions and extended results for stochastic comparisons, which exploit the structure of stochastic chemical reaction networks. Furthermore, unlike Massey, we do not require a uniform bound on the rates of leaving each state. In addition, under the latter assumption, we explicitly construct versions of the stochastic processes on the same probability space that have comparable sample paths. More detail on the relationship of our work to that of Massey is given in Remark \ref{rem:MasseyRemark}. In contrast to work on sensitivity analysis of distributions at a finite set of times and which considers only local changes in parameters (see for example Gunawan et al. \cite{2005Gunawan}, Gupta \& Khammash \cite{GuptaKhammash} and references therein), our work provides a sample path comparison between stochastic processes for global changes in their parameters.

\subsection{Notation and Terminology}
\label{sec:PreliminariesNotation}

Denote by $\Z_+ = \{0,1,2, \ldots \}$ the set of non-negative integers. For an integer $d \geq 1$ we denote by $\Z_+^d$ the set of $d$-dimensional vectors with entries in $\Z_+$. For any integer $d \geq 1$, let $\R^d$ denote the $d$-dimensional Euclidean space. We usually write $\R$ for $\R^1$. We denote by $\R^d_+$ the set of vectors $x \in \R^d$ such that $x_i \geq 0$ for every $1 \leq i \leq d$. For $x \in \R^d$, let $\norm{x}_{\infty} = \sup_{1 \leq i \leq d} |x_i|$ be the supremum norm. In this paper, the sum over the empty set is considered to be $0$.

A binary relation $\preccurlyeq$ on a set $\X$ will be called reflexive if $x \preccurlyeq x$ for every $x \in \X$, transitive if $x \preccurlyeq y$ and $y \preccurlyeq z$ implies $x \preccurlyeq z$ for every $x,y,z \in \X$ and antisymmetric if $x \preccurlyeq y$ and  $y \preccurlyeq x$ implies $x=y$ for every $x,y \in \X$. A preorder is a binary relation that is reflexive and transitive. A partial order is a preorder that is antisymmetric.

In this paper, a probability space $(\Omega,\F,\PP)$ will consist of a sample space $\Omega$, a $\sigma$-algebra of events $\F$ and a probability measure $\PP$ on $(\Omega, \F)$. We will say that two real-valued random variables $Y,Y'$ (defined on possibly different probability spaces) are equal in distribution, denoted as $Y' \overset{\text{dist}}{=} Y$, if their cumulative distribution functions agree. All stochastic processes considered in this paper will have right-continuous sample paths that also have finite left-limits. 
	
\section{Stochastic Chemical Reaction Networks (SCRNs)}
\label{Sec2:prelim}

In this section we provide necessary background on Stochastic Chemical Reaction Networks. The reader is referred to Anderson \& Kurtz \cite{AndersonKurtzBook} for an introduction to this subject. 

We assume there is a finite non-empty set $\Ss = \{\mathrm{S}_1,\ldots,\mathrm{S}_d\}$ of $d$ \textbf{species}, and a finite non-empty set $\Rs \subseteq \Z_+^d \times \Z_+^d$ that represents chemical \textbf{reactions}. We assume that $(w,w) \notin \Rs$ for every $w \in \Z^d_+$. The set $\Ss$ represents $d$ different molecular species in a system subject to reactions $\Rs$ which change the number of molecules of each species. For each $(v^{-},v^+) \in \Rs$, the $d$-dimensional vector $v^{-}$ (the \textbf{reactant vector}) counts how many molecules of each species are consumed in the reaction, while $v^{+}$ (the \textbf{product vector}) counts how many molecules of each species are produced. The reaction is usually written as
\begin{equation}
	\label{eq:ReactionNotation}
	\sum_{i=1}^d (v^{-})_{i}\mathrm{S}_i \longrightarrow  \sum_{i=1}^d (v^{+})_{i}\mathrm{S}_i.
\end{equation}
To avoid the use of unnecessary symbols, we will assume that for each $1 \leq i \leq d$, there exists a vector $w=(a_1, \ldots, a_d)^T \in \Z_+^d$ with $a_i >0$ such that $(w,v)$ or $(v,w)$ is in $\Rs$ for some $v \in \Z^d_+$, i.e., each species is either a reactant or a product in some reaction.

The net change in the quantity of molecules of each species due to a reaction $(v^{-},v^{+}) \in \Rs$ is described by $v^{+}-v^{-}$ and it is called the associated \textbf{reaction vector}. We denote the set of reaction vectors $\V := \{ v \in \Z^d \:|\: v = v^{+}- v^{-} \text{ for some } (v^{-},v^{+}) \in \Rs\}$, let $n := |\V|$ the size of $\V$ and enumerate the members of $\V$ as $\{v_1,\ldots,v_n\}$. Note that $\V$ does not contain the zero vector because $\Rs$ has no elements of the form $(w,w)$. Different reactions might have the same reaction vector. For each $v_j \in \V$ we consider the set $\Rs_{v_j} := \{(v^{-},v^{+}) \in \Rs \:|\: v_j =v^{+}-v^{-} \}$. The reaction vectors generate the \textbf{stoichiometric subspace} $\L:=  \vspan(\V)$. For $z \in \R^d$, we call $z + \L$ a \textbf{stoichiometric compatibility class}.

Given $(\Ss,\Rs)$ we will consider an associated continuous-time Markov chain $X=(X_1,\ldots,X_d)$, with a state space $\X$ contained in $\Z^d_+$, which tracks the number of molecules of each species over time. Roughly speaking, the dynamics of $X$ will be given by the following: given a current state  $x=(x_1,\ldots,x_d) \in \X \subseteq \Z_+^{d}$, for each reaction $(v^{-},v^{+}) \in \Rs$, there is a clock which will ring at an exponentially distributed time (with rate $\Lambda_{(v^{-},v^{+})}(x)$). The clocks for distinct reactions are independent of one another. If the clock corresponding to $(v^{-},v^{+})\in \Rs$ rings first, the system moves from $x$ to $x+v^{+}- v^{-}$ at that time, and then the process repeats. We now define the Markov chain in more detail.

Consider a set of species $\Ss$ and of reactions $\Rs$, a set $\X \subseteq \Z^d_+$ and a collection of functions $\{\prp_{(v^{-},v^{+})}:\X \longrightarrow \R_+\}_{(v^{-},v^{+}) \in \Rs}$ such that for each $x \in \X$ and $(v^{-},v^{+}) \in \Rs$, if $x+v^{+}-v^{-} \notin \X$, then $\Lambda_{(v^{-},v^{+})}(x)=0$. Now, for $1 \leq j \leq n$, $v_j \in \V$, define
\begin{equation}
	\rate_j(x) := \sum_{(v^{-},v^{+}) \in \Rs_{v_j}} \prp_{(v^{-},v^{+})}(x).
\end{equation}
Note that for each $x \in \X$ and $1 \leq j \leq n$, if $x +v_j \notin \X$, then $\rate_j(x) = 0$. A \textbf{stochastic chemical reaction network (SCRN)} is a Markov chain $X$ with state space $\X$ and infinitesimal generator\footnote{Note that $Q$ is sometimes called an infinitesimal transition matrix although it may have countably many ``rows" and ``columns". The entries $Q_{x,y}$ for $x \neq y$ are the infinitesimal transition rates of going from $x$ to $y$: $\PP[X(t+h) = y | X(t)=x] = Q_{x,y}h + o(h)$ as $h \to 0$.} $Q$ given for $x,y \in \X$ by
\begin{equation}
	\label{eq:TransitionMatrixQ}
	Q_{x,y} = \begin{cases}
		\rate_j(x) & \text{ if } y-x = v_j \text{ for some } 1 \leq j \leq n, \\
		- \sum_{j=1}^n\rate_j(x) & \text{ if } y = x, \\
		0 & \text{ otherwise.}
	\end{cases}
\end{equation}
The functions $\{\prp_{(v^{-},v^{+})}:\X \longrightarrow \R_+\}_{(v^{-},v^{+}) \in \Rs}$ are called \textbf{propensity} or \textbf{intensity} functions. A common form for the propensity functions is the following associated with \textbf{mass action kinetics}:
\begin{equation}
	\label{eq:PropensityFunctions}
	\prp_{(v^{-},v^{+})}(x) = \kappa_{(v^{-},v^{+})}\prod_{i=1}^{d}(x_i)_{(v^{-})_i},
\end{equation}
where $\{\kappa_{(v^{-},v^{+})}\}_{(v^{-},v^{+}) \in \Rs}$ are positive  constants and for $m,\l \in \Z_+$, the quantity $(m)_\l$ is the falling factorial, i.e., $(m)_0 := 1$ and $(m)_\l := m(m-1)\ldots(m-\l+1)$.

\begin{remark}
	Our definition of SCRN allows for some model flexibility. Notice that the propensity functions are not necessarily defined on the whole lattice $\Z^d_+$ and they are not necessarily of the form \eqref{eq:PropensityFunctions}. Indeed, in some of our applications, mass-conservation laws restrict the possible values that $X$ may take (see Example \ref{ex:HistoneModification}). In addition, there may be other types of kinetics, such as those described by Hill functions (see Example \ref{ex:HistoneModificationAndProtein}).
\end{remark}

A convenient way to represent such a Markov chain is given in Theorem 6.4.1 of Ethier \& Kurtz \cite{EthierKurtz}. For this, consider a probability space $(\Omega, \F, \PP)$ equipped with independent unit rate Poisson processes $N_1, \ldots,N_n$. There is a version of $X$ defined on $(\Omega,\F,\PP)$ such that
\begin{equation}
	\label{eq:PossionProcessRepresentation}
	X(t) = X(0) + \sum_{j=1}^{n} v_jN_j\left(\int_{0}^{t}\rate_j\left(X(s)\right)ds \right),
\end{equation}
for every $0 \leq t < \tau$, where $\tau$ is the explosion time for $X$ (which may be $+\infty$). From \eqref{eq:PossionProcessRepresentation}, it is easy to see that for a SCRN $X$ with initial state $z \in \X$, $X(t)$ will stay in the stoichiometric compatibility class $z +\L$ intersected with $\Z^d_+$ for all time $0 \leq t < \tau$, with probability one. For this reason, sometimes it will be convenient to choose $\X = (z+\L) \cap \Z^d_+$, for a fixed $z \in \Z^d_+$.

While our work was initially motivated by questions for SCRNs, we will first develop our results in a more general context of continuous-time Markov chains, for which the state space is contained in $\Z^d_+$ and the set of all possible transition vectors is a finite set, and then illustrate them for SCRNs.

\section{Main Results}
\label{sec:MainResults}

The general stochastic ordering results provided in this paper are relative to a preorder relation on a state space $\X \subseteq \Z^d_+ \subseteq \R^d$. We will define the preorder on all of $\R^d$ and then restrict it to various subsets. We introduce this notation and related notation in Section \ref{sec:PreordersRd}. In Section \ref{sec:CouplingOfCRNs} we present the main results of this article, and in Section \ref{sec:MonotonicPropMFPTandSD} we discuss relevant consequences for the comparison of (mean) first passage times and stationary distributions.

\subsection{Preorders in $\R^d$}
\label{sec:PreordersRd}

Let $m,d \geq 1$ be integers. Denote by $\leq$ the usual componentwise partial order on $\R^d$, i.e., for $x,y \in \R^d$, $x \leq y$ whenever $x_i \leq y_i$ for every $1 \leq i \leq d$. Additionally, we write $x < y$ whenever $x_i < y_i$ for every $1 \leq i \leq d$. For the rest of the paper, we consider a matrix $A \in \R^{m \times d}$, where no row of $A$ is identically zero.

\begin{definition}
	\label{def:preorderA}
	For $x,y \in \R^d$, we say that $x \preccurlyeq_A y$ whenever $A(y-x) \geq 0$.
\end{definition}

For the matrix $A$, consider the convex cone $K_A := \{ x \in \R^d \:|\: Ax \geq 0\}$. Note that $x \preccurlyeq_A y$ holds if and only if $y-x \in K_A$. Moreover, the relation $\preccurlyeq_A$ is reflexive and transitive, and therefore a preorder on $\R^d$. Also, for this relation,
\begin{equation}
	\label{eq:translations_precurly}
	\text{if } x \preccurlyeq_A y, \text{ then } x +z \preccurlyeq_A y +z \text{ for any } z \in \R^d.     
\end{equation} 
For any $x \in \R^d$ consider the set
\begin{equation*}
	K_A +x =\{ y \in \R^d \:|\: A(y-x) \geq 0 \} = \{ y \in \R^d \:|\: x \preccurlyeq_A y \}.
\end{equation*}

In the coming sections, we will consider the notions of increasing and decreasing sets with respect to $\preccurlyeq_A$ in  a given subset of $\Z_+^d$. More concretely, consider a non-empty set $\X \subseteq \Z_+^d$. We will say that a set $\Gamma \subseteq \X$ is \textbf{increasing} in $\X$ with respect to $\preccurlyeq_A$ if for every $x \in \Gamma$ and $y \in \X$, $x \preccurlyeq_A y$ implies that $y \in \Gamma$. We observe that, for $x \in \X$, the set
\begin{equation}
	(K_A +x) \cap \X = \{ y \in \X \:|\: x \preccurlyeq_A y \}
\end{equation}
is increasing in $\X$ by the transitivity property of $\preccurlyeq_A$. On the other hand, we will say that a set $\Gamma \subseteq \X$ is \textbf{decreasing} in $\X$ with respect to $\preccurlyeq_A$ if for every $x \in \Gamma$ and $y \in \X$, $y \preccurlyeq_A x$ implies that $y \in \Gamma$. We will say that a point $x$ is \textbf{maximal} (resp. \textbf{minimal}) in $\X$ if for every $y \in \X$, $x \preccurlyeq_A y$ (resp. $y \preccurlyeq_A x$) implies that $x=y$. In this case, the set $\Gamma = \{x\}$ would be increasing (resp. decreasing) in $\X$.

\begin{remark}
	If $\rank(A)=d$, then the relation $\preccurlyeq_A$ will be antisymmetric and therefore a partial order on $\R^d$. 
	Indeed, if $\rank(A)=d$, then $A(y-x)=0$ implies that $x=y$. In addition, $\preccurlyeq_A$ will then be a partial order when restricted to $\X\subset \Z_+^d$. Throughout this article, we will not assume that $\rank(A)=d$ and therefore, the relation $\preccurlyeq_A$ might not be a partial order on $\X$ (see Examples \ref{exEK}, \ref{exEK2}, and \ref{ex:Braess}). 
\end{remark}

\subsection{Stochastic Comparison Theorems}
\label{sec:CouplingOfCRNs}

The fundamental objects in the following results are a non-empty set $\X \subseteq \Z_+^d$ and a pair of continuous-time Markov chains $X$ and $\breve{X}$ with the same state space $\X$ and where it is assumed that the set of all possible transition vectors for $X$ or $\Breve{X}$ is a finite set. We denote the size of this set by $n$. A primary example of this setup is two stochastic chemical reaction networks as described in Section \ref{Sec2:prelim} with different propensity functions. We will now formally introduce the notation for stating our results.

\medskip
Consider a non-empty set $\X \subseteq \Z_+^d$, an integer $n \geq 1$ and a collection of distinct vectors $v_1,\ldots,v_n$ in $\Z^d \setminus \{0\}$, where $0$ is the origin in $\Z^d$. Consider two collections of functions $\rate=(\rate_1, \dots,\rate_n)$ and $\Breve{\rate}= (\Breve{\rate}_1, \dots,\Breve{\rate}_n)$ defined on $\X$ and taking values in $\R_+$, such that for every $1 \leq j \leq n$ and $x \in \X$:
\begin{equation}
\label{eq:LambdaNotOutofSpace}
\text{if } x +v_j \notin \X, \text{ then } \rate_j(x) = \Breve{\rate}_j(x)= 0.   
\end{equation}
Consider a continuous-time Markov chain $X$ on the state space $\X$ with infinitesimal generator $Q=(Q_{x,y})_{x,y \in \X}$ defined for $x,y \in \X$ by
\begin{equation}
\label{eq:GeneratorAStochasticOrder}
Q_{x,y} := \begin{cases}
	\rate_j(x) & \text{ if } y-x = v_j \text{ for some } 1 \leq j \leq n, \\
	- \sum_{j=1}^n\rate_j(x) & \text{ if } x = y, \\
	0 & \text{ otherwise. }
\end{cases}
\end{equation}
Consider the analogous continuous-time Markov chain $\Breve{X}$ with infinitesimal generator $\breve{Q}$ as in \eqref{eq:GeneratorAStochasticOrder} but with functions $\Breve{\rate}_1, \dots,\Breve{\rate}_n$ instead of $\rate_1, \dots,\rate_n$. We call $X$ and $\Breve{X}$ the continuous-time Markov chains associated with $\rate$ and $\Breve{\rate}$ respectively. We will assume that $X$ and $\Breve{X}$ do not explode in finite time. The following is our main result.

	\begin{figure}[t]
		\centering
		\includegraphics[scale=0.5]{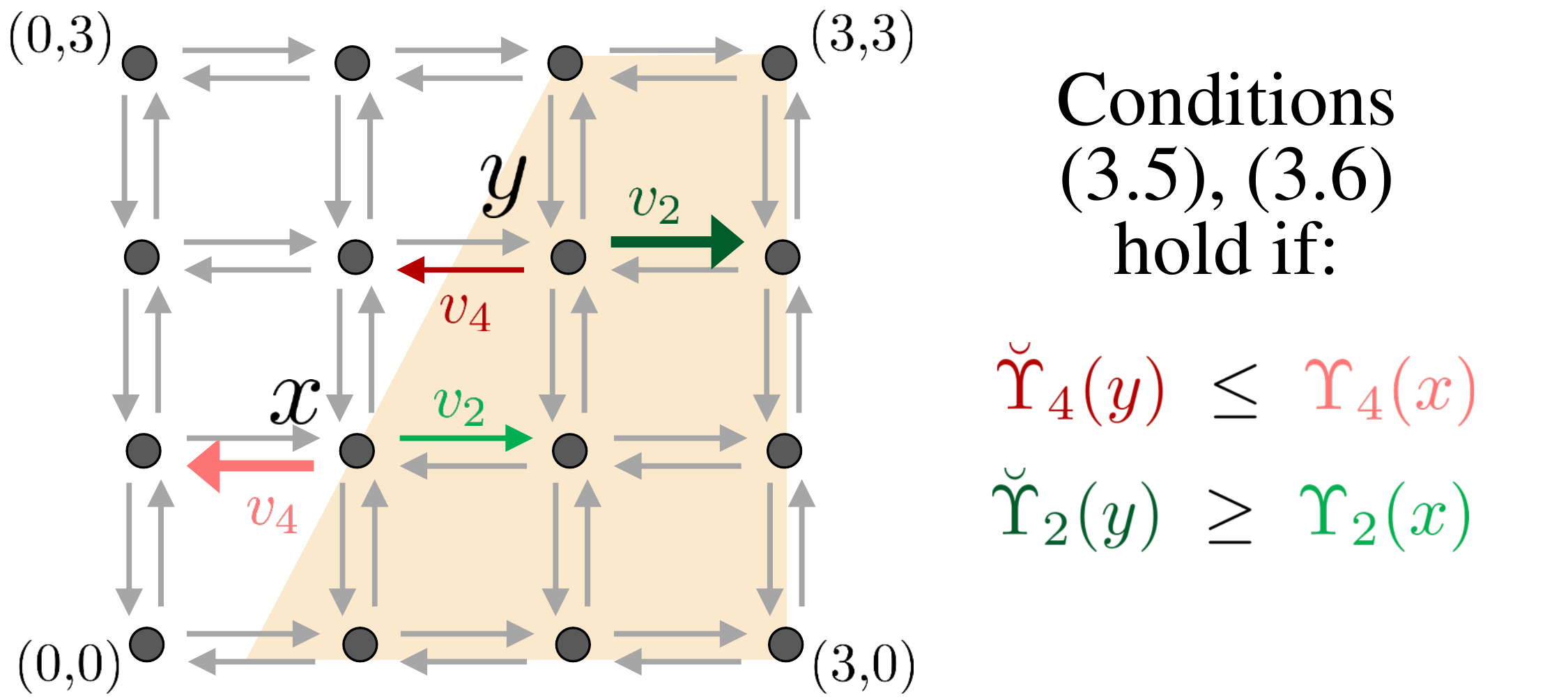}
		\caption{\small { \bf Pictorial representation of conditions (\ref{eq:CouplingCondition1}), (\ref{eq:CouplingCondition2}) for a certain $\intervalx$ in a two-dimensional lattice.} Here, $\X=\{0,1, 2, 3\} \times \{0,1, 2, 3\}$, $n=4$, $v_1=(0,1)^T$, $v_2=(1,0)^T$, $v_3=(0, -1)^T$, $v_4=(-1,0)^T$, where $T$ denotes transpose, $A=[2\; -1]$, and $\intervalx= \{ w \in \mathcal{X} \:|\: [2\; -1](w - x) \geq 0 \}$. In the graph, $\intervalx$ consists of the states (black dots) that lie in the light orange region and the arrows represent possible transitions along $v_1,v_2,v_3,v_4$ between states. For the exhibited states $x,y \in \X$ with $x \preccurlyeq_A y$, the light green (dark green) and light red (dark red) arrows represent the transitions with rates $\rate_2(x)$ ($\Breve{\rate}_2 (y)$) and $\rate_4(x)$ ($\Breve{\rate}_4 (y)$) for the Markov chain $X$ ($\Breve{X}$). Higher transitions rates are associated with thicker arrows. To check the conditions (\ref{eq:CouplingCondition1}) and (\ref{eq:CouplingCondition2}), since $y+v_4 \notin K_A + x$ and $y \notin K_A +x+v_2$, we need to check that $\Breve{\rate}_4(y) \leq \rate_4(x)$ and $\Breve{\rate}_2(y) \geq \rate_2(x)$.}
		\label{fig:conditions}
	\end{figure}

	\begin{theorem}
		\label{thm:MainResult}
		Consider a non-empty set $\X \subseteq \Z_+^d$, a collection of distinct vectors $v_1,\ldots,v_n$ in $\Z^d \setminus \{0\}$ and two collections of non-negative functions on $\X$, $\rate=(\rate_1, \dots,\rate_n)$ and  $\Breve{\rate}= (\Breve{\rate}_1, \dots,\Breve{\rate}_n)$, such that \eqref{eq:LambdaNotOutofSpace} holds and the associated continuous-time Markov chains do not explode in finite time. Consider a matrix $A \in \R^{m \times d}$ with non-zero rows and suppose that for every $x,y \in \X$ such that $x \preccurlyeq_A y$ the following hold:
		\begin{equation}
			\label{eq:CouplingCondition1}
			\Breve{\rate}_j(y) \leq \rate_j(x), \quad \text{for each } 1 \leq j \leq n \text{ such that } y +v_j \in \X \setminus (K_A +x),
		\end{equation}
		and
		\begin{equation}
			\label{eq:CouplingCondition2}
			\Breve{\rate}_j(y)\geq \rate_j(x), \quad \text{for each } 1 \leq j \leq n \text{ such that } x +v_j \in \X \text{ and } y \notin K_A + x + v_j.
		\end{equation}
		Then, for each pair $\initialx,\initialxbreve \in \X$ such that $\initialx \preccurlyeq_A \initialxbreve$, there exists a probability space $(\Omega,\F,\PP)$ with two continuous-time Markov chains $X = \{X(t), \: t \geq 0\}$ and $\Breve{X}=\{\Breve{X}(t), \: t \geq 0\}$ defined there, each having state space $\X \subseteq \Z^d_+$, with infinitesimal generators $Q$ and $\Breve{Q}$, associated with $\rate$ and $\Breve{\rate}$, respectively, with initial conditions $X(0)=\initialx$ and $\Breve{X}(0)=\initialxbreve$ and such that:
		\begin{equation}
			\label{eq:CouplingCondition}
			\PP\left[X(t) \preccurlyeq_A \Breve{X}(t) \text{ for every } t \geq 0 \right]=1.   
		\end{equation}
	\end{theorem}
	
	An example of checking conditions \eqref{eq:CouplingCondition1} and \eqref{eq:CouplingCondition2} is given in  Figure \ref{fig:conditions}. The proof of Theorem \ref{thm:MainResult} is given in Section \ref{PROOFtheoremMonAppendix_I}. The main idea in the construction of the processes $X$ and $\Breve{X}$ is \textit{uniformization} (see Chapter 2 in Keilson \cite{KeilsonRarityExponentiality}) together with a suitable coupling. Our proof uses a single Poisson process together with a sequence of i.i.d. uniform random variables to determine potential jumps for the two continuous-time Markov chains, where for $X$ and $\Breve{X}$, potential jumps in the same direction $v_j$ are coupled together, and their probabilities of acceptance are given by normalized versions of their infinitesimal transition rates $\rate_j$ and $\Breve{\rate}_j$. Uniformization can be done provided the diagonal terms of the infinitesimal generators are uniformly bounded in size. In the proof of Theorem \ref{thm:MainResult}, we initially make this assumption on $Q$ and $\Breve{Q}$ in order to construct $X$ and $\Breve{X}$. We then generalize the result to Markov chains that do not explode in finite time by using a truncation and limiting procedure. The construction mentioned here, for the case where the diagonal terms of the infinitesimal generators are uniformly bounded in size, besides playing a key role in our proofs, is also the basis for an algorithm described in SI - Section S.4, which provides a way to simultaneously simulate the processes $X$ and $\Breve{X}$ in a comparable manner.
		
	\begin{remark}
		\label{rem:MasseyRemark}
		In Theorem 5.3 of Massey \cite{MasseyStochasticOrdering}, the author provides a necessary and sufficient condition for stochastic comparison of continuous-time Markov chains at each fixed time for all partially ordered initial conditions. By the work of Kamae et al. \cite{KamaeKrengelOBrien}, the conditions in \cite{MasseyStochasticOrdering} imply the existence of a coupling of continuous-time Markov chains so that a relation such as \eqref{eq:CouplingCondition} holds. Massey's condition requires that $\sum_{w \in \Gamma} Q_{xw} \leq \sum_{w \in \Gamma}\Breve{Q}_{yw}$ for every $x \preccurlyeq_A y$ and every set $\Gamma \subseteq \X$ that is increasing in $\X$ with respect to $\preccurlyeq_A$ and such that either $x \in \Gamma$ or $y \notin \Gamma$. These inequalities can often be difficult to check since first, they involve computing sums of terms in the infinitesimal generators and second, the form of all increasing sets can be hard to determine. In Theorem \ref{thm:MainResult} we overcome these obstacles by providing simplified sufficient conditions that involve only pointwise comparison of entries in the infinitesimal generators associated to each of the transition vectors $v_j$. Besides this practical value, in our context, our results go beyond the work of Massey \cite{MasseyStochasticOrdering}, since he assumes that $\preccurlyeq_A$ is a partial order (we only assume preorder) and he assumes that the diagonal entries of the infinitesimal generators are bounded (we generalize to non-exploding Markov chains). Our proof has a commonality with the work of Massey in the sense that we also use uniformization. It is different in the sense that, when infinitesimal transition rates are bounded, we construct an explicit coupling for all time, exploiting the simplified nature of our conditions, while Massey does not provide an explicit coupling. Instead, he proves existence of a stochastic comparison for each fixed time, using a semigroup approach.
	\end{remark}
	
	Conditions \eqref{eq:CouplingCondition1} and \eqref{eq:CouplingCondition2} may be simplified if we consider a particular relation between the matrix $A$ and the vectors $v_1,\ldots,v_n$ in which $A \in \Z^{m \times d}$ and $Av_j$ has entries taking values only in $\{-1,0,1\}$ for every $1\leq j \leq n$. More concretely, let us consider a class of continuous-time Markov chains such that, for a given matrix $A$ with non-zero rows, if the Markov chain starts within the set $K_A + x$, then to go outside of it, the process will necessarily hit its boundary. In this case, we can derive a theorem whose conditions must be checked only on the boundary of $K_A + x$ because the only transitions that can lead the Markov chain outside or inside the set $K_A + x$ are ones starting on the boundary of $K_A + x$. Before stating the theorem, let us introduce the sets $\partial_i(K_A+x) := \{ y \in K_A +x \:|\: \inn{A_{i\bullet},y} = \inn{A_{i\bullet},x} \}$ \footnote{Here, for convenience of notation, let $A_{i\bullet}$ denote the row vector corresponding to the $i$-th row of $A$, for $1 \leq i \leq m$. In this article we will adopt the convention of considering the inner product $\inn{\cdot,\cdot}$ as a function of a row vector in its first entry and as a function of a column vector in the second entry. In particular, $\inn{A_{i\bullet},x} = \sum_{k=1}^d A_{ik}x_k$.} for each $1 \leq i \leq m$. We can then characterize\footnote{The fact that $A$ does not contain zero rows allows for equation \eqref{eq:BoundaryKAx} to hold. In fact, let $A \in \R^{m \times d}$ be a matrix that is not identically zero and let $\bar{A}$ be the matrix obtained from $A$ by erasing any rows that contain all zeros. Then, for $x,y \in \R^d$, $A(y-x) \geq 0$ if and only if $\bar{A}(y-x) \geq 0$, and so $K_A+x = K_{\bar{A}}+x$ and $x \preccurlyeq_A y$ if and only if $x \preccurlyeq_{\bar{A}} y$. However, if $A$ contains a row $A_{i\bullet}$ such that $A_{i\bullet}=0$, then $\partial_i(K_A +x) = K_A +x$ and if $K_A +x$ has nonempty interior, then equation \eqref{eq:BoundaryKAx} will not hold. Consequently, we have made the assumption that $A$ has no zero rows.} the boundary of $K_A + x$ as follows:
	\begin{equation}
		\label{eq:BoundaryKAx}
		\partial(K_A +x) = \bigcup_{i=1}^{m}  \partial_i(K_A+x).  
	\end{equation}

	\begin{theorem}
		\label{thm:InnerProductTheorem}
		Consider a non-empty set $\X \subseteq \Z_+^d$, a collection of distinct vectors $v_1,\ldots,v_n$ in $\Z^d \setminus \{0\}$ and two collections of non-negative functions on $\X$, $\rate=(\rate_1, \dots,\rate_n)$ and  $\Breve{\rate}= (\Breve{\rate}_1, \dots,\Breve{\rate}_n)$ such that \eqref{eq:LambdaNotOutofSpace} holds and the associated continuous-time Markov chains do not explode in finite time. Consider a matrix $A \in \Z^{m \times d}$ with non-zero rows and suppose that both of the following conditions hold:
		\begin{enumerate}
			\item[(i)]
			For each $1\leq j \leq n$, the vector $Av_j$ has entries in $\{-1,0,1\}$ only.
			\item[(ii)]
			For each $x \in \X$, $1 \leq i \leq m$ and $y \in \partial_i(K_A+x) \cap \X$ we have that
			\begin{equation}
				\label{eq:CouplingConditions_InnerProduct_I}
				\Breve{\rate}_j(y) \leq \rate_j(x), \quad \text{for each } 1 \leq j \leq n \text{ such that } \inn{A_{i\bullet},v_j} < 0,
			\end{equation}
			and  
			\begin{equation}
				\label{eq:CouplingConditions_InnerProduct_II}
				\Breve{\rate}_j(y) \geq \rate_j(x), \quad \text{for each } 1 \leq j \leq n \text{ such that } \inn{A_{i\bullet},v_j} > 0.
			\end{equation}
		\end{enumerate}
		Then, for each pair $\initialx,\initialxbreve \in \X$ such that $\initialx \preccurlyeq_A \initialxbreve$, there exists a probability space $(\Omega,\F,\PP)$ with two continuous-time Markov chains $X = \{X(t), \: t \geq 0\}$ and $\Breve{X}=\{\Breve{X}(t), \: t \geq 0\}$ defined there, each having state space $\X \subseteq \Z^d_+$, with infinitesimal generators given by $Q$ and $\Breve{Q}$, associated with $\rate$ and $\Breve{\rate}$ respectively, with initial conditions $X(0)=\initialx$ and $\Breve{X}(0)=\initialxbreve$, and such that:
		\begin{equation}
			\label{eq:DominanceProperty_InnerProduct}
			\PP\left[X(t) \preccurlyeq_A \Breve{X}(t) \text{ for every } t \geq 0 \right]=1.   
		\end{equation}
	\end{theorem}    
	
	The proof of this theorem is given in Section \ref{sec:ProofInnerProductTheorem} and involves checking that \eqref{eq:CouplingCondition1} and \eqref{eq:CouplingCondition2} of Theorem \ref{thm:MainResult} hold, using conditions \it (i) \rm and \it (ii) \rm of Theorem \ref{thm:InnerProductTheorem}.
	
	\begin{remark}
		\label{rem:EqualRates}
		In the context of Theorem \ref{thm:InnerProductTheorem}, it is possible that for $x \in \X$, and $y \in \partial_{i_1}(K_A+x) \cap \partial_{i_2}(K_A+x) \cap \X$ with $i_1 \neq i_2$, it happens that $\inn{A_{i_1\bullet},v_j} < 0$ and $\inn{A_{i_2\bullet},v_j} > 0$ for some $1 \leq j \leq n$. For condition $(ii)$ to hold, we must then have $\Breve{\rate}_j(y) = \rate_j(x)$.
	\end{remark}
	
	When there are multiple vectors $v_j$ with a common value for $Av_j$, the pointwise comparison in $j$, for $1 \leq j \leq n$, in conditions \eqref{eq:CouplingConditions_InnerProduct_I} and \eqref{eq:CouplingConditions_InnerProduct_II} in Theorem \ref{thm:InnerProductTheorem}, can be weakened. To this end, let us introduce the set of distinct vectors $\{\eta^1,\dots,\eta^s\}$ formed by $Av_j$, for $1 \leq j \leq n$, where $s$ denotes the cardinality of this set. Consider the subsets of indices
	\begin{equation}
		\label{eq:sets_of_dot_product+-13}
		G^{k} :=  \{ j  \:|\: 1\leq j \leq n \hbox{ and } Av_j =\eta^k \}, \quad \text{for } 1 \leq k \leq s.
	\end{equation}  
		Then we have the following theorem.
		\begin{theorem}
			\label{thm:SumTheorem2}
			Consider a non-empty set $\X \subseteq \Z_+^d$, a collection of distinct vectors $v_1,\ldots,v_n$ in $\Z^d \setminus \{0\}$ and two collections of non-negative functions on $\X$, $\rate=(\rate_1, \dots,\rate_n)$ and  $\Breve{\rate}= (\Breve{\rate}_1, \dots,\Breve{\rate}_n)$ such that \eqref{eq:LambdaNotOutofSpace} holds and the associated continuous-time Markov chains do not explode in finite time. Consider a matrix $A \in \Z^{m \times d}$ with non-zero rows and suppose that both of the following conditions hold:
			\begin{enumerate}
				\item[(i)]
				For each $1\leq j \leq n$, the vector $Av_j$ has entries in $\{-1,0,1\}$ only.
				\item[(ii)]
				For each $x \in \X$, $1 \leq i \leq m$ and $y \in \partial_i(K_A+x) \cap \X$ we have that
				\begin{equation}
					\label{eq:CouplingConditionBIS32}
					\sum_{j \in G^{k}} \Breve{\rate}_j(y) \leq \sum_{j \in G^{k}} \rate_j(x), \quad \text{for each } k \text{ such that } \eta^k_i <0,
				\end{equation}
				and
				\begin{equation}
					\label{eq:CouplingConditionBIS31}
					\sum_{j \in G^{k} } \Breve{\rate}_j(y) \geq \sum_{j \in G^{k}} \rate_j(x), \quad \text{for each } k \text{ such that } \eta^k_i >0.
				\end{equation}
			\end{enumerate}
			Then, for each pair $\initialx,\initialxbreve \in \X$ such that $\initialx \preccurlyeq_A \initialxbreve$, there exists a probability space $(\Omega,\F,\PP)$ with two continuous-time Markov chains $X = \{X(t), \: t \geq 0\}$ and $\Breve{X}=\{\Breve{X}(t), \: t \geq 0\}$ defined there, each having state space $\X \subseteq \Z^d_+$, with infinitesimal generators $Q$ and $\Breve{Q}$, associated with $\rate$ and $\Breve{\rate}$ respectively, with initial conditions $X(0)=\initialx$ and $\Breve{X}(0)=\initialxbreve$ and such that:
			\begin{equation}
				\label{eq:CouplingConditionBIS}
				\PP\left[X(t) \preccurlyeq_A \Breve{X}(t) \text{ for every } t \geq 0 \right]=1.   
			\end{equation}
		\end{theorem}

		The proof of this theorem is given in Section \ref{sec:ProofSumTheorem2}.
		
		\begin{remark}
			\label{remark:incond}
			If $\Upsilon = \Breve{\Upsilon}$, Theorems \ref{thm:MainResult}, \ref{thm:InnerProductTheorem} and  \ref{thm:SumTheorem2} give sufficient conditions for monotonic dependence of the stochastic dynamic behavior on the initial condition. In the sense of Massey \cite{MasseyStochasticOrdering}, this notion corresponds to constructing a \textit{strongly monotone Markov chain}.
		\end{remark}
		
		\begin{remark}
			For \textit{deterministic} dynamical systems, there is a considerable literature giving monotonicity conditions with respect to initial conditions (see e.g., Hirsch \& Smith \cite{HirschSmith2004}). Furthermore, Angeli \& Sontag \cite{2003Sontag} extended the concept of monotone systems to systems having external inputs (i.e., $\dot x = f(x,u)$, with $x$ representing the state and $u$ representing the input). More precisely, they developed tools to prove monotonic dependence of the deterministic dynamic behavior on the initial condition and external input, provided that certain sign conditions on the first partial derivatives of the function $f(x,u)$ are satisfied on the entire state and input space. These theoretical tools can be used also to study how changing a system parameter affects the deterministic behavior of the system, by viewing $u$ as the system parameter of interest.
		\end{remark}

		\begin{remark}
			Checking the conditions in Theorems \ref{thm:InnerProductTheorem} and \ref{thm:SumTheorem2} (if they hold) is less cumbersome than checking the conditions in Theorem \ref{thm:MainResult}. In fact, compared to Theorem \ref{thm:MainResult}, for Theorems \ref{thm:InnerProductTheorem} and \ref{thm:SumTheorem2}, the conditions must be checked only on the boundaries of $K_A+x$, given that condition $(i)$ there is assumed to hold. Furthermore, Theorem \ref{thm:SumTheorem2} has less restrictive conditions (i.e., comparing sums of infinitesimal rates associated with transitions inwards or outwards with respect to the hyperplanes $\{ z \in \R^d \:|\:\inn{A_{i\bullet},z} =\inn{A_{i\bullet},x} = \inn{A_{i\bullet},y} \}$, $1 \leq i \leq m$, instead of comparing transition rates one-by-one for $1 \leq j \leq n$).
		\end{remark}

		\subsection{Monotonicity properties for (Mean) First Passage Times and Stationary Distributions}
		\label{sec:MonotonicPropMFPTandSD}
		
		The first consequence of our main results is for first passage times and it is related to stochastic orderings of real-valued random variables. Let $Y$ and $Z$ be real-valued random variables with cumulative distribution functions $F_Y$ and $F_Z$ respectively. We say that $Y$ is smaller than $Z$ in the \textbf{usual stochastic order}, written $Y \preccurlyeq_{st} Z$ if $F_Y(t) \geq F_Z(t)$ for every $t \in \R$. The relation $Y \preccurlyeq_{st} Z$ is equivalent to the existence of a probability space $(\Omega,\F,\PP)$ with random variables $Y' \overset{\text{dist}}{=} Y$ and $Z' \overset{\text{dist}}{=} Z$ defined there such that $\PP(Y' \leq Z') =1$. Furthermore, it is equivalent to the condition: $\int_{-\infty}^{\infty} f(x) dF_Y(x) \leq\int_{-\infty}^{\infty} f(x) dF_Z(x)$ for every bounded increasing function $f:\R \longrightarrow \R$. The reader may consult Chapter 1 in Muller \& Stoyan \cite{MullerStoyan} for the corresponding proofs and further properties of this notion.

		\begin{theorem}
			\label{thm:comparison_MFPT}
			Consider a non-empty set $\X \subseteq \Z_+^d$, a collection of distinct vectors $v_1,\ldots,v_n$ in $\Z^d \setminus \{0\}$ and two collections of non-negative functions on $\X$, $\rate=(\rate_1, \dots,\rate_n)$ and  $\Breve{\rate}= (\Breve{\rate}_1, \dots,\Breve{\rate}_n)$, such that \eqref{eq:LambdaNotOutofSpace} holds and the associated continuous-time Markov chains do not explode in finite time. Consider a matrix $A \in \R^{m \times d}$ with non-zero rows and suppose that at least one of the following holds:
			\begin{enumerate}
				\item[(i)]
				For every $x,y \in \X$ such that $x \preccurlyeq_A y$, conditions \eqref{eq:CouplingCondition1} and \eqref{eq:CouplingCondition2} are satisfied.
				\item[(ii)]
				The matrix $A$ has integer-valued entries and conditions $(i)$ and $(ii)$ in Theorem \ref{thm:InnerProductTheorem} are satisfied.
				\item[(iii)]
				The matrix $A$ has integer-valued entries and conditions $(i)$ and $(ii)$ in Theorem \ref{thm:SumTheorem2} are satisfied.
			\end{enumerate}
			Let $\initialx,\initialxbreve \in \X$ be such that $\initialx \preccurlyeq_A \initialxbreve$ and let $X = \{X(t), \: t \geq 0\}$ and $\Breve{X}=\{\Breve{X}(t), \: t \geq 0\}$ be two continuous-time Markov chains (possibly defined on different probability spaces), each having state space $\X \subseteq \Z^d_+$, with infinitesimal generators $Q$ and $\Breve{Q}$, associated with $\rate$ and $\Breve{\rate}$ respectively, and with initial conditions $X(0)=\initialx$ and $\Breve{X}(0)=\initialxbreve$. For a non-empty set $\Gamma \subseteq \X$, consider $T_{\Gamma} := \inf\{ t \geq 0 \:|\: X(t) \in \Gamma \}$ and $\Breve{T}_{\Gamma} := \inf\{ t \geq 0 \:|\: \Breve{X}(t) \in \Gamma \}$. If $\Gamma$ is increasing in $\X$ with respect to the relation $\preccurlyeq_A$, then
			\begin{equation}
				\label{eq:StochasticOrderingFPT_increasing}
				\Breve{T}_{\Gamma} \preccurlyeq_{st} T_\Gamma,
			\end{equation}
			and the mean first passage time of $\Breve{X}$ from $\initialxbreve$ to $\Gamma$ is dominated by the mean first passage time of $X$ from $\initialx$ to $\Gamma$. If $\Gamma$ is decreasing in $\X$ with respect to the relation $\preccurlyeq_A$, then
			\begin{equation}
				\label{eq:StochasticOrderingFPT_decreasing}
				T_\Gamma \preccurlyeq_{st} \Breve{T}_{\Gamma},
			\end{equation}
			and the mean first passage time of $X$ from $\initialx$ to $\Gamma$ is dominated by the mean first passage time of $\Breve{X}$ from $\initialxbreve$ to $\Gamma$.
		\end{theorem}
		
		\begin{proof}
			By Theorem \ref{thm:MainResult}, \ref{thm:InnerProductTheorem} or \ref{thm:SumTheorem2}, we can construct two versions of the processes $X$ and $\Breve{X}$ on a common probability space $(\Omega,\F,\PP)$ with initial conditions $\initialx$ and $\initialxbreve$, respectively and such that \eqref{eq:CouplingCondition} or  \eqref{eq:DominanceProperty_InnerProduct} or \eqref{eq:CouplingConditionBIS} hold. We denote these versions again by $X$ and $\Breve{X}$, and we observe that to show \eqref{eq:StochasticOrderingFPT_increasing}, it suffices to show that for an increasing set $\Gamma$, $\PP[\Breve{T}_{\Gamma} \leq T_\Gamma ]=1$ for $T_{\Gamma}$ and $\Breve{T}_{\Gamma}$ associated with these versions of $X$ and $\Breve{X}$. To see that this holds, let $\tilde{\Omega}$ be a set of probability one on which
			\begin{equation}
				\label{eq:ComparisonProofMFPT}
				X(t) \preccurlyeq_A \Breve{X}(t), \quad \text{ for all } t \geq 0    
			\end{equation}
			(this exists by \eqref{eq:CouplingCondition}, \eqref{eq:DominanceProperty_InnerProduct} or \eqref{eq:CouplingConditionBIS}). On $\{T_{\Gamma}=+\infty\}$, it is clear that $\Breve{T}_{\Gamma} \leq T_{\Gamma}$. For each $\omega \in \{T_{\Gamma} < + \infty\} \cap \tilde{\Omega}$ and $\eps > 0$ there is $\tau_{\eps}(\omega) \in [T_{\Gamma}(\omega),T_{\Gamma}(\omega) + \eps)$ such that $X(\tau_{\eps}(\omega)) \in \Gamma$ and by \eqref{eq:ComparisonProofMFPT}, $X(\tau_{\eps}(\omega)) \preccurlyeq_A \Breve{X}(\tau_{\eps}(\omega))$. And then, since $\Gamma$ is increasing, $\Breve{X}(\tau_{\eps}(\omega)) \in \Gamma$. It follows that $\Breve{T}_{\Gamma}(\omega) \leq T_{\Gamma}(\omega) + \eps$ and letting $\eps \to 0$ we obtain that $\Breve{T}_{\Gamma}(\omega) \leq T_{\Gamma}(\omega)$. It follows that $\PP[\Breve{T}_{\Gamma} \leq T_\Gamma]=1$. For the result on mean first passage times, let $\overline{F}_{T_{\Gamma}} := 1 - F_{T_\Gamma}$ and $\overline{F}_{\Breve{T}_{\Gamma}} := 1 - F_{\Breve{T}_\Gamma}$ represent the complementary cumulative distribution functions for $T_{\Gamma}$ and $\Breve{T}_{\Gamma}$, respectively. Observe that \eqref{eq:StochasticOrderingFPT_increasing} implies that $\overline{F}_{\Breve{T}_{\Gamma}} \leq \overline{F}_{T_\Gamma}$. For a non-negative random variable, the mean of the random variable is given by the Lebesgue integral of the complementary cumulative distribution function. Consequently, the mean first passage time for $\Breve{X}$ from $\initialxbreve$ to $\Gamma$ is given by $\int_{0}^{\infty} \overline{F}_{\Breve{T}_{\Gamma}}(t)dt \leq \int_0^{\infty} \overline{F}_{T_\Gamma}(t) dt$, where the latter is the mean first passage time for $X$ from $\initialx$ to $\Gamma$. If $\Gamma$ is decreasing, analogous arguments yield the results stated for that case.
		\end{proof}

		The second consequence of our results provides a comparison result for stationary distributions.
		
		\begin{theorem}
			\label{comparStationaryDistributionMAIN}
			Consider a non-empty set $\X \subseteq \Z_+^d$, a collection of distinct vectors $v_1,\ldots,v_n$ in $\Z^d \setminus \{0\}$ and two collections of non-negative functions on $\X$, $\rate=(\rate_1, \dots,\rate_n)$ and  $\Breve{\rate}= (\Breve{\rate}_1, \dots,\Breve{\rate}_n)$, such that \eqref{eq:LambdaNotOutofSpace} holds and the associated continuous-time Markov chains do not explode in finite time. Consider a matrix $A \in \R^{m \times d}$ with non-zero rows and suppose that at least one of the following holds:
			\begin{enumerate}
				\item[(i)]
				For every $x,y \in \X$ such that $x \preccurlyeq_A y$, conditions \eqref{eq:CouplingCondition1} and \eqref{eq:CouplingCondition2} are satisfied.
				\item[(ii)]
				The matrix $A$ has integer-valued entries and conditions $(i)$ and $(ii)$ in Theorem \ref{thm:InnerProductTheorem} are satisfied.
				\item[(iii)]
				The matrix $A$ has integer-valued entries and conditions $(i)$ and $(ii)$ in Theorem \ref{thm:SumTheorem2} are satisfied.
			\end{enumerate}
			Assume that the two continuous-time Markov chains on the set $\X$ with infinitesimal generators $Q$ and $\Breve{Q}$, associated with $\rate$ and $\Breve{\rate}$ respectively, are irreducible and positive recurrent on $\X$, and denote the associated stationary distributions by $\pi$ and $\Breve{\pi}$, respectively. If $\Gamma \subseteq \X$ is a non-empty set that is increasing in $\X$ with respect to $\preccurlyeq_A$, then
			\begin{equation}
				\label{eq:MonotonicityResultStatDist}
				\sum_{x \in \Gamma} \pi_x \leq  \sum_{x \in \Gamma} \Breve{\pi}_x.
			\end{equation}
			If $\Gamma \subseteq \X$ is a non-empty set that is decreasing in $\X$ with respect to $\preccurlyeq_A$, then
			\begin{equation}
				\label{eq:MonotonicityResultStatDist_II}
				\sum_{x \in \Gamma} \Breve{\pi}_x \leq  \sum_{x \in \Gamma} \pi_x.
			\end{equation}
		\end{theorem}
		
		\begin{proof}
			As in the proof of Theorem \ref{thm:comparison_MFPT}, we can construct two versions of the processes $X$ and $\Breve{X}$ on a common probability space $(\Omega,\F,\PP)$ for some pair of initial conditions $\initialx \preccurlyeq_A \initialxbreve$. If $\Gamma \subseteq \X$ is increasing, equation \eqref{eq:CouplingCondition} or \eqref{eq:DominanceProperty_InnerProduct} or \eqref{eq:CouplingConditionBIS} yields that $\PP(X(t) \in \Gamma) \leq \PP(\Breve{X}(t) \in \Gamma)$ for every $t \geq 0$. By letting $t \to \infty$ and observing that the stationary distribution is the steady-state distribution under our assumptions of irreducibility and positive recurrence, we obtain \eqref{eq:MonotonicityResultStatDist}. If $\Gamma$ is decreasing, an analogous argument yields  \eqref{eq:MonotonicityResultStatDist_II}.  
		\end{proof}
		
		\begin{remark}
			A special case of Theorems \ref{thm:comparison_MFPT} and \ref{comparStationaryDistributionMAIN} is when $\Gamma = \{x\}$ for some maximal or minimal element $x \in \X$.
		\end{remark}
		
		In the next section, we give examples which illustrate Theorem \ref{thm:InnerProductTheorem} (see Examples \ref{exEK}, \ref{exEK2}, \ref{ex:HistoneModification} and \ref{ex:HistoneModificationAndProtein}), Theorem \ref{thm:SumTheorem2} (see Example \ref{ex:Braess}), Theorem \ref{thm:comparison_MFPT} and Theorem \ref{comparStationaryDistributionMAIN} for continuous-time Markov chains that are stochastic chemical reaction networks. For Examples \ref{exEK}, \ref{exEK2} and \ref{ex:Braess}, the state space $\X$ will be a stoichiometric compatibility class $z + \L$ intersected with $\Z^d_+$. For Examples \ref{ex:HistoneModification} and \ref{ex:HistoneModificationAndProtein}, we work with reduced Markov chains and the state space $\X$ will be a projection of a suitable higher dimensional stoichiometric compatibility class $z + \L$ intersected with $\Z^d_+$.

	\section{Examples}
	\label{sec:examples}
	
	In this section, we apply the theoretical tools developed in the paper to several examples. While in Examples \ref{exEK}, \ref{ex:Braess} and \ref{ex:HistoneModification} the Markov chains analyzed have a finite state space, in Examples \ref{exEK2} and \ref{ex:HistoneModificationAndProtein} the Markov chains have a countably infinite state space, but it is straightforward to verify that they do not explode (see SI - Sections S.1.2 and S.1.3, respectively). The choice of matrix $A$ in each example is based on the specific monotonicity relationship of interest. While for simpler cases the choice of $A$ is straightforward, for more complicated systems the choice can be more subtle. In many cases, in order to study the monotonicity properties for the stochastic behavior of our system, we can rely on Theorem \ref{thm:InnerProductTheorem}, which provides a reasonable approach to narrow down the choices for suitable $A$. The approach consists in solving, for each row $i$, the system of equations $\sum_{k=1}^d A_{ik}(v_j)_k=b_{ij}$, with $b_{ij}$ equal to $1,-1,$ or $0$ depending, based on the monotonicity relationship of interest, whether we expect that the Markov chain transition in the direction $v_j$ leads inside, outside, or is parallel to the boundary of the region $K_A+x$. Finally, it is worth noticing that, while all the following examples compare two identical reaction networks with different rate constants, our theory can also be applied to compare two different reaction networks as long as they have the same net reaction vectors $\{v_j\}_{j=1}^n$.

		\begin{example} \textbf{Enzyme kinetics I}\label{exEK}\\
			\begin{figure}[t!]
				\centering
				\includegraphics[scale=0.44]{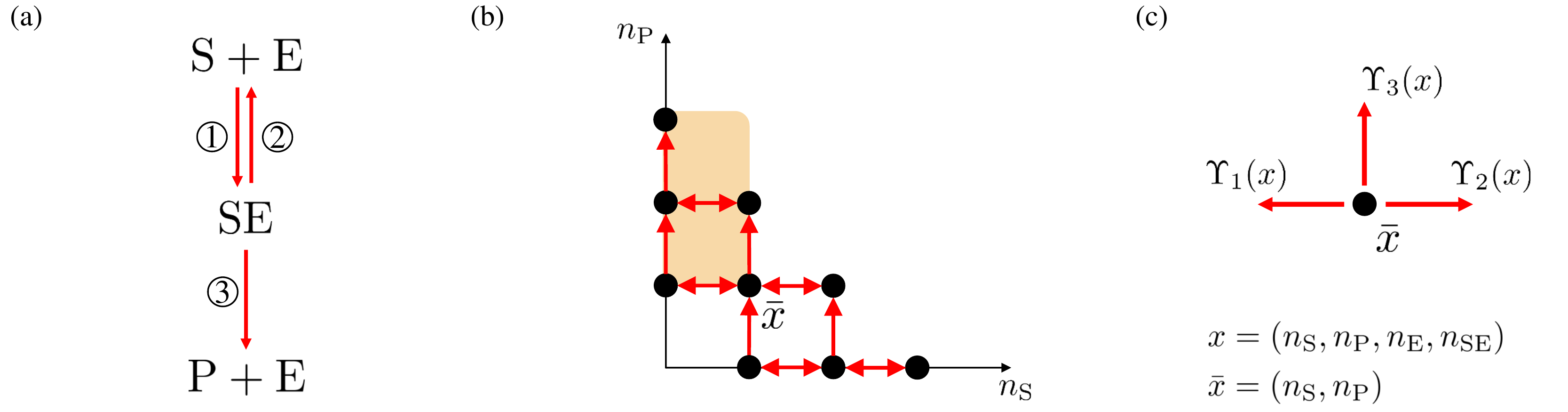}
				\caption{\small { \bf Reaction model and corresponding Markov chain for enzymatic kinetics I example.} (a) Chemical reaction system. The numbers on the arrows correspond to the associated reactions. (b) Projected Markov chain graph for one stoichiometric compatibility class with two conservation laws $n_{\mathrm{S}}+n_{\mathrm{P}}+n_{\mathrm{SE}}=\mathrm{S_{tot}}=3$ and $n_{\mathrm{E}}+n_{\mathrm{SE}}=\mathrm{E_{tot}}=2$. The projection takes a state $x=(n_{\mathrm{S}},n_{\mathrm{P}},n_{\mathrm{E}},n_{\mathrm{SE}})=(n_{\mathrm{S}},n_{\mathrm{P}},\mathrm{E_{tot}}-\mathrm{S_{tot}}+n_{\mathrm{S}}+n_{\mathrm{P}},\mathrm{S_{tot}}-n_{\mathrm{S}}-n_{\mathrm{P}})$ to $\bar{x} = (n_{\mathrm{S}},n_{\mathrm{P}})$. We use black dots to represent the states, red double-ended arrows to represent transitions in both directions and red single-ended arrows to represent transitions in one direction. Note that $\bar{x} = (0,0)$ is not a vertex in the graph because $0 \leq n_{\mathrm{E}} = \mathrm{E_{tot}}-\mathrm{S_{tot}}+n_{\mathrm{S}}+n_{\mathrm{P}}$, and so $n_{\mathrm{S}}+n_{\mathrm{P}}\ge 3-2=1$. We use orange to highlight the projection of the region $K_A+x$ intersected with the stoichiometric compatibility class, where $A$ is defined in  \eqref{matrixAEK}. (c) The projections of the directions of the possible transitions of the Markov chain. The transition rates $\Upsilon_1(x)$, $\Upsilon_2(x)$, and $\Upsilon_3(x)$ are defined in \eqref{ratesEK}.  }
				\label{fig:EK}
			\end{figure}
			Let us consider a classic model of enzyme kinetics (see Michaelis \& Menten \cite{MM1913} and Kang et al. \cite{KangKhudaBukhshKoepplRempala}), where an enzyme catalyzes the conversion of a substrate to a product. The species considered here are substrate (S), enzyme (E), intermediate enzyme-substrate complex (SE), and product (P), and the chemical reaction system is depicted in Fig. \ref{fig:EK}(a). We are interested in how the rate constant $\kappa_3$ affects the time to convert the substrate to the final product.

			To this end, let us first introduce the set of species $\Ss =\{\mathrm{S},\mathrm{P},\mathrm{E},\mathrm{SE}\}$, and the set of reactions $\Rs=\{(v^-_1,v^+_1),(v^-_2,v^+_2),(v^-_3,v^+_3)\}$, where $v^-_1= v^+_2= (1,0,1,0)^T$, $v^+_1=v^-_2=v^-_3=(0,0,0,1)^T$, $v^+_3=(0,1,1,0)^T$, where $T$ denotes transpose. At a given time, let the counts of each of the species S, P, E and SE be denoted by $n_{\mathrm{S}}$, $n_{\mathrm{P}}$, $n_{\mathrm{E}}$ and $n_{\mathrm{SE}}$, respectively. The state of the associated Markov chain is $(n_{\mathrm{S}},n_{\mathrm{P}},n_{\mathrm{E}},n_{\mathrm{SE}})$. The potential transitions of the Markov chain are in three possible directions:
			$$v_1=v^+_1-v^-_1= (-1,0,-1,1)^T, v_2=v^+_2-v^-_2= (1,0,1,-1)^T, v_3=v^+_3-v^-_3=(0,1,1,-1)^T.$$
			Fixing integers $\mathrm{S_{tot}},\mathrm{E_{tot}}>0$, we have a stoichiometric compatibility class $z+\L$ with $z = (\mathrm{S_{tot}},0,\mathrm{E_{tot}},0)$ and $\L := \vspan\{v_1,v_2,v_3\}$, which is contained in a two-dimensional affine subspace of four dimensional space. Then, the state space of the Markov chain is
			$$\X = (z+\L) \cap \Z_+^4 = \{(x_1,x_2,x_3,x_4) \in \Z_+^4 | x_1+x_2+x_4=\mathrm{S_{tot}}, x_3+x_4=\mathrm{E_{tot}} \}.$$
			The two constraints described in the last expression for $\X$ characterize the two linearly independent conservation laws for this chemical reaction system: $n_{\mathrm{S}}+n_{\mathrm{P}}+n_{\mathrm{SE}}=\mathrm{S_{tot}}$ and $n_{\mathrm{E}}+n_{\mathrm{SE}}=\mathrm{E_{tot}}$.

			Given a state $x=(x_1,x_2,x_3,x_4) \in \X$, following mass-action kinetics, the infinitesimal transition rates are
			\begin{equation}\label{ratesEK}
				\begin{aligned}
					&\rate_1(x)= \kappa_1 x_1 x_3,\;\;\rate_2(x)= \kappa_2 x_4,\;\;
					&\rate_3(x)= \kappa_3 x_4,
				\end{aligned}
			\end{equation}
			for constants $\kappa_1,\kappa_2,\kappa_3 > 0$. Here, we have used $\kappa_j$ as an abbreviation for $\kappa_{(v_j^-,v_j^+)}$, $j=1,2,3$. We will use similar abbreviations in the other examples too.

			We note that the projected process $(X_1,X_2)(\cdot)$ is still a continuous-time Markov chain, and we could apply our theory to it. However, when the functions $\Upsilon_j$, $j=1,2,3$, are written in terms of these two components, they will have a more complex, non-mass action form. Here we apply our theory directly to our four dimensional Markov chain. For the purpose of visualization, Fig. \ref{fig:EK}(b) shows the two dimensional projection of the four dimensional Markov chain graph for one stoichiometric compatibility class. In Examples \ref{exEK2} and \ref{ex:Braess}, we also analyze Markov chains without projections, and in Examples \ref{ex:HistoneModification} and \ref{ex:HistoneModificationAndProtein}, we analyze projected Markov chains.
			
			In order to study how the rate constant $\kappa_3$ affects the time to convert the substrate to the final product, let us define the state $(0,\mathrm{S_{tot}},\mathrm{E_{tot}},0)$ associated with $n_{\mathrm{P}}=\mathrm{S_{tot}}$ as $p$, the state $(\mathrm{S_{tot}},0,\mathrm{E_{tot}},0)$
			associated with $n_{\mathrm{S}}=\mathrm{S_{tot}}$ as $s$, and the mean first passage time to reach the state $p$, starting from $s$, as $\E_{s} [T_p]$. We will verify that the assumptions of Theorems \ref{thm:InnerProductTheorem},  \ref{thm:comparison_MFPT}  hold and exploit them to determine how $\kappa_3$ affects $\E_{s} [T_p]$. To this end, define the matrix 
			\begin{equation}\label{matrixAEK}
				A= \begin{bmatrix}
					-1 & 0 & 0 & 0\\
					0 & 1 & 0 & 0
				\end{bmatrix}
			\end{equation}
			and consider the preorder $x \preccurlyeq_A y$, defined by $A(y-x)\ge 0$, and the set $K_A +x = \{ w \in \R^4 \:|\: x \preccurlyeq_A w \}$. Let us also consider the infinitesimal transition rates $\Breve{\rate}_1(x),\Breve{\rate}_2(x)$ and $\Breve{\rate}_3(x)$ defined as for $\rate_1(x),\rate_2(x)$ and $\rate_3(x)$, but with $\Breve{\kappa}_1=\kappa_1$, $\Breve{\kappa}_2=\kappa_2$, $\Breve{\kappa}_3>\kappa_3$ in place of $\kappa_1$, $\kappa_2$, $\kappa_3$, respectively. Condition $(i)$ of Theorem \ref{thm:InnerProductTheorem} (i.e., for every $1\leq j \leq n$, the vector $Av_j$ has entries in $\{-1,0,1\}$) holds since $Av_1=(1,0)^T, Av_2=(-1,0)^T$ and $Av_3=(0,1)^T$. Condition $(ii)$ of Theorem \ref{thm:InnerProductTheorem} also holds, as shown in the paragraph below.
			
			\textbf{Verification of condition $(ii)$ of Theorem \ref{thm:InnerProductTheorem}.} We first consider $x \in \X$ and $y \in \partial_1(K_A+x) \cap \X$, where
			\small
			\begin{eqnarray*}
				&& \partial_1(K_A+x) \cap \X\\
				&=& \{ w \in \Z^4_+ \:|\: x_1 = w_1, x_2 \leq w_2 \} \cap \X \\
				&=& \{ w \in \Z^4_+ \:|\: x_1 = w_1, x_2 \leq w_2, x_1 + x_2 + x_4 = w_1 + w_2 + w_4=\mathrm{S_{tot}}, x_3 + x_4 = w_3 + w_4=\mathrm{E_{tot}} \} \\
				&=& \{ w \in \Z_+^4 \:|\: x_1 = w_1, x_2 \le w_2, x_3 \le w_3, x_4 \ge w_4,  w_1 + w_2 + w_4=\mathrm{S_{tot}}, w_3 + w_4=\mathrm{E_{tot}}\}\\
				&=& \{ w \in \X \:|\: x_1 = w_1, x_2 \le w_2, x_3 \le w_3, x_4 \ge w_4\}.
			\end{eqnarray*}
			\normalsize
			Since $\inn{A_{1\bullet},v_1}=1$, $\inn{A_{1\bullet},v_2}=-1$, we need to check that $\rate_1(x)\le\Breve{\rate}_1(y)$ and $ \rate_2(x)\ge\Breve{\rate}_2(y)$. The first inequality holds because $y \in \partial_1(K_A+x) \cap \X$ implies $x_1=y_1$ and $x_3 \leq y_3$ so that $\rate_1(x) = \kappa_1 x_1 x_3 \le \kappa_1 y_1 y_3 = \Breve{\kappa}_1 y_1 y_3 = \Breve{\rate}_1(y)$. The second inequality holds because $y \in \partial_1(K_A+x) \cap \X$ implies $x_4 \geq y_4$ so that $\rate_2(x) = \kappa_2 x_4 \ge \kappa_2 y_4 = \Breve{\kappa}_2 y_4 = \Breve{\rate}_2(y)$.
			Secondly, we consider $x \in \X$, $y \in \partial_2(K_A+x) \cap \X = \{ w \in \X \:|\: x_1 \ge w_1, x_2 = w_2, x_3 \ge w_3, x_4 \le w_4\}$. Then, since $\inn{A_{2\bullet},v_3}=1$, we need to check that $\rate_3(x)\le\Breve{\rate}_3(y)$. This holds because $y \in \partial_2(K_A+x) \cap \X$ implies $x_4 \leq y_4$ so that $\rate_3(x) = \kappa_3 x_4 \le \kappa_3 y_4 \le \Breve{\kappa}_3 y_4 = \Breve{\rate}_3(y)$.

			Since all of the hypotheses of Theorem \ref{thm:InnerProductTheorem} hold, we can conclude that, for each $\initialx,\initialxbreve \in \X$ with $\initialx \preccurlyeq_A \initialxbreve$, there exists a probability space $(\Omega,\F,\PP)$ with two Markov chains $X = \{X(t), \: t \geq 0\}$ and $\Breve{X}=\{\Breve{X}(t), \: t \geq 0\}$ associated with $\rate$ and $\Breve{\rate}$, respectively, such that $X(0)=\initialx$, $\Breve{X}(0)=\initialxbreve$ and $$\PP\left[X(t) \preccurlyeq_A \Breve{X}(t) \text{ for every } t \geq 0 \right]=1.$$ Furthermore, applying Theorem \ref{thm:comparison_MFPT} with the set $\Gamma = \{p\}=\{ (0,\mathrm{S_{tot}},\mathrm{E_{tot}},0) \}$, which is increasing in $\X$ with respect to $\preccurlyeq_A$, we see that the mean first passage time from $s$ to $p$, $\E_{s} [T_p]$, is a decreasing function of $\kappa_3$.
			
			Because the Markov chain has one absorbing state, $p$, per stoichiometric compatibility class, the stationary distribution on a given stoichiometric compatibility class is trivial, and hence so too are its monotonicity properties.			
			
		\end{example}

		\begin{example} \textbf{Enzyme kinetics II}\\
			\label{exEK2}
			\begin{figure}[t!]
				\centering
				\includegraphics[scale=0.44]{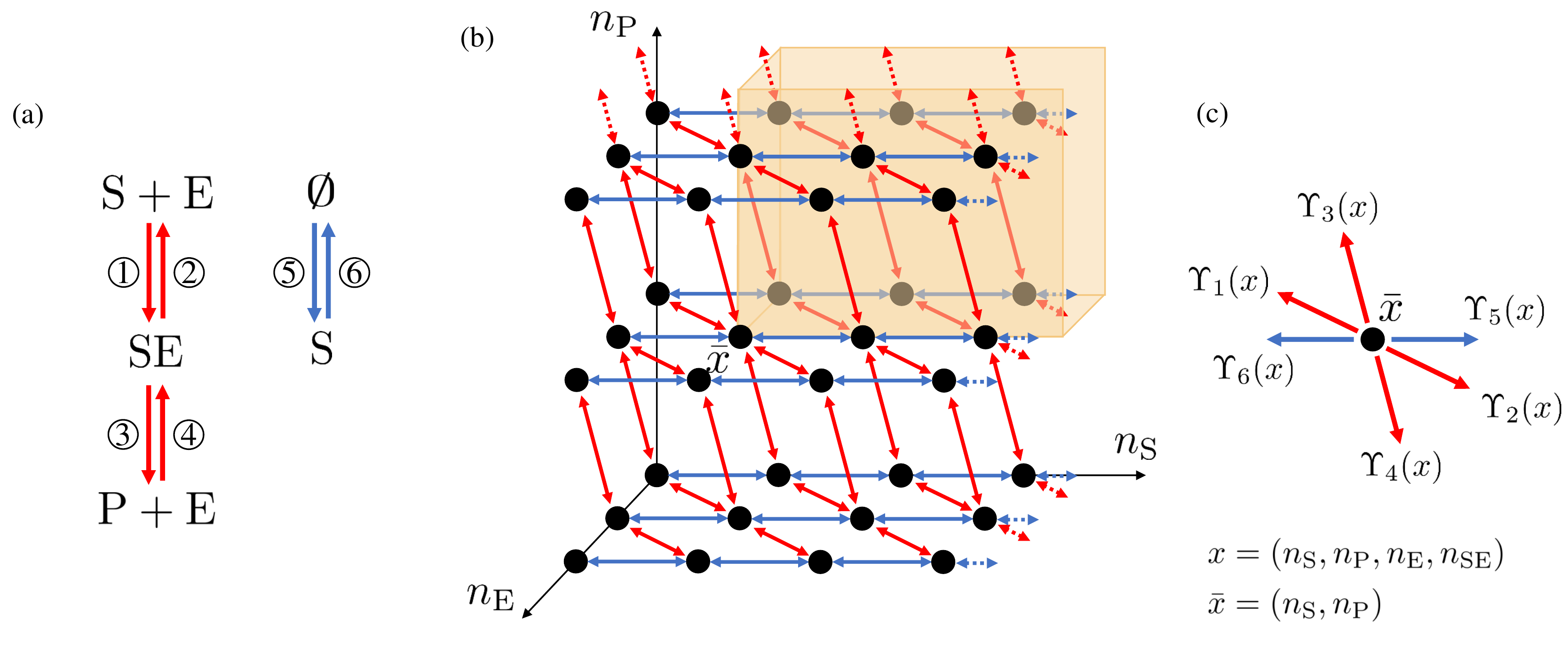}
				\caption{\small { \bf Reaction model and corresponding Markov chain for enzymatic kinetics II example.} (a) Chemical reaction system. The numbers on the arrows correspond to the associated reactions. (b) Projected Markov chain graph for one stoichiometric compatibility class with the conservation law $n_{\mathrm{E}}+n_{\mathrm{SE}}=\mathrm{E_{tot}}=2$. The projection takes a state $x=(n_{\mathrm{S}},n_{\mathrm{P}},n_{\mathrm{E}},n_{\mathrm{SE}})=(n_{\mathrm{S}},n_{\mathrm{P}},n_{\mathrm{E}},\mathrm{E_{tot}}-n_{\mathrm{E}}) \in \X$ to $\bar{x} = (n_{\mathrm{S}},n_{\mathrm{P}},n_{\mathrm{E}})\in \Z_+^3:0\le n_{\mathrm{E}}\le 2$. Here, we use black dots to represent the states, red double-ended arrows to represent transitions in both directions associated with the reactions represented by the red arrows in (a) and blue double-ended arrows to represent transitions in both directions associated with the reactions represented by the blue arrows in (a). We use dotted arrowed-lines to indicate that the pattern of Markov chain transitions extends to infinity. We use orange to highlight the projections of the region $K_A+x$ intersected with the stoichiometric compatibility class, where $A$ is defined in \eqref{matrixAEK2}. (c) The projections of the directions of the possible transitions of the Markov chain within a stoichiometric compatibility class. The transition rates $\Upsilon_i(x)$, $i=1,2,3,4,5,6$, are defined in (\ref{ratesEK2}).
				}
				\label{fig:EK2}
			\end{figure}
			Let us consider an extension of the enzymatic kinetics model introduced in the previous example, in which the substrate S can enter and leave the system and the product can revert to the substrate. This is a simplified version of the enzymatic kinetics considered by Anderson et al. \cite{2010Anderson}. 
			The chemical reaction system is depicted in Fig. \ref{fig:EK2}(a).
			Now, for this case study, we first determine how the reaction rate constant $\kappa_5$ affects the stochastic behavior of the system and then we will study properties of the system with respect to initial conditions. To this end, let us introduce the set of species $\Ss =\{\mathrm{S},\mathrm{P},\mathrm{E},\mathrm{SE}\}$, and, similar to Example \ref{exEK}, we let $(n_{\mathrm{S}},n_{\mathrm{P}},n_{\mathrm{E}},n_{\mathrm{SE}})$ be the state of the Markov chain that records the number of molecules of each species. The potential transitions of the Markov chain are in six possible directions, $v_j$ for $j=1,...,6$, where $v_1=-v_2=(-1,0,-1,1)^T$, $v_3=-v_4= (0,1,1,-1)^T$, and $v_5=-v_6=(1,0,0,0)^T$ (see SI-Section S.2.1 for the derivation of the $v_j$, $j=1,...,6$). Since there is one linearly independent conservation law in this chemical reaction system: $n_{\mathrm{E}}+n_{\mathrm{SE}}=\mathrm{E_{tot}}$, each stoichiometric compatibility class is contained in a three-dimensional affine subspace of four dimensional space, denoted as $z+\L$, where $z=(0,0,\mathrm{E_{tot}},0)$ and $\L := \vspan\{v_1,v_3,v_5\}$, with fixed integer $\mathrm{E_{tot}}>0$. Then, we can choose the state space of the Markov chain to be $\X = (z+\L) \cap \Z_+^4 = \{(x_1,x_2,x_3,x_4) \in \Z_+^4 |x_3+x_4=\mathrm{E_{tot}} \}$. Furthermore, given a state $x=(x_1,x_2,x_3,x_4) \in \X$, following mass-action kinetics, the associated infinitesimal transition rates are given by
			\begin{equation}\label{ratesEK2}
				\begin{aligned}
					&\rate_1(x)= \kappa_1 x_1 x_3,\;\;\rate_2(x)= \kappa_2 x_4,\;\;\rate_3(x)= \kappa_3 x_4,\\
					&\rate_4(x)= \kappa_4 x_2 x_3,\;\;\rate_5(x)= \kappa_5,\;\;\rate_6(x)= \kappa_6 x_1,
				\end{aligned}
			\end{equation}
			for $\kappa_1,\kappa_2,\kappa_3,\kappa_4,\kappa_5,\kappa_6 > 0$. As in Example \ref{exEK}, we apply our theory directly to our four dimensional Markov chain, but, for the purpose of illustration, Fig. \ref{fig:EK2}(b) shows the three dimensional projection of the Markov chain graph for one stoichiometric compatibility class.
			Now, for the first analysis (determining how $\kappa_5$ affects the stochastic behavior of the system), we verify that the assumptions of Theorems \ref{thm:InnerProductTheorem} and \ref{comparStationaryDistributionMAIN} hold and use them to determine how $\kappa_5$ affects the stationary distribution. 
			
			To this end, define the matrix
			\begin{equation}\label{matrixAEK2}
				A= \begin{bmatrix}
					1 & 0 & 0 & 0\\
					0 & 1 & 0 & 0\\
					0 & 0 & -1 & 0
				\end{bmatrix}
			\end{equation}
			and consider the preorder $x \preccurlyeq_A y$, defined by $A(y-x)\ge 0$. For $x \in \X$, $K_A +x = \{ w \in \R^4 \:|\: x \preccurlyeq_A w \}$. Furthermore, let us consider the infinitesimal transition rates $\Breve{\rate}_1(x)$, $\Breve{\rate}_2(x)$, $\Breve{\rate}_3(x)$, $\Breve{\rate}_4(x)$, $\Breve{\rate}_5(x)$ and $\Breve{\rate}_6(x)$ defined as for $\rate_1(x),\rate_2(x),\rate_3(x),\rate_4(x),\rate_5(x)$ and $\rate_6(x)$, but with $\Breve{\kappa}_i$ in place of $\kappa_i$, where $\Breve{\kappa}_i=\kappa_i$, for $i = 1,2,3,4,6$, and $\Breve{\kappa}_5\ge \kappa_5$. Given that $Av_1=(-1,0,1)^T$, $Av_2=(1,0,-1)^T$, $Av_3=(0,1,-1)^T$, $Av_4=(0,-1,1)^T$, $Av_5=(1,0,0)^T$ and $Av_6=(-1,0,0)^T$, we have that condition $(i)$ of Theorem \ref{thm:InnerProductTheorem} holds. Condition $(ii)$ of that theorem also holds, as shown in the next paragraph.
			
			\textbf{Verification of condition $(ii)$ of Theorem \ref{thm:InnerProductTheorem}.} First consider $x \in \X$ and $y \in \partial_1(K_A+x) \cap \X$, where $\partial_1(K_A+x) \cap \X = \{ w \in \X \:|\: x_1= w_1, x_2\le w_2, x_3\ge w_3, x_4\le w_4 \}$. Since $\inn{A_{1\bullet},v_2}=\inn{A_{1\bullet},v_5}=1$ and $\inn{A_{1\bullet},v_1}=\inn{A_{1\bullet},v_6}=-1$, we need to check that $\rate_1(x)\ge\Breve{\rate}_1(y),\rate_6(x) \ge \Breve{\rate}_6(y),\rate_2(x)\le\Breve{\rate}_2(y),$ and $\rate_5(x)\le\Breve{\rate}_5(y)$. Given that $y \in \partial_1(K_A+x) \cap \X$, the first inequality holds because $\rate_1(x)=\kappa_1 x_1 x_3\ge \kappa_1 y_1 y_3=\Breve{\kappa}_1 y_1 y_3= \Breve{\rate}_1(y)$, the second inequality holds because $\rate_6(x) =\kappa_6 x_1=\kappa_6 y_1=\Breve{\kappa}_6 y_1 =\Breve{\rate}_6(y)$, the third inequality holds because $\rate_2(x)=\kappa_2 x_4\le \kappa_2 y_4=\Breve{\kappa}_2 y_4=\Breve{\rate}_2(y)$, and the fourth inequality holds because $\rate_5(x)= \kappa_5\le\Breve{ \kappa}_5=\Breve{\rate}_5(y)$. 
			Secondly, we consider $x \in \X$ and $y \in \partial_2(K_A+x) \cap \X = \{ w \in \X \:|\: x_1\le w_1, x_2=w_2, x_3\ge w_3, x_4\le w_4 \}$. Given that $\inn{A_{3\bullet},v_3}=1$ and $\inn{A_{3\bullet},v_4}=-1$, we need to check that $\rate_4(x)\ge\Breve{\rate}_4(y)$ and $\rate_3(x)\le\Breve{\rate}_3(y)$. The first inequality holds because $\rate_4(x)=\kappa_4 x_2 x_3\ge \kappa_4 y_2 y_3=\Breve{\kappa}_4 y_2 y_3=\Breve{\rate}_4(y)$ and the second inequality holds because $\rate_3(x)=\kappa_3 x_4 \le \kappa_3 y_4 = \Breve{\kappa}_3 y_4=\Breve{\rate}_3(y)$. 
			Finally, consider $x \in \X$ and $y \in \partial_3(K_A+x) \cap \X = \{ w \in \X \:|\: x_1\le w_1, x_2 \le w_2, x_3=w_3, x_4= w_4 \}$. Since $\inn{A_{3\bullet},v_1}=\inn{A_{3\bullet},v_4}=1$ and $\inn{A_{3\bullet},v_2}=\inn{A_{3\bullet},v_3}=-1$, we need to check that $\rate_2(x)\ge \Breve{\rate}_2(y)$, $\rate_3(x)\ge \Breve{\rate}_3(y)$, $\rate_1(x) \le \Breve{\rate}_1(y)$, and $\rate_4(x)\le\Breve{\rate}_4(y)$. Indeed, we have that $\rate_2(x)=\kappa_2 x_4 = \kappa_2 y_4= \Breve{\kappa}_2 y_4=\Breve{\rate}_2(y)$, $\rate_3(x) = \kappa_3 x_4 = \kappa_3 y_4 =\Breve{\kappa}_3 y_4 =  \Breve{\rate}_3(y)$, $\rate_1(x) = \kappa_1 x_1 x_3 \le \kappa_1 y_1 y_3 = \Breve{\kappa}_1 y_1 y_3 = \Breve{\rate}_1(y)$, and $\rate_4(x)=\kappa_4 x_2x_3\le \kappa_4 y_2y_3= \Breve{\kappa}_4 y_2y_3 = \Breve{\rate}_4(y)$.

			Thus, all of the hypothesis of Theorem \ref{thm:InnerProductTheorem} are verified, and so, for each pair $\initialx,\initialxbreve \in \X$ satisfying  $\initialx \preccurlyeq_A \initialxbreve$, there exists a probability space $(\Omega,\F,\PP)$ with two Markov chains $X = \{X(t), \: t \geq 0\}$ and $\Breve{X}=\{\Breve{X}(t), \: t \geq 0\}$ associated with $\rate$ and $\Breve{\rate}$, respectively, such that $X(0)=\initialx$, $\Breve{X}(0)=\initialxbreve$ and $\PP\left[X(t) \preccurlyeq_A \Breve{X}(t) \text{ for every } t \geq 0 \right]=1$.
			The Markov chains $X,\Breve{X}$ are irreducible and positive recurrent  (see SI - Section S.1.1).
			Furthermore, for the increasing set in $\X$ with respect to $\preccurlyeq_A$ defined as $\Gamma(x) = \{ w \in \X \:|\: x_1 \le w_1, x_2 \le w_2, x_3 \ge w_3, x_4 \le w_4 \}$, we can apply Theorem \ref{comparStationaryDistributionMAIN} and obtain that $\sum_{w \in \Gamma(x)} \pi_w \leq  \sum_{w \in \Gamma(x)} \Breve{\pi}_w$. 
			Loosely speaking, this means that increasing $\kappa_5$ causes the stationary distribution $\pi(x)$ to shift mass towards states characterized by lower $x_3$ and higher $x_1,x_2$ and $x_4$. 
			
			For this specific case, in which we have a stochastic chemical reaction network associated with a complex balanced dynamical system, an explicit expression for the stationary distribution can be obtained by applying Theorem 4.1 in Anderson et al. \cite{2010Anderson}. Analysis of this formula would provide results in agreement with the ones obtained by applying the theoretical tools developed in this paper. Specifically, $\pi_x$ can be written as a product of two Poisson distributions and a binomial distribution, i.e.,
		\begin{equation}\label{statdistr}
			\pi_x=\left(e^{-c_1}\frac{c_1^{x_1}}{x_1!}\right)\left(e^{-c_2}\frac{c_2^{x_2}}{x_2!}\right)\left(\mathrm{E_{tot}}!\frac{c_3^{x_3}}{x_3!}\frac{c_4^{x_4}}{x_4!}\right),\;\;\;x\in \X,
		\end{equation}
	in which $(c_1,c_2,c_3,c_4)$ represents the complex balanced equilibrium for the deterministic model, where
	\begin{equation}\label{cs}
		c_1=\frac{\kappa_5}{\kappa_6},\;\;c_2=\frac{\kappa_1\kappa_3\kappa_5}{\kappa_2\kappa_4\kappa_6},\;\;c_3=\frac{1}{1+\frac{\kappa_1\kappa_5}{\kappa_2\kappa_6}},\; \mathrm{and}\; c_4=\frac{\frac{\kappa_1\kappa_5}{\kappa_2\kappa_6}}{1+\frac{\kappa_1\kappa_5}{\kappa_2\kappa_6}}.
	\end{equation}
	In most cases, it is not possible to derive an analytical formula for the stationary distribution, but our theorems can still be applied and then monotonicity properties for $\pi$ can still be determined even without an explicit expression for $\pi$. For instance, in the context of the above example, if the infinitesimal transition rates $\rate_i$ do not follow mass-action kinetics, the deficiency zero theorem and Theorem 4.1 in Anderson et al. \cite{2010Anderson} do not apply. Nevertheless, our theory can still be easily applied to study monotonicity properties for sample paths and stationary distributions.
	
	As pointed out in Remark \ref{remark:incond}, we can also exploit our theoretical tools to determine monotonicity properties of the system with respect to the initial conditions. 
	For this, suppose that $\Breve{\kappa}_i=\kappa_i$ for $i = 1,2,3,4,5,6$. Then, by the analysis above, Theorem \ref{thm:InnerProductTheorem} holds and yields monotonically (with preorder induced by the matrix $A$) with respect to the initial conditions.

\end{example}   

\begin{example}{\bf A network topology arising in Braess' paradox}
\label{ex:Braess}
A natural question in synthetic biology may involve the prediction of whether an engineered biological circuit with additional reactions will lead to the desired effect of accelerating the process or unexpected behaviors. Now, we consider an example inspired by Braess' paradox, which arises from transportation networks, where adding one or more roads to a road network can slow down overall traffic flow through the network (see Braess \cite{braess} and see also a related state-dependent queuing network model in Calvert et al. \cite{braessQueuing}). A simple network of this type is one where there are two routes to get from the start to the final destination, and adding a linkage road between the routes can in some cases increase travel times. Fig. \ref{fig:braess}(a) shows a reaction network analogue of the Braess' network topology. Of course, our chemical reaction network is a little different from a road network since there is no congestion nor competition between molecules and pathways are chosen randomly with certain probabilities instead of routing decisions being based on the number of cars on the routes. Nevertheless, the example considered here is interesting because adding a reaction to cross-link two pathways might intuitively be interpreted as a detour and be expected to increase the time to the final destination, while this is sometimes not the case in this example. 
\begin{figure}[t!]
	\centering
	\includegraphics[scale=0.44]{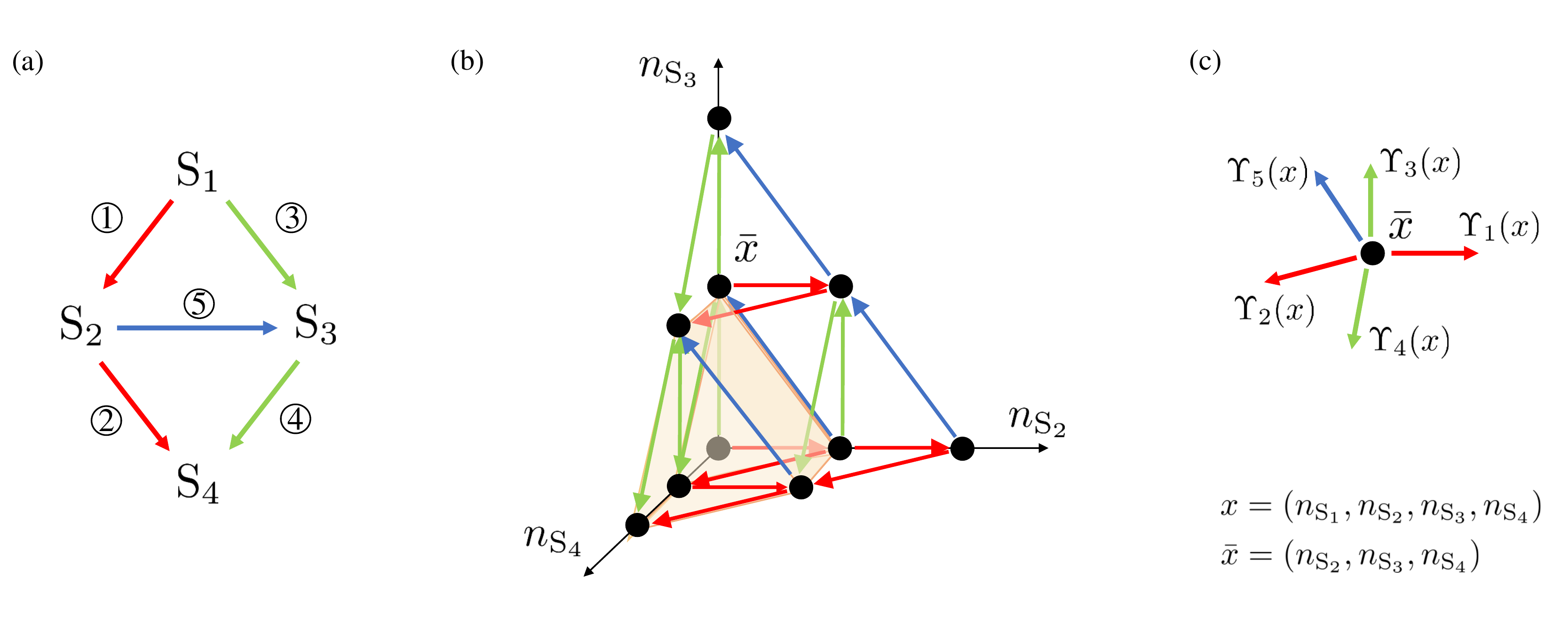}
	\caption{\small { \bf Circuit inspired by Braess' paradox and corresponding Markov chain.} (a) Chemical reaction system. The numbers on the arrows correspond to the associated reactions. (b) Projected Markov chain graph for one stoichiometric compatibility class with the conservation law $n_{\mathrm{S_1}}+n_{\mathrm{S_2}}+n_{\mathrm{S_3}}+n_{\mathrm{S_4}}=\mathrm{S_{tot}}=2$. The projection takes a state $x=(n_{\mathrm{S_1}},n_{\mathrm{S_2}},n_{\mathrm{S_3}},n_{\mathrm{S_4}})=(\mathrm{S_{tot}}-n_{\mathrm{S_2}}-n_{\mathrm{S_3}}-n_{\mathrm{S_4}},n_{\mathrm{S_2}},n_{\mathrm{S_3}},n_{\mathrm{S_4}}) \in \X$ to $\bar{x} = (n_{\mathrm{S_2}},n_{\mathrm{S_3}},n_{\mathrm{S_4}})$. Here, we use black dots to represent the states and red (blue, green) arrows to represent transitions in directions associated with the reactions represented by the red (blue, green) arrows in (a). We use orange to highlight the projection of the region $K_A+x$ intersected with the stoichiometric compatibility class, where $A$ is defined in  \eqref{matrixAbraess}. (c) The projections of the directions of the possible transitions of the Markov chain within a stoichiometric compatibility class. The transition rates $\Upsilon_i(x)$, $i=1,2,3,4,5$, are given in (\ref{ratesBraess}).
	}
	\label{fig:braess}
\end{figure} 

The chemical reaction system is depicted in Fig. \ref{fig:braess}(a), which involves four species $\Ss =\{\mathrm{S_1}, \mathrm{S_2}, \mathrm{S_3}, \mathrm{S_4}\}$. The state of the Markov chain is $(n_{\mathrm{S_1}}, n_{\mathrm{S_2}}, n_{\mathrm{S_3}}, n_{\mathrm{S_4}})$ where $n_\mathrm{S_i}$ is the number of copies of $\mathrm{S_i}$ for $i=1,2,3,4$. The potential transitions of the Markov chain are in five possible directions, $v_j$, $j=1,...,5$, where $v_1=(-1,1,0,0)^T$, $v_2=(0,-1,0,1)^T$, $v_3=(-1,0,1,0)^T$, $v_4=(0,0,-1,1)^T$ and $v_5=(0,-1,1,0)^T$ (see SI-Section S.2.2 for the derivation of the $v_j$, $j=1,...,5$). Fixing an integer $\mathrm{S_{tot}}>0$, the associated stoichiometric compatibility class is $z+\L$ with $z = (\mathrm{S_{tot}},0,0,0)$ and $\L := \vspan\{v_1,v_2,v_3,v_4,v_5\}$. The set $z+\L$ is a three-dimensional affine subspace of four dimensional space. We choose the state space of our Markov chain to be $\X = (z+\L) \cap \Z_+^4 = \{(x_1,x_2,x_3,x_4) \in \Z_+^4 | x_1+x_2+x_3+x_4=\mathrm{S_{tot}}\}$. The constraint introduced in the last expression for $\X$ follows from the conservation law in this chemical reaction system, that is $n_{\mathrm{S_1}}+n_{\mathrm{S_2}}+n_{\mathrm{S_3}}+n_{\mathrm{S_4}}=\mathrm{S_{tot}}$. 
Given a generic state $x=(x_1,x_2,x_3,x_4)$, following mass-action kinetics, the infinitesimal transition rates are
\begin{equation}\label{ratesBraess}
	\rate_1(x)= \kappa_1 x_1,\;\;\rate_2(x)= \kappa_2 x_2,\;\;\rate_3(x)= \kappa_3 x_1,\;\;\rate_4(x)= \kappa_4 x_3,\;\;\rate_5(x)= \kappa_5 x_2.
\end{equation}
For the purpose of illustration, Fig. \ref{fig:braess}(b) shows the three dimensional projection of the Markov chain graph for one stoichiometric compatibility class.

A natural question is how the time $T_{(0,0,0,\mathrm{S_{tot}})}$ to reach the state $(0,0,0,\mathrm{S_{tot}})$ from $(\mathrm{S_{tot}},0,0,0)$ depends on the rate constants $\kappa_1$,$\kappa_2$,$\kappa_3$,$\kappa_4$ and $\kappa_5$. For this, we use Theorem \ref{thm:comparison_MFPT}. Let
\begin{equation*}\label{matrixAbraess}
	A= \begin{bmatrix}
		-1 & 0 & 0 & 0\\
		0 & -1 & -1 & 0
	\end{bmatrix}.
\end{equation*}
The matrix $A$ here defines a preorder that is not a partial order of $\X$. For $x \in \X$, consider infinitesimal transition rates $\Breve{\rate}_1(x),\Breve{\rate}_2(x), \Breve{\rate}_3(x), \Breve{\rate}_4(x)$ and $\Breve{\rate}_5(x)$ defined as for $\rate_1(x),\rate_2(x),\rate_3(x),\rate_4(x)$ and $\rate_5(x)$, but with $\Breve{\kappa}_i$ in place of $\kappa_i$ where $\Breve{\kappa}_i=\kappa_i$, for $i = 1,2,3,4$, and $\Breve{\kappa}_5 \neq \kappa_5$. Suppose that $\kappa_2 = \kappa_4$. Now, let us verify that the assumptions of Theorem \ref{thm:SumTheorem2} hold. Condition $(i)$ holds since $Av_1=(1,-1)^T$, $Av_2=(0,1)^T$, $Av_3=(1,-1)^T$, $Av_4=(0,1)^T$ and $Av_5=(0,0)^T$. Condition $(ii)$ of Theorem \ref{thm:SumTheorem2} also holds, as shown in the paragraph below.

\textbf{Verification of condition $(ii)$ of Theorem \ref{thm:SumTheorem2}.} Let $x \in \X$, and first consider $x\in \X$ and $y\in \partial_1(K_A+x) \cap \X $, where $\partial_1(K_A+x) \cap \X = \{ w \in \X \:|\: x_1 = w_1, x_2 + x_3 \geq w_2 + w_3, x_4 \le w_4\}$. Given that $Av_2=Av_4, Av_1=Av_3$, and $\inn{A_{1\bullet},v_1}=\inn{A_{1\bullet},v_3}=1$, we need to check that $\rate_1(x) + \rate_3(x) \le \Breve{\rate}_1(y) + \Breve{\rate}_3(y)$. Since $y \in \partial_1(K_A+x) \cap \X$, then $\rate_1(x) = \kappa_1 x_1 =\kappa_1 y_1 = \Breve{\kappa}_1 y_1 = \Breve{\rate}_1(y)$ and $\rate_3(x) = \kappa_3 x_1 =\kappa_3 y_1 = \Breve{\kappa}_3 y_1 = \Breve{\rate}_3(y)$, and so the desired inequality holds with equality. Secondly, consider $y \in \partial_2(K_A+x) \cap \X = \{ w \in \X \:|\: x_1 \geq w_1, x_2 + x_3 = w_2 + w_3, x_4\le w_4\}$. Given that $Av_2=Av_4, Av_1=Av_3$, and $\inn{A_{1\bullet},v_1}=\inn{A_{1\bullet},v_3}=-1$ and $\inn{A_{1\bullet},v_2}=\inn{A_{1\bullet},v_4}=1$, we need to check that $\rate_2(x) + \rate_4(x) \le \Breve{\rate}_2(y) + \Breve{\rate}_4(y)$ and $\rate_1(x) +\rate_3(x) \ge \Breve{\rate}_1(y) + \Breve{\rate}_3(y)$. For $x\in \X$ and $y\in \partial_2(K_A+x) \cap \X $, we have that $\rate_2(x) + \rate_4(x) =  \kappa_2 x_2 + \kappa_4 x_3 = \kappa_2 (x_2 + x_3) \le \kappa_2 (y_2 + y_3) = \Breve{\kappa}_2 (y_2 + y_3) = \Breve{\rate}_2(y) + \Breve{\rate}_4(y)$ and $\rate_1(x) = \kappa_1 x_1 \ge \kappa_1 y_1 = \Breve{\kappa}_1 y_1  = \Breve{\rate}_1(y)$, $\rate_3(x) = \kappa_3 x_1 \ge \kappa_3 y_1 = \Breve{\kappa}_3 y_1 = \Breve{\rate}_3(y)$.

Thus, all hypotheses of Theorem \ref{thm:SumTheorem2} hold, and so for every $\initialx,\initialxbreve \in \X$ where $\initialx \preccurlyeq_A \initialxbreve$ there there exists a probability space $(\Omega,\F,\PP)$ with two Markov chains $X = \{X(t), \: t \geq 0\}$ and $\Breve{X}=\{\Breve{X}(t), \: t \geq 0\}$ associated with $\rate$ and $\Breve{\rate}$, respectively, such that $X(0)=\initialx$, $\Breve{X}(0)=\initialxbreve$ and $\PP\left[X(t) \preccurlyeq_A \Breve{X}(t) \text{ for every } t \geq 0 \right]=1$. Let $\Gamma = \{(0,0,0,\mathrm{S_{tot}})\}$. This is an increasing set in $\X$ with respect to the relation $\preccurlyeq_A$. Let $T_{(0,0,0,\mathrm{S_{tot}})}$, respectively $\Breve{T}_{(0,0,0,\mathrm{S_{tot}})}$ be the first time that the Markov chain $X$, respectively $\Breve{X}$, reaches the set $\Gamma$. Then, by Theorem \ref{thm:comparison_MFPT}, if $X(0)=\Breve{X}(0)=(\mathrm{S_{tot}},0,0,0)$, we have that $\Breve{T}_{(0,0,0,\mathrm{S_{tot}})} \preccurlyeq_{st} T_{(0,0,0,\mathrm{S_{tot}})}$. By interchanging $\Breve{\rate}_5$ and $\kappa_5$, we can conclude that $\Breve{T}_{(0,0,0,\mathrm{S_{tot}})}$ and $T_{(0,0,0,\mathrm{S_{tot}})}$ are stochastically equivalent (equal in distribution). It follows that the mean first passage time from $(\mathrm{S_{tot}},0,0,0)$ to $(0,0,0,\mathrm{S_{tot}})$ is insensitive to $\kappa_5$ when $\kappa_2 = \kappa_4$. This is naively counter-intuitive: since the fifth reaction re-routes some samples to another state where the last reaction has the same rate constant as the final reaction without re-routing, it should take a longer expected time since re-routing also takes some time. However, in reality, the presence of the fifth reaction also fastens the rate to transition from $\mathrm{S_2}$, and this balances the time of re-routing. Most importantly, our theorem is able to capture this result without explicitly calculating the mean first passage time and allows us to reach the conclusion easily. We expect that in more complex situations, our method will be a valuable tool to establish monotonicity and insensitivity results.

Given that the Markov chain has one absorbing state per stoichiometric compatibility class, the stationary distribution for a given stoichiometric compatibility class is trivial, and hence so too are its monotonicity properties.

Theorem S.2 allows us to conclude further interesting properties for this network. Using two other $A$ matrices (see SI - Section S.3.2), we can conclude that adding reaction ${\large \textcircled{\small 5}}$ (changing from $\kappa_5 = 0$ to $\kappa_5 > 0$) causes the mean first passage time from $(\mathrm{S_{tot}},0,0,0)$ to $(0,0,0,\mathrm{S_{tot}})$ to increase if $\kappa_2 > \kappa_4$ or to decrease if $\kappa_2 < \kappa_4$. More explicitly, this shows that there can be opposing effects on the mean first passage time with different choices of $\kappa_2$ and $\kappa_4$ when reaction ${\large \textcircled{\small 5}}$ is added.

\end{example}
\begin{example}\textbf{Epigenetic regulation by chromation modifications}
\label{ex:HistoneModification}\\
\begin{figure}[t!]
	\centering
	\includegraphics[scale=0.44]{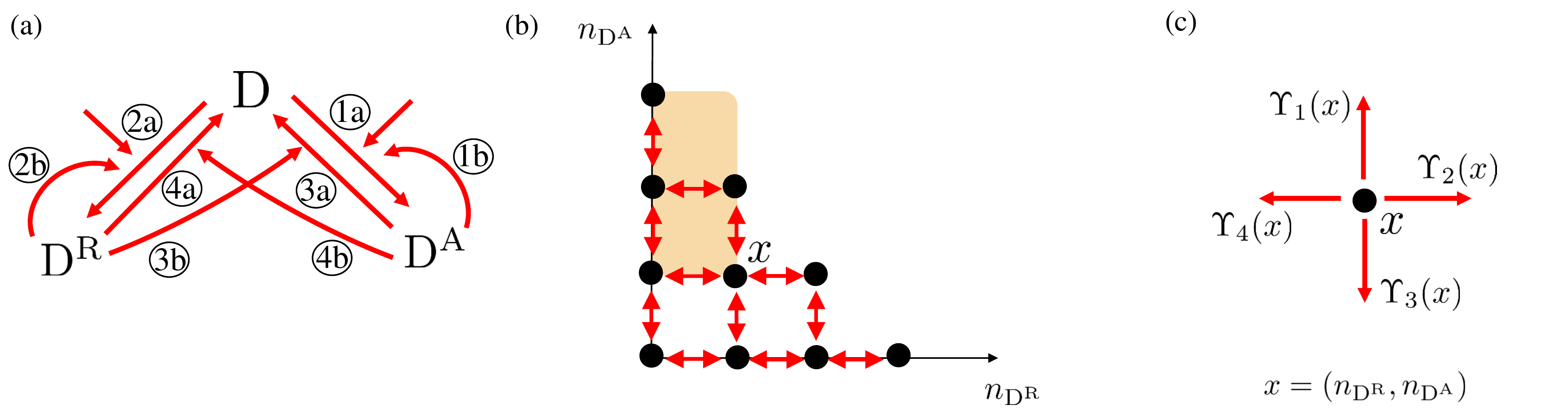}
	\caption{\small { \bf Histone modification circuit and corresponding Markov chain.} (a) Original chemical reaction system.  The numbers on the arrows correspond to the associated reactions. (b) Markov chain graph associated with the reduced chemical reaction system. Here, we consider $\Dtot=3$, we use black dots to represent the states and red double-ended arrows to represent transitions in both directions. We use orange to highlight the region $(K_A +x) \cap \X $, with $A$ defined in (\ref{matrixA2D}). (c) Direction of the possible transitions of the Markov chains, whose rates are given in equation (\ref{rates2D}).
	}
	\label{fig:2Dmodel}
\end{figure} 
Epigenetic regulation is the modification of the DNA structure, due to chromatin modifications, that determines if a gene is active or repressed. There are several chromatin modifications that can affect the DNA structure. Here, we will focus only on histone modifications. More precisely, we consider a ubiquitous model for a histone modification circuit (see Dodd et al. \cite{2007DoddCellPaper} and Bruno et al. \cite{BrunoDelVecchio}). The species considered are nucleosomes that are unmodified (D), modified with repressive modifications ($\mathrm{D^R}$), and modified with activating modifications ($\mathrm{D^A}$), and, in terms of molecular interactions, each histone modification autocatalyzes itself and promotes the erasure of the other one. The chemical reaction system considered is depicted in Fig. \ref{fig:2Dmodel}(a). The amount of each species is represented by $n_{\mathrm{D}}$, $n_{\mathrm{D^R}}$ and $n_{\mathrm{D^A}}$, respectively, and their sum is conserved, that is $n_{\mathrm{D}}+n_{\mathrm{D^R}}+n_{\mathrm{D^A}}=\Dtot$, with $\Dtot$ representing the total number of nucleosomes within the gene.

By fixing an integer $\Dtot>0$, we fix one stoichiometric compatibility class. The projected process $(X_1,X_2)(\cdot) = (n_{\mathrm{D^R}},n_{\mathrm{D^A}})$ is still a continuous-time Markov chain, and in this example we choose to apply our theory to this reduced system. This is the same as studying the reduced chemical reaction system defined as follows:
\begin{equation}\label{reacs2DBIS}
	\begin{aligned}
		&{\large \textcircled{\small $1$}}\;\ce{$\emptyset$ -> D^A },\;\;{\large \textcircled{\small $2$}}\;\ce{$\emptyset$ -> D^R },\;\;{\large \textcircled{\small $3$}}\;\ce{D^A -> $\emptyset$ },\;\;{\large \textcircled{\small $4$}}\;\ce{D^R -> $\emptyset$},\\
	\end{aligned}
\end{equation} 
with two species $\Ss =\{\mathrm{D^R}, \mathrm{D^A}\}$ and four reactions $\Rs = \{(v^-_{1},v^+_{1}), (v^-_{2},v^+_{2}),$ $(v^-_{3},v^+_{3}), (v^-_{4},v^+_{4})\}$, where $v^-_{1} =v^-_{2} =v^+_{3} =v^+_{4}=(0,0)^T$, $v^+_{2}= v^-_{4}= (1,0)^T$, $v^+_{1}= v^-_{3}= (0,1)^T$, and with associated propensity functions of non mass-action type defined as follows:
\begin{align}
	\notag &\prp_{(v^-_{1},v^+_{1})}(x)= (\Dtot -(x_1+x_2))\left(\kappa_{1a} + \kappa_{1b}x_2\right),\\
	\label{propfunct1} &\prp_{(v^-_{2},v^+_{2})}(x)= (\Dtot -(x_1+x_2))\left(\kappa_{2a} + \kappa_{2b} x_1\right),\\
	\notag &\prp_{(v^-_{3},v^+_{3})}(x)= x_2\left(\kappa_{3a} + x_1\kappa_{3b}\right),\;\;\prp_{(v^-_{4},v^+_{4})}(x) = x_1\mu\left(c\kappa_{3a} + x_2 \kappa_{3b}\right),
\end{align}
in which $\kappa_{1a}$, $\kappa_{1b}$, $\kappa_{3a}$, $\kappa_{3b}$, $\kappa_{2a}$, $\kappa_{2b}$, $\kappa_{4a} = \mu c\kappa_{3a}$, $\kappa_{4b} = \mu \kappa_{3b}$ are the rate constants that go with each of the reactions shown in Fig. \ref{fig:2Dmodel}(a), respectively. 

The state space for the Markov chain is $\mathcal{X}= \{(x_1,x_2) \in \Z_+^2 \:|\: x_1 + x_2 \leq \Dtot \}$. Given a generic state $x=(x_1,x_2)\in \X$, the potential transitions of the Markov chain are in four possible directions $v_j=v^+_j-v^-_j$, $j=1, 2, 3, 4$, that can be written as $v_{1}= (0,1)^T, v_2=(1,0)^T, v_3=(0,-1)^T$ and $v_4=(-1,0)^T$, with associated infinitesimal transition rates
\begin{equation}\label{rates2D}
	\begin{aligned}
		&\rate_{1}(x)= \prp_{(v^-_{1},v^+_{1})}(x),\;\;\rate_{2}(x)= \prp_{(v^-_{2},v^+_{2})}(x),\;\;\rate_{3}(x)= \prp_{(v^-_{3},v^+_{3})}(x),\;\;\rate_{4}(x)= \prp_{(v^-_{4},v^+_{4})}(x).
	\end{aligned}
\end{equation}
We are interested in determining how the asymmetry of the system, represented by the parameter $\mu$ affects the stochastic behavior of the system. In particular, we will focus on studying the stationary distribution and the \textit{time to memory loss} of the active and repressed state, defined as the mean first passage time to reach the fully repressed state ($r=(n_{\mathrm{D^R}},n_{\mathrm{D^A}})=(\Dtot,0)$), starting from the fully active state ($a=(n_{\mathrm{D^R}},n_{\mathrm{D^A}})=(0,\Dtot)$), and vice versa (i.e., $h_{a,r}= \E_{a}[T_{r}]$ and $h_{r,a}= \E_{r}[T_{a}]$). To this end, we first verify that we can apply Theorem \ref{thm:InnerProductTheorem}.

Let
\begin{equation}\label{matrixA2D}
	A= \begin{bmatrix}
		-1 & 0\\
		0 & 1
	\end{bmatrix}.
\end{equation}
For $x\in \X$, $K_A +x = \{ w \in \R^2 \:|\: x \preccurlyeq_A w \}$ and $ (K_A +x)\cap \X = \{ w \in \X \:|\: x \preccurlyeq_A w \}$. See Fig. \ref{fig:2Dmodel}(b) for an example of $\X$ and $(K_A +x) \cap \X $ for $\Dtot=3$. We introduce infinitesimal transition rates $\Breve{\rate}_{1}(x),\Breve{\rate}_{2}(x), \Breve{\rate}_{3}(x)$ and $\Breve{\rate}_{4}(x)$ defined as for $\rate_{1}(x),\rate_{2}(x),\rate_{3}(x)$ and $\rate_{4}(x)$, with all the parameters having the same values except that $\mu$ is replaced by $\Breve{\mu}$, where $\Breve{\mu}\ge \mu$. Since $Av_1=(0,1)^T$, $Av_2=(-1,0)^T$, $Av_3=(0,-1)^T$ and $Av_4=(1,0)^T$, we have that condition $(i)$ of Theorem \ref{thm:InnerProductTheorem} holds. Condition $(ii)$ also holds, as shown in the paragraph below.

\textbf{Verification of condition $(ii)$ of Theorem \ref{thm:InnerProductTheorem}.} Consider $x\in \X$ and $y \in \partial_1(K_A+x) \cap \X$, where $\partial_1(K_A+x) \cap \X=\{ w \in \X \:|\: x_1=w_1, x_2\le w_2\}$. Since $\inn{A_{1\bullet},v_{4}}=1$ and $\inn{A_{1\bullet},v_{2}}=-1$, we must check that $\rate_{2}(x)\ge\Breve{\rate}_{2}(y)$ and $\rate_{4}(x)\le\Breve{\rate}_{4}(y)$. Since $y \in \partial_1(K_A+x) \cap \X$ implies $x_1=y_1$ and $x_2\le y_2$, we have $\rate_{2}(x) = (\Dtot -(x_1+x_2))\left(\kappa_{2a} + \kappa_{2b} x_1\right) \ge (\Dtot -(y_1+y_2))\left(\kappa_{2a} + \kappa_{2b} y_1\right) = \Breve{\rate}_{2}(y)$ and $\rate_{4}(x) =x_1\mu\left(c\kappa_{3a} + x_2 \kappa_{3b}\right) \le y_1\mu\left(c\kappa_{3a} + y_2 \kappa_{3b}\right) \le y_1\Breve{\mu}\left(c\kappa_{3a} + y_2 \kappa_{3b}\right) = \Breve{\rate}_{4}(y)$, and so both inequalities hold. Similarly, for $x\in \X$ and $y \in \partial_2(K_A+x) \cap \X=\{ w \in \X \:|\: x_1\ge w_1, x_2=w_2 \}$, since $\inn{A_{2\bullet},v_{1}}=1$ and $\inn{A_{2\bullet},v_{3}}=-1$, we need to check that $\rate_{1}(x)\le\Breve{\rate}_{1}(y)$ and $\rate_{3}(x)\ge\Breve{\rate}_{3}(y)$. Indeed, $\rate_{1}(x) = (\Dtot -(x_1+x_2))\left(\kappa_{1a} + \kappa_{1b}x_2\right) \le (\Dtot -(y_1+y_2))\left(\kappa_{1a} + \kappa_{1b} y_2\right) = \Breve{\rate}_{1}(y)$ and $\rate_{3}(x) = x_2\left(\kappa_{3a} + x_1\kappa_{3b}\right) \ge y_2\left(\kappa_{3a} + y_1\kappa_{3b}\right) = \Breve{\rate}_{3}(y)$. 

Since all of the hypotheses of Theorem \ref{thm:InnerProductTheorem} hold, for each pair $\initialx,\initialxbreve \in \X$ satisfying $\initialx \preccurlyeq_A \initialxbreve$, there exists a probability space $(\Omega,\F,\PP)$ with two Markov chains $X = \{X(t), \: t \geq 0\}$ and $\Breve{X}=\{\Breve{X}(t), \: t \geq 0\}$ associated with $\rate$ and $\Breve{\rate}$, respectively, such that $X(0)=\initialx$, $\Breve{X}(0)=\initialxbreve$ and $\PP\left[X(t) \preccurlyeq_A \Breve{X}(t) \text{ for every } t \geq 0 \right]=1$.

We can also apply Theorem \ref{comparStationaryDistributionMAIN}. The Markov chains $X$ and $\Breve{X}$ are irreducible and, having only finitely many states, are positive recurrent. Based on the order $\preccurlyeq_A$ we introduced, the fully active state $a=(0,\Dtot)$ is maximal in $\X$ and the fully repressed state $r=(\Dtot,0)$ is minimal in $\X$. Then, by Theorem \ref{comparStationaryDistributionMAIN}, we can conclude that  $\pi_a \leq  \Breve{\pi}_a$ and $\pi_r \geq  \Breve{\pi}_r$. This implies that increasing $\mu$ increases the probability of the system in steady-state to be in the active state $a$ to the detriment of the repressed state $r$ (and vice versa for decreasing $\mu$). We can also apply Theorem \ref{thm:comparison_MFPT}. Since $\{a\}$ is increasing and $\{r\}$ is decreasing, then by Theorem \ref{thm:comparison_MFPT}, $\Breve{h}_{r,a}=\E_{r}[\Breve{T}_{a}]\le \E_{r}[T_{a}]=h_{r,a}$ and $h_{a,r}=\E_{a}[T_{r}]\le \E_{a}[\Breve{T}_{r}]=\Breve{h}_{a,r}$. Since the only difference between the two systems was that $\mu \le\Breve{\mu}$, these results imply that the time to memory loss of the active state increases for higher values of $\mu$, while the time to memory loss of the repressed state decreases for higher values of $\mu$.

\end{example}

\begin{example}\textbf{Epigenetic regulation by chromatin modifications with positive TF-enabled autoregulation}
\label{ex:HistoneModificationAndProtein}\\
\begin{figure}[t!]
	\centering
	\includegraphics[scale=0.44]{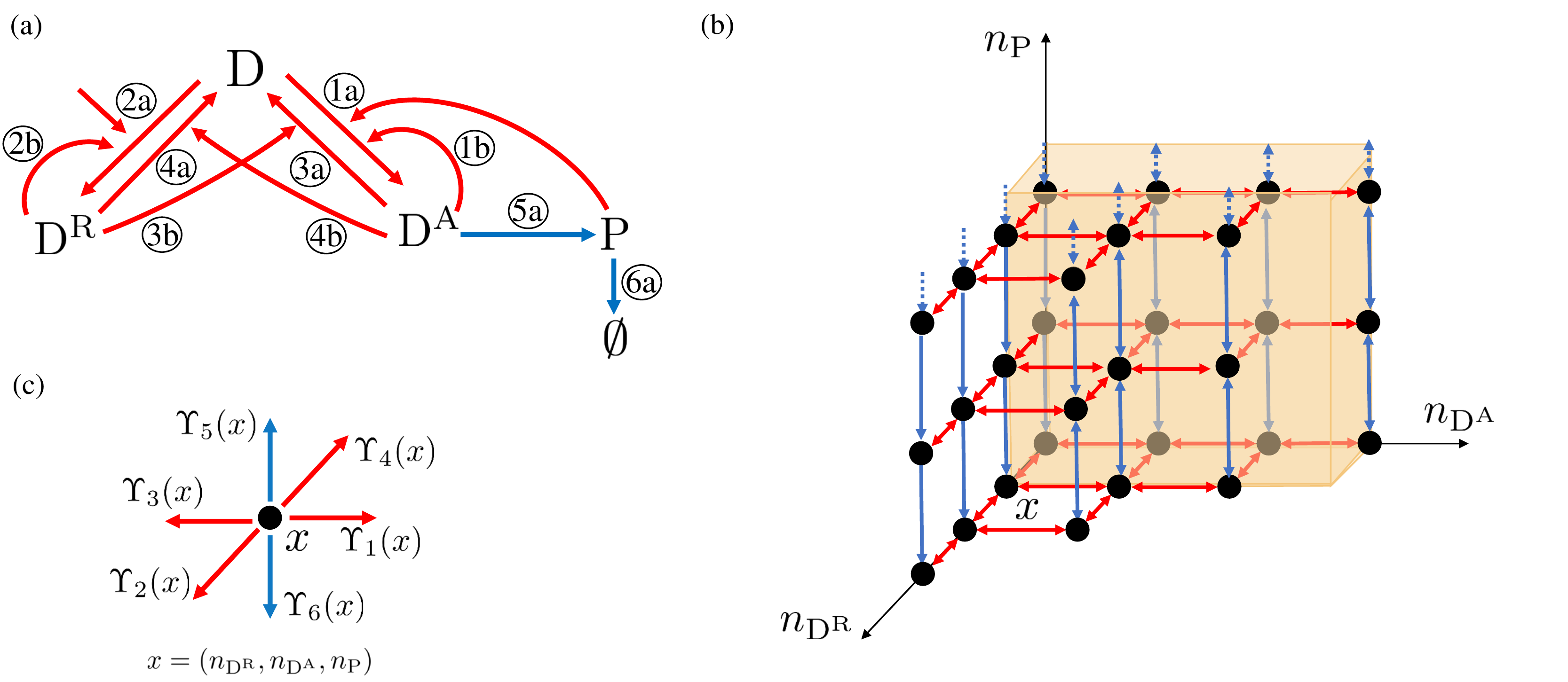}
	\caption{\small { \bf Histone modification circuit with positive TF-enabled autoregulation and corresponding Markov chain.} (a) Original chemical reaction system.  The numbers on the arrows correspond to the associated reactions. (b) Markov chain graph. Here, we consider $\Dtot=3$, we use black dots to represent the states and red double-ended arrows to represent transitions in both directions associated with the reactions represented by the red arrows in (a). Similarly we use blue double-ended (single-ended) arrows to represent transitions in both directions (in one direction) associated with the reactions represented by the blue arrows in (a). We use blue dotted lines to show that, in the vertical direction, the Markov chain has countably infinitely many states, connected by transitions in both directions. Finally, we use orange to highlight the region $K_A+x$ intersected with the state space $\X$, with $A$ defined in (\ref{matrixA2D+1}). (c) Direction of the possible transitions of the Markov chain starting from a state $x$, whose rates are defined in equation (\ref{rates2D+1}).
	}
	\label{fig:2Dmodel+P}
\end{figure}
Now, we consider the histone modification circuit considered in the previous example with an additional positive autoregulation loop. For this, we assume that a protein expressed by the gene of interest recruits writers for the activating histone modifications. Consequently, we introduce the gene product P as an additional species for our system and add the following reactions to the ones shown in Fig. \ref{fig:2Dmodel}(a):
\begin{equation}\label{autoreg}
	\begin{aligned}
		&{\large \textcircled{\small 5a}}\;\ce{D^A  ->[ ] D^A + P },\;\;{\large \textcircled{\small 6a}}\; \ce{P ->[ ] \emptyset }.\\
	\end{aligned}
\end{equation}
Furthermore, given the P-enabled autoregulation loop (Fig. \ref{fig:2Dmodel+P}(a)), let us consider the rate constant that goes with  ${\large \textcircled{\small 1a}}$ in Fig. \ref{fig:2Dmodel}(a) as
$\kappa_{1a}=\kappa_{1a}^0+\kappa_{1a}^1g(n_{\mathrm{P}})$, with
$\kappa_{1a}^0$ and $\kappa_{1a}^1$ representing the rate constants that go with the $\mathrm{D^A}$ basal \textit{de-novo} establishment process and with the $\mathrm{D^A}$ \textit{de-novo} establishment process enhanced by $n_{\mathrm{P}}$, respectively, and $g(n_{\mathrm{P}})$ representing a non-negative, bounded, monotonically increasing function of $n_{\mathrm{P}}$ (see Bruno et al. \cite{BrunoDelVecchio}, Section 3.4). 

Here, we are interested in determining how the reaction rate constant $\kappa_{5a}$ affects the reactivation time of the gene. As before, we have the conservation law $n_{\mathrm{D}}+n_{\mathrm{D^R}}+n_{\mathrm{D^A}}=\Dtot$, with $\Dtot$ representing the total number of nucleosomes within the gene, and by fixing $\Dtot>0$, we fix one stoichiometric compatibility class and the projected process $(X_1,X_2,X_3)(\cdot) = (n_{\mathrm{D^R}},n_{\mathrm{D^A}},n_{\mathrm{P}})$ is a continuous-time Markov chain. This is the same as studying the reduced chemical reaction system:
\begin{equation}\label{reacs2DBIS+1}
	\begin{aligned}
		&{\large \textcircled{\small $1$}}\;\ce{$\emptyset$ ->[] D^A },\;\;{\large \textcircled{\small $2$}}\;\ce{$\emptyset$ ->[] D^R },\;\;{\large \textcircled{\small $3$}}\;\ce{D^A ->[] $\emptyset$ },\\
		&{\large \textcircled{\small $4$}}\;\ce{D^R ->[] $\emptyset$},\;\;{\large \textcircled{\small $5$}}\;\ce{D^A  ->[] D^A + P },\;\;{\large \textcircled{\small $6$}}\;\ce{P ->[] \emptyset },\\
	\end{aligned}
\end{equation} 
with set of species $\Ss =\{\mathrm{D^R}, \mathrm{D^A}, \mathrm{P}\}$, set of reactions $\Rs = \{(v^-_{1},v^+_{1})$, $(v^-_{2},v^+_{2})$, $(v^-_{3},v^+_{3})$, $(v^-_{4},v^+_{4})$, $(v^-_{5},v^+_{5})$, $(v^-_{6},v^+_{6})\}$, where $v^-_{1} =v^-_{2} =v^+_{3} =v^+_{4}=v^+_{6} =(0,0,0)^T$, $v^+_{2}= v^-_{4}=(1,0,0)^T$, $v^+_{1}= v^-_{3}= v^-_{5}=(0,1,0)^T$, $v^+_{5}=(0,1,1)^T$, $v^-_{6} =(0,0,1)^T$, and with associated propensity functions of non mass-action type defined as follows:
\begin{equation*}
	\begin{aligned}
		&\prp_{(v^-_{1},v^+_{1})}(x)= (\Dtot -(x_1+x_2))\left(\kappa_{1a}^0+\kappa_{1a}^1g(x_3) + \kappa_{1b}x_2\right),\\
		&\prp_{(v^-_{2},v^+_{2})}(x)= (\Dtot -(x_1+x_2))\left(\kappa_{2a} + \kappa_{2b}x_1\right),\;\;\prp_{(v^-_{3},v^+_{3})}(x)= x_2\left(\kappa_{3a} + x_1 \kappa_{3b}\right),\\ &\prp_{(v^-_{4},v^+_{4})}(x) = x_1\mu\left(c\kappa_{3a} + x_2\kappa_{3b}\right),\;\;\prp_{(v^-_{5},v^+_{5})}(x) = \kappa_{5a} x_2,\;\;\prp_{(v^-_{6},v^+_{6})}(x) = \kappa_{6a} x_3,
	\end{aligned}
\end{equation*}
in which $\kappa_{5a}$ and $\kappa_{6a}$ are the rate constants that go with reactions ${\large \textcircled{\small 5a}}$ and ${\large \textcircled{\small 6a}}$ in (\ref{autoreg}), respectively, and all the other rate constants are defined as for (\ref{propfunct1}).

The state space for the Markov chain is $\mathcal{X}= \{(x_1,x_2,x_3) \in \Z_+^3 \:|\: x_1 + x_2 \leq \Dtot \}$. Given a generic state $x=(x_1,x_2,x_3)$, the transitions of the Markov chain are in six possible directions  $v_j=v^+_j-v^-_j$, $j\in\{1, ..., 6\}$, that can be written as $v_1=(0,1,0)^T$, $v_2=(1,0,0)^T$, $v_3=(0,-1,0)^T$, $v_4=(-1,0,0)^T$, $v_5=(0,0,1)^T$, $v_6=(0,0,-1)^T$, with associated infinitesimal transition rates:
\begin{equation}\label{rates2D+1}
	\begin{aligned}
		&\rate_{1}(x)= \prp_{(v^-_{1},v^+_{1})}(x),\;\;\rate_{2}(x)= \prp_{(v^-_{2},v^+_{2})}(x),\;\;\rate_{3}(x)= \prp_{(v^-_{3},v^+_{3})}(x),\\
		&\rate_{4}(x)= \prp_{(v^-_{4},v^+_{4})}(x),\;\;\rate_{5}(x)= \prp_{(v^-_{5},v^+_{5})}(x),\;\;\rate_{6}(x)= \prp_{(v^-_{6},v^+_{6})}(x).
	\end{aligned}
\end{equation}
As mentioned before, we are interested in determining how the protein production rate $\kappa_{5a}$ affects the reactivation time of the gene, defined as $h_{r,\Theta}= \E_{r}[T_{\Theta}]$, where $r=(\Dtot,0,0)$ and $\Theta=\{w\in \X| w=(0,\Dtot,i), i\in\Z_{+}\}$ corresponds to the set of states characterized by the fully active state $n_{\mathrm{D^A}}=\Dtot$. We first check that the assumptions of Theorem \ref{thm:InnerProductTheorem} hold.
Let
\begin{equation}\label{matrixA2D+1}
	A= \begin{bmatrix}
		-1 & 0 & 0\\
		0 & 1 & 0\\
		0 & 0 & 1
	\end{bmatrix}.
\end{equation}
For $x\in \mathcal{X}$, $x \preccurlyeq_A y$ and the set $K_A +x = \{ y \in \R^3 \:|\: x \preccurlyeq_A y \}$. For our example, the region $(K_A +x) \cap \X$ is depicted in orange in Fig. \ref{fig:2Dmodel+P}(b).
We introduce infinitesimal transition rates $\Breve{\rate}_{1}(x),\Breve{\rate}_{2}(x), \Breve{\rate}_{3}(x), \Breve{\rate}_{4}(x), \Breve{\rate}_{5}(x)$ and $\Breve{\rate}_{6}(x)$ defined as for $\rate_{1}(x),\rate_{2}(x),\rate_{3}(x),\rate_{4}(x),\rate_{5}(x)$ and $\rate_{6}(x)$, with all the parameters having the same values except that $\kappa_{5a}$ is replaced by $\Breve{\kappa}_{5a}> \kappa_{5a}$. Condition $(i)$ of Theorem \ref{thm:InnerProductTheorem} holds since $Av_{1}=(0,1,0)^T$, $Av_{2}=(-1,0,0)^T$, $Av_{3}=(0,-1,0)^T$, $Av_4=(1,0,0)^T$, $Av_5=(0,0,1)^T$, $Av_6=(0,0,-1)^T$. Condition $(ii)$ also holds, as shown in the paragraph below.

\textbf{Verification of condition $(ii)$ of Theorem \ref{thm:InnerProductTheorem}.} First consider $x\in \X$ and $y \in \partial_1(K_A+x) \cap \X=\{ w \in \X \:|\: x_1=w_1, x_2\le w_2, x_3\le w_3\}$. Since $\inn{A_{1\bullet},v_{4}}=1$ and $\inn{A_{1\bullet},v_{2}}=-1$, we need to check that $\rate_{4}(x)\le\Breve{\rate}_{4}(y)$ and $\rate_{2}(x)\ge\Breve{\rate}_{2}(y)$. Since $x_1=y_1, x_2\le y_2, x_3\le y_3$, we have that $\rate_{4}(x)=x_1\mu\left(c\kappa_{3a} + x_2\kappa_{3b}\right) \le y_1\mu\left(c\kappa_{3a} + y_2\kappa_{3b}\right)  = \Breve{\rate}_{4}(y)$ and $\rate_{2}(x)= (\Dtot -(x_1+x_2))\left(\kappa_{2a} + \kappa_{2b}x_1\right) \ge (\Dtot -(y_1+y_2))\left(\kappa_{2a} + \kappa_{2b}y_1\right) = \Breve{\rate}_{2}(y)$.
Secondly, consider $x\in \X$ and $y \in \partial_2(K_A+x) \cap \X=\{ w \in \X \:|\: x_1\ge w_1, x_2=w_2,  x_3\le w_3\}$. Since 
$\inn{A_{2\bullet},v_{1}}=1$ and $\inn{A_{2\bullet},v_{3}}=-1$, we need to check that $\rate_{1}(x)\le\Breve{\rate}_{1}(y)$ and $\rate_{3}(x)\ge\Breve{\rate}_{3}(y)$. Since $x_1\ge y_1, x_2=y_2,  x_3\le y_3$, we have $\rate_{1}(x)= (\Dtot -(x_1+x_2))\left(\kappa_{1a}^0+\kappa_{1a}^1g(x_3) + \kappa_{1b}x_2\right) \le (\Dtot -(y_1+y_2))\left(\kappa_{1a}^0+\kappa_{1a}^1g(y_3) + \kappa_{1b} y_2\right) =  \Breve{\rate}_{1}(y)$ and $\rate_{3}(x)=x_2\left(\kappa_{3a} + x_1 \kappa_{3b}\right) \ge y_2\left(\kappa_{3a} + y_1 \kappa_{3b}\right) = \Breve{\rate}_{3}(y)$.
Finally, consider $x\in \X$ and $y \in \partial_3(K_A+x) \cap \X=\{ w \in \X \:|\: x_1 \ge w_1, x_2\le w_2,  x_3=w_3\}$. Since $\inn{A_{3\bullet},v_{5}}=1$ and $\inn{A_{3\bullet},v_{6}}=-1$, we must check that $\rate_{5}(x)\le\Breve{\rate}_{5}(y)$ and $\rate_{6}(x)\ge\Breve{\rate}_{6}(y)$. Since $x_1 \ge y_1, x_2\le y_2,  x_3=y_3$, we obtain $\rate_{5}(x)=\kappa_{5a} x_2 \le \kappa_{5a} y_2 \le \Breve{\kappa}_{5a} y_2 =  \Breve{\rate}_{5}(y)$ and $\rate_{6}(x)=\kappa_{6a} x_3 =\kappa_{6a} y_3 =\Breve{\rate}_{6}(y)$.

Since all the hypotheses of Theorem \ref{thm:InnerProductTheorem} hold, for each $\initialx,\initialxbreve \in \X$ satisfying $\initialx \preccurlyeq_A \initialxbreve$, there exists a probability space $(\Omega,\F,\PP)$ with two Markov chains $X = \{X(t), \: t \geq 0\}$ and $\Breve{X}=\{\Breve{X}(t), \: t \geq 0\}$ associated with $\rate$ and $\Breve{\rate}$, respectively, such that $X(0)=\initialx$, $\Breve{X}(0)=\initialxbreve$ and $\PP\left[X(t) \preccurlyeq_A \Breve{X}(t) \text{ for every } t \geq 0 \right]=1$.

Furthermore, since the hypotheses of Theorem \ref{thm:InnerProductTheorem} hold, we can also apply Theorem \ref{thm:comparison_MFPT}. Specifically, for $r=(\Dtot,0,0)$ and $\Theta=\{y\in \X| y=(0,\Dtot,i), i\in\Z_{+}\}$, since $\Theta$ is an increasing set in $\X$ with respect to the relation $\preccurlyeq_A$, then $h_{r,\Theta}\ge \Breve{h}_{r,\Theta}$. This implies that, assuming that the only difference between the two systems is in the value of the protein production rate parameter, $\kappa_{5a}$, higher protein production rates reduce the mean reaction time for the gene.

\end{example}

	\section{Proofs of the Main Results}
	\label{sec:proofmainresults}
	
	\subsection{Proof of Theorem \ref{thm:MainResult}.}
	\label{PROOFtheoremMonAppendix_I}
	
	Consider a non-empty set $\X \subseteq \Z_+^d$, a collection of distinct vectors $v_1,\ldots,v_n$ in $\Z^d \setminus \{0\}$ and two collections of non-negative functions on $\X$, $\rate=(\rate_1, \dots,\rate_n)$ and  $\Breve{\rate}= (\Breve{\rate}_1, \dots,\Breve{\rate}_n)$, such that \eqref{eq:LambdaNotOutofSpace} holds. Let $Q=(Q_{x,y})_{x,y \in \X}$ and $\Breve{Q}=(\Breve{Q}_{x,y})_{x,y \in \X}$ denote the infinitesimal generators for the continuous-time Markov chains associated with $\rate$ and $\Breve{\rate}$, respectively. In the following, let $A \in \R^{m \times d}$ be a matrix with non-zero rows and consider the relation $\preccurlyeq_A$ as defined in Definition \ref{def:preorderA}.
	
	For the proof of Theorem \ref{thm:MainResult}, we first assume that 
	\begin{equation}
		\label{eq:uniformization_condition_intensity}
		\sup_{x \in \X} \rate_j(x) < \infty \quad \text{ and } \quad \sup_{x \in \X} \Breve{\rate}_j(x) < \infty  \qquad \text{ for every } 1 \leq j \leq n.
	\end{equation}
	This restriction will be relaxed later. Then, we define a constant $\lambda > 0$ and a pair of functions $\Phi_{\lambda}$ and $\Breve{\Phi}_{\lambda}$, which will be key to our construction of the coupled processes $X$ and $\Breve{X}$. Let $\lambda > 0$ such that:
	\begin{equation}
		\label{eq:LambdaCondition}
		\lambda > n\max\left\{\sup_{x \in \X} \sum_{j=1}^n\rate_j(x), \sup_{x \in \X} \sum_{j=1}^n\Breve{\rate}_{j}(x) \right\}.
	\end{equation}
	Note that both $\frac{\rate_j(x)}{\lambda}$ and $\frac{\Breve{\rate}_j(x)}{\lambda}$ are less than $\frac{1}{n}$ for every $x \in \X$ and $1 \leq j \leq n$. For $x \in \X$, consider the sets
	\begin{equation}
		\label{eq:IntervalForPhi}
		I_j(x) := \left[\frac{j-1}{n},\frac{j-1}{n} + \frac{\rate_j(x)}{\lambda}\right), \qquad 1 \leq j \leq n.
	\end{equation}
	If $\rate_j(x) = 0$, then $I_j(x)$ is the empty set. On the other hand, if $\rate_j(x) >0$, then $I_j(x)$ is an interval that is a strict subset of $[\frac{j-1}{n},\frac{j}{n})$. Define the function $\Phi_{\lambda}(\cdot,\cdot): \X \times [0,1] \longrightarrow \X$ by 
	\begin{equation}
		\label{eq:DefinitionPhiLambda}
		\Phi_{\lambda}(x,u) := x + \sum_{j=1}^n v_j\one_{I_j(x)}(u), \qquad x \in \X, \; u \in [0,1].
	\end{equation}
	For $x \in \X$, the sets $I_1(x),\ldots,I_n(x)$ are mutually disjoint and so for any $u \in [0,1]$ either $\Phi_{\lambda}(x,u) =x$ or $\Phi_{\lambda}(x,u) =x + v_j$ for some $1 \leq j \leq n$. In the second case, this will happen if and only if $u \in I_j(x)$ for the corresponding index $j$. The latter condition implies that $I_j(x) \neq \emptyset$, hence by \eqref{eq:IntervalForPhi}, $\rate_{j}(x) > 0$ and by \eqref{eq:LambdaNotOutofSpace}, $x + v_{j} \in \X$.
	
	This shows that $\Phi_{\lambda}(\cdot,\cdot)$ is well-defined as an $\X$-valued function. We define intervals $\Breve{I}_j(x), \: 1 \leq j \leq n, \: x \in \X$ and a function $\Breve{\Phi}_{\lambda}:\X \times [0,1] \longrightarrow \X$ in an analogous manner to that above, where $\Breve{\Phi}_{\lambda}$ is defined as in \eqref{eq:DefinitionPhiLambda}, but with the intervals $I_j(x)$ replaced by $\Breve{I}_j(x)$, where these are defined as in \eqref{eq:IntervalForPhi}, but with $\rate_j(x)$ replaced by $\Breve{\rate}_j(x)$.
	
	\begin{lemma}
		\label{lem:PhiLambda_ineq}
		Suppose that $x,y \in \X$ are such that $x \preccurlyeq_A y$ and the following hold:
		\begin{equation}
			\label{eq:PropensityFunctions_ineq_I}
			\Breve{\rate}_j(y) \leq \rate_j(x), \quad \text{for each } 1 \leq j \leq n \text{ such that } y +v_j \in \X \setminus (K_A +x),
		\end{equation}
		and
		\begin{equation}
			\label{eq:PropensityFunctions_ineq_II}
			\Breve{\rate}_j(y)\geq \rate_j(x), \quad \text{for each } 1 \leq j \leq n \text{ such that } x +v_j \in \X \text{ and } y \notin K_A + x+ v_j.
		\end{equation}
		Then, for each $u \in [0,1]$,
		\begin{equation}
			\label{eq:ClaimPhi_LambdaMAIN}
			\Phi_{\lambda}(x,u) \preccurlyeq_A  \Breve{\Phi}_{\lambda}(y,u).
		\end{equation}
	\end{lemma}
	
	\begin{proof}
		First, we note that $\Phi_{\lambda}, \Breve{\Phi}_{\lambda}$ have the following property: for every $u \in [0,1]$ and $1 \leq j \leq n$,
		\begin{equation}
			\label{eq:PhiLambdaCondition1}
			\text{ if } \; \Phi_{\lambda}(x,u) = x +v_j, \: \text{ then } \; \Breve{\Phi}_{\lambda}(y,u) \in \{y,y +v_j\},
		\end{equation}
		since $I_j(x), \Breve{I}_j(y) \subseteq [\frac{j-1}{n},\frac{j}{n})$. Similarly,
		\begin{equation}
			\label{eq:PhiLambdaCondition2}
			\text{ if } \; \Breve{\Phi}_{\lambda}(y,u) = y +v_j, \: \text{ then } \; \Phi_{\lambda}(x,u) \in \{x,x +v_j\}.
		\end{equation}
		Furthermore, if $\Breve{\rate}_j(y) \geq \rate_j(x)$, then
		\begin{equation}
			\label{eq:PhiLambdaImplicationI}
			\Phi_{\lambda}(x,u) = x +v_j \: \text{ implies that } \; \Breve{\Phi}_{\lambda}(y,u) =y +v_j,
		\end{equation}
		since under this condition, $I_j(x) \subseteq \Breve{I}_j(y)$. Similarly, if $\Breve{\rate}_j(y) \leq \rate_j(x)$, then
		\begin{equation}
			\label{eq:PhiLambdaImplicationII}
			\Breve{\Phi}_{\lambda}(y,u) = y +v_j \: \text{ implies that } \; \Phi_{\lambda}(x,u) =x +v_j.
		\end{equation}
		
		Now, to prove \eqref{eq:ClaimPhi_LambdaMAIN}, fix $u \in [0,1]$. We consider two cases. 
		
		\medskip
		\textbf{Case 1:} $\Breve{\Phi}_{\lambda}(y,u) = y +v_j$ for some $1 \leq j \leq n$.
		
		\medskip
		Fix such an index $j$. Then, by \eqref{eq:PhiLambdaCondition2}, either $\Phi_{\lambda}(x,u) = x + v_j$ or $\Phi_{\lambda}(x,u) = x$. 
		\begin{enumerate}
			\item[a)]
			If $\Phi_{\lambda}(x,u) = x + v_j$, then, by \eqref{eq:translations_precurly}, $x + v_j \preccurlyeq_A y +v_j$ and therefore $\Phi_{\lambda}(x,u) \preccurlyeq_A  \Breve{\Phi}_{\lambda}(y,u)$. 
			\item[b)]
			If $\Phi_{\lambda}(x,u) = x$, then $y + v_j \in K_A +x$. To see this, we note that $y + v_j \in \X$ by \eqref{eq:LambdaNotOutofSpace} and since $\Breve{\rate}_j(y) > 0$ because $\Breve{I}_j(y) \neq \emptyset$. Then, if $y + v_j \notin K_A + x$, by \eqref{eq:PropensityFunctions_ineq_I}, we would have $\Breve{\rate}_j(y) \leq \rate_j(x)$, which would imply that $\Phi_{\lambda}(x,u) = x +v_j$ by \eqref{eq:PhiLambdaImplicationII}. But this contradicts the assumption that $\Phi_{\lambda}(x,u)=x$. Thus, $y + v_j \in K_A + x$ and so $\Phi_{\lambda}(x,u)=x \preccurlyeq_A y+v_j=\Breve{\Phi}_{\lambda}(y,u)$.
		\end{enumerate}
		
		\medskip
		\textbf{Case 2:} $\Breve{\Phi}_{\lambda}(y,u) = y$. Again, we consider two subcases.
		\begin{enumerate}
			\item[a)]
			If $\Phi_{\lambda}(x,u)=x$, then \eqref{eq:ClaimPhi_LambdaMAIN} holds, since $x \preccurlyeq_A y$ by assumption.
			\item[b)]
			If $\Phi_{\lambda}(x,u) = x +v_j$ for some $1 \leq j \leq n$, then $y \in K_A + x + v_j$ for the corresponding value of $j$. To see this, fix the value of $j$ for which $\Phi_{\lambda}(x,u)=x+v_j$ and notice that $x +v_j \in \X$ by \eqref{eq:LambdaNotOutofSpace} and since $\rate_j(x) >0$. If $y \notin K_A +x+v_j$, then by \eqref{eq:PropensityFunctions_ineq_II} we would have $\rate_j(x) \leq \Breve{\rate}_j(y)$, which would imply that $\Breve{\Phi}_{\lambda}(y,u) = y +v_j$. This contradicts the assumption that $\Breve{\Phi}_{\lambda}(y,u)=y$. Thus, we must have $y \in K_A +x +v_j$ and then $\Phi_{\lambda}(x,u)=x+v_j \preccurlyeq_A y=\Breve{\Phi}_{\lambda}(y,u)$.
		\end{enumerate}
		\vskip -0.2in    
	\end{proof}
	
	Now that all these preliminaries have been established under assumption \eqref{eq:uniformization_condition_intensity}, we proceed with the main part of the proof of Theorem \ref{thm:MainResult} with this assumption. For this proof, we assume that all of the conditions of Theorem \ref{thm:MainResult} hold and in addition that condition \eqref{eq:uniformization_condition_intensity} holds. The latter ensures that the pair of continuous-time Markov chains with infinitesimal generators $Q$ and $\Breve{Q}$ are \textit{uniformizable} (see Chapter 2 in Keilson \cite{KeilsonRarityExponentiality}). With $\lambda > 0$ as in \eqref{eq:LambdaCondition}, the (possibly infinite) matrices \footnote{These ``matrices" may have countably many rows and columns, in which case they could be considered as operators on $\l^{\infty}$. For convenience, we still call them matrices here.} $P_{\lambda}(Q) := \frac{1}{\lambda}Q +I$ and $P_{\lambda}(\Breve{Q}) := \frac{1}{\lambda}\Breve{Q} +I$ are stochastic \footnote{Stochastic here means that all entries take values in $[0,1]$ and all row sums equal one.}, where $I =(I_{x,y})_{x,y \in \X}$ is the identity matrix. Indeed, for $x \in \X$, $(P_{\lambda}(Q))_{x,x} = \frac{Q_{x,x}}{\lambda} + 1 = 1 - \frac{|Q_{x,x}|}{\lambda} \in [1-\frac{1}{n},1]$, for $y \neq x$, $(P_{\lambda}(Q))_{x,y} = \frac{Q_{x,y}}{\lambda} \in [0, \frac{1}{n}]$ and $\sum_{y \in \X} (P_{\lambda}(Q))_{x,y} = \sum_{y \in \X} \frac{1}{\lambda}Q_{x,y} + 1 = 1$.  
	
	\medskip
	Now, let $\initialx,\initialxbreve \in \X$ be such that $\initialx \preccurlyeq_A \initialxbreve$. Consider a probability space $(\Omega,\F,\PP)$ where the following are defined:
	\begin{enumerate}
		\item[(i)]
		A Poisson process $N = \{N(t), \: 0 \leq t < \infty\}$ of rate $\lambda > 0$.
		\item[(ii)]
		A sequence of independent and identically distributed (i.i.d.) random variables $U =(U_k)_{k \geq 1}$ where each $U_k$ has the uniform distribution on $[0,1]$. \end{enumerate}
	Additionally, choose $N$ to be independent of $U$. We construct two discrete-time processes, $Y=(Y_k)_{k \geq 0}$ and $\Breve{Y}=(\Breve{Y}_k)_{k \geq 0}$, by defining $Y_0 := \initialx$, $\Breve{Y}_0 := \initialxbreve$, and for $k \geq 0$,
	\begin{equation}
		\label{eq:Y_and_YbreveMAIN}
		Y_{k+1} := \Phi_{\lambda}(Y_k,U_{k+1}), \qquad  \Breve{Y}_{k+1} := \Breve{\Phi}_{\lambda}(\Breve{Y}_k,U_{k+1}).
	\end{equation}
	Then $Y$ and $\Breve{Y}$ are discrete-time Markov chains with transition matrices $P_{\lambda}(Q)$ and $P_{\lambda}(\Breve{Q})$, respectively. Now, define the processes 
	\begin{equation}
		\label{def:X_and_XbreveMAIN}
		X(t) := Y_{N(t)}, \qquad \Breve{X}(t) := \Breve{Y}_{N(t)}, \qquad t \geq 0.
	\end{equation}
	According to Section 2.1 in Keilson \cite{KeilsonRarityExponentiality} (see the discussion around Equation 2.1.6), $X$ and $\Breve{X}$ are continuous-time Markov chains with infinitesimal generators $Q$ and $\Breve{Q}$ respectively, and with initial conditions $X(0)=\initialx$ and $\Breve{X}(0)=\initialxbreve$.
	
	In order to prove \eqref{eq:CouplingCondition}, it suffices to check that the following holds:
	\begin{equation}
		\label{eq:OrderingDTMC_CRNMAIN}
		\PP[Y_k \preccurlyeq_A \Breve{Y}_k] =1, \qquad \text{ for every } k \geq 0.
	\end{equation}
	Indeed, if this is true, then  $\PP[Y_k \preccurlyeq_A \Breve{Y}_k  \text{ for every } k \geq 0] =1$ and therefore $\PP[Y_{N(t)} \preccurlyeq_A \Breve{Y}_{N(t)} \text{ for every }$ $t \geq 0] =1$. We will prove \eqref{eq:OrderingDTMC_CRNMAIN} by induction on $k$. We already know that $\initialx \preccurlyeq_A \initialxbreve$ and so \eqref{eq:OrderingDTMC_CRNMAIN} holds for $k=0$. Now, assume $\PP[Y_k \preccurlyeq_A \Breve{Y}_k] =1$ for some $k \geq 0$. Since conditions \eqref{eq:CouplingCondition1} and \eqref{eq:CouplingCondition2} hold for every $x,y \in \X$ such that $x \preccurlyeq_A y$, by Lemma \ref{lem:PhiLambda_ineq} we obtain that on a set of probability one, on which $Y_k \preccurlyeq_A \Breve{Y}_{k}$,
	\begin{equation}
		Y_{k+1} =\Phi_{\lambda}(Y_k,U_{k+1}) \preccurlyeq_A \Breve{\Phi}_{\lambda}(\Breve{Y}_k,U_{k+1}) = \Breve{Y}_{k+1},
	\end{equation}
	and so \eqref{eq:OrderingDTMC_CRNMAIN} holds with $k+1$ in place of $k$. This completes the induction step and so Theorem \ref{thm:MainResult} is proved whenever \eqref{eq:uniformization_condition_intensity} holds.
	
	\medskip
	For the case where \eqref{eq:uniformization_condition_intensity} does not hold, we construct the corresponding continuous-time Markov chains as a limit in distribution of appropriately coupled continuous-time Markov chains with truncated propensity functions for which \eqref{eq:uniformization_condition_intensity} holds. Many elements for this case are similar to the previous case, although the use of Lemma \ref{lem:PhiLambda_ineq} is different. We provide the details below, where we assume that the hypotheses of Theorem \ref{thm:MainResult} hold.
	
	We consider truncations of the propensity functions $\rate$ and $\Breve{\rate}$. More concretely, for $\initialx,\initialxbreve \in \X$ such that $\initialx \preccurlyeq_A \initialxbreve$, let $M_0 \geq 1$ be an integer such that $\norm{\initialx}_{\infty},\norm{\initialxbreve}_{\infty} \leq M_0$. For every integer $M \geq M_0$, consider the finite set $\X_{M} := \{ x \in \X \:|\: \norm{x}_{\infty} \leq M\}$, together with the functions $\rate^M_j,\Breve{\rate}^M_j : \X \longrightarrow \R_+$ defined by $\rate^M_j(x) := \rate_j(x)\one_{\X_{M}}(x)$ and $\Breve{\rate}^M_j(x) := \Breve{\rate}_j(x)\one_{\X_{M}}(x)$ for $1 \leq j \leq n$ and $x \in \X$. We see that for every $M \geq M_0$, \eqref{eq:LambdaNotOutofSpace} holds with the functions $\rate^M=(\rate^M_1, \dots,\rate^M_n)$ and  $\Breve{\rate}^M= (\Breve{\rate}^M_1, \dots,\Breve{\rate}^M_n)$ in place of $\rate$ and $\Breve{\rate}$. Also, since $\X_M$ is a finite set, $\sup_{x \in \X} \rate^M_j(x) = \sup_{x \in \X_M} \rate_j(x)  < \infty$ and $\sup_{x \in \X} \Breve{\rate}^M_j(x)  = \sup_{x \in \X_M} \Breve{\rate}_j(x) < \infty$ for every $1 \leq j \leq n$. Furthermore, by \eqref{eq:CouplingCondition1} and \eqref{eq:CouplingCondition2}, we have that for every pair $x,y \in \X_M$ such that $x \preccurlyeq_A y$,
	\begin{equation}
		\label{eq:PropensityFunctionsM}
		\begin{split}
			\Breve{\rate}^M_j(y) \leq \rate^M_j(x), \quad \text{for every } 1 \leq j \leq n \text{ such that } y +v_j \in \X \setminus (K_A + x), \text{ and } \\
			\Breve{\rate}^M_j(y)\geq \rate^M_j(x), \quad \text{for every } 1 \leq j \leq n \text{ such that } x +v_j \in \X \text{ and } y \notin K_A + x + v_j.
		\end{split}
	\end{equation}
	Let $Q^M$ and $\Breve{Q}^M$ denote the infinitesimal generators associated with $\rate^M$ and $\Breve{\rate}^M$ respectively. We define an increasing sequence $\{\lambda_M\}_{M \geq M_0}$ of positive numbers such that $\lambda_{M} \longrightarrow \infty$ as $M \longrightarrow \infty$ and $\lambda_M > n\max\left\{\sup_{x \in \X} \sum_{j=1}^n\rate^M_j(x), \sup_{x \in \X} \sum_{j=1}^n\Breve{\rate}^M_{j}(x) \right\}$ for every $M \geq M_0$. Define $\Phi_{\lambda_M}(\cdot,\cdot),\Breve{\Phi}_{\lambda_M}(\cdot,\cdot): \X \times [0,1] \longrightarrow \X$ as in \eqref{eq:DefinitionPhiLambda}, but with $\rate^M$ and $\Breve{\rate}^M$ in place of $\rate$ and $\Breve{\rate}$, respectively. Since \eqref{eq:PropensityFunctionsM} holds, applying Lemma \ref{lem:PhiLambda_ineq} with $\rate^M,\Breve{\rate}^M,\lambda_M,\Phi_{\lambda_M},\Breve{\Phi}_{\lambda_M}$ in place of $\rate,\Breve{\rate},\lambda,\Phi_{\lambda},\Breve{\Phi}_{\lambda}$ yields that
	\begin{equation}
		\label{eq:ComparisonOfPhiLambdaM}
		\Phi_{\lambda_M}(x,u) \preccurlyeq_A  \Breve{\Phi}_{\lambda_M}(y,u) \quad \text{ for every } x,y \in \X_M \text{ such that } x \preccurlyeq_A y \text{ and } u \in [0,1].
	\end{equation}
	
	Now, for each $M \geq M_0$ consider a probability space $(\Omega^M,\F^M,\PP^M)$ where the following are defined:
	\begin{enumerate}
		\item[(i)]
		A Poisson process $N^{M} = \{N^M(t), \: 0 \leq t < \infty\}$ of rate $\lambda_M > 0$.
		\item[(ii)]
		An i.i.d. sequence $U^M =(U^M_{k})_{k \geq 1}$ of uniform $[0,1]$ random variables.
	\end{enumerate}
	Additionally, choose $N^M$ to be independent of $U^M$. For every $M \geq M_0$, we construct two discrete-time processes, $Y^{M}=(Y^M_k)_{k \geq 0}$ and $\Breve{Y}^M=(\Breve{Y}^M_k)_{k \geq 0}$, by defining $Y^M_0 := \initialx$, $\Breve{Y}^M_0 := \initialxbreve$ and for $k \geq 0$,
	\begin{equation}
		Y^M_{k+1} := \Phi_{\lambda_M}(Y^M_k,U^M_{k+1}), \qquad  \Breve{Y}^M_{k+1} := \Breve{\Phi}_{\lambda_M}(\Breve{Y}^M_k,U^M_{k+1}).
	\end{equation}
	Similarly to the previous case, $Y^M$ and $\Breve{Y}^M$ are discrete-time Markov chains with transition matrices $P_{\lambda_M}(Q^{M}) := \frac{1}{\lambda_M}Q^{M} + I$ and $P_{\lambda_M}(\Breve{Q}^M) := \frac{1}{\lambda_M}\Breve{Q}^{M} + I$, respectively.
	
	Now, we claim that for each $M \geq M_0$:
	\begin{equation}
		\label{eq:InequalitieYMMAIN}
		\PP^M\left[Y^M_{k \wedge S^M} \preccurlyeq_A \Breve{Y}^M_{k \wedge S^M} \text{ for every } k \geq 0 \right] =1.
	\end{equation}
	where $S^M := \inf\{ k \geq 0 \:|\: Y^M_k \notin \X_M \text{ or } \Breve{Y}^M_k \notin \X_M\}$. In fact, \eqref{eq:InequalitieYMMAIN} is equivalent to proving that $\PP^M\left[Y^M_{k \wedge S^M} \preccurlyeq_A \Breve{Y}^M_{k \wedge S^M} \right] =1$ for every $k \geq 0$, which we do by induction. We already know that $Y^M_0 \preccurlyeq_A \Breve{Y}^M_0$. Assuming the statement is true for some $k \geq 0$, to establish it for $k+1$ we distinguish between two cases. First, on $\{S^M \leq k\}$, $Y^M_{(k+1) \wedge S^M} = Y^M_{k \wedge S^M} \preccurlyeq_A \Breve{Y}^M_{k \wedge S^M} =\Breve{Y}^M_{(k+1) \wedge S^M}$, $\PP^M$-a.s.. Second, on $\{S^M > k\}$, $Y^M_k \in \X_M, \Breve{Y}^M_k \in \X_M$, and by the induction assumption, $Y^M_k \preccurlyeq_A \Breve{Y}^M_k$, $\PP^M$-a.s.. Applying Lemma \ref{lem:PhiLambda_ineq}, we obtain $\PP^M$-a.s. on $\{S^M > k\}$ that
	\begin{equation}
		Y^M_{(k+1) \wedge S^M} = Y^M_{k+1} =\Phi_{\lambda_M}(Y^M_k,U^M_{k+1}) \preccurlyeq_A \Breve{\Phi}_{\lambda_M}(\Breve{Y}^M_k,U^M_{k+1}) =\Breve{Y}^M_{(k+1) \wedge S^M},
	\end{equation}
	where we have used \eqref{eq:ComparisonOfPhiLambdaM}.
	
	Now, for each $M \geq M_0$, we define the processes
	\begin{equation}
		X^M(t) := Y^M_{N^M(t)}, \qquad \Breve{X}^{M}(t) := \Breve{Y}^M_{N^M(t)}, \qquad t \geq 0.
	\end{equation}
	Then, $X^{M}$ and $\Breve{X}^{M}$ are continuous-time Markov chains with infinitesimal generators $Q^M$ and $\Breve{Q}^M$ respectively, and with initial conditions $X^{M}(0)=\initialx$ and $\Breve{X}^{M}(0)=\initialxbreve$. Define $T^M := \inf\{ t \geq 0 \:|\: X^{M}(t) \notin \X_M \text{ or } \Breve{X}^{M}(t) \notin \X_M\}$ and, because $Y^M$ and $\Breve{Y}^M$ are the discrete time skeletons for $X^M$ and $\Breve{X}^M$, we have that $\PP^M$-a.s.
	\begin{equation}
		T^M = \inf\{ t \geq 0 \:|\: N^M(t) =S^M \}.
	\end{equation}
	Then, it follows from \eqref{eq:InequalitieYMMAIN} that
	\begin{equation}
		\label{eq:StoppedXaredominated}
		\PP^M\left[X^{M}(t \wedge T^M) \preccurlyeq_A \Breve{X}^{M}(t \wedge T^M) \text{ for every } t \geq 0 \right] =1. 
	\end{equation}
	
	We now prove that for every $t \geq 0$,
	\begin{equation}
		\label{eq:TM_grows_to_infty}
		\PP^M[T^M < t ] \longrightarrow 0, \quad \text{ as } M \longrightarrow \infty.
	\end{equation}
	For this, let $T^M_{X^M} := \inf\{ t \geq 0 \:|\: X^{M}(t) \notin \X_M \}$  and $T^M_{\Breve{X}^M} := \inf\{ t \geq 0 \:|\: \Breve{X}^{M}(t) \notin \X_M \}$. Since $T^M =T^M_{X^M} \wedge T^M_{\Breve{X}^M}$, then
	\begin{equation}
		\label{eq:BoundT^M}
		\PP^M[T^M < t ] \leq  \PP^M[T^M_{X^M} < t ] +  \PP^M[T^M_{\Breve{X}^M} < t], \qquad \text{ for every } t \geq 0.
	\end{equation}
	Now, since $Q^M_{x,y}=Q_{x,y}$ for $x \in \X_M$ and $y \in \X$, $X^M(\cdot \wedge T^M_{X^M})$ will have the same distribution as a Markov chain with infinitesimal generator $Q$ and initial condition $\initialx$, stopped at the first time it leaves $\X_M$. Because of this, $T^M_{X^M}$ has the same distribution as the first time a continuous-time Markov chain with infinitesimal generator $Q$ leaves $\X_M$. Since a continuous-time Markov chain with infinitesimal generator $Q$ has been assumed to not explode in finite time, we obtain that $\PP^M[T^M_{X^M} < t ] \longrightarrow 0$ as $M \to \infty$. Similar reasoning holds for $T^M_{\Breve{X}^M}$. Combining with \eqref{eq:BoundT^M}, we obtain \eqref{eq:TM_grows_to_infty}.
	
	Denote by $\D([0,\infty),\X^2)$ the space of right-continuous functions from $[0,\infty)$ into $\X^2$ that also have finite left-limits. As usual, this space is endowed with Skorokhod's $J_1$ topology. The pair $(X^M,\Breve{X}^M)$ have paths in $\D([0,\infty),\X^2)$ and we obtain $(X,\Breve{X})$ as a limit in distribution of $(X^M,\Breve{X}^M)$ as $M \to \infty$. We first verify that the sequence of processes $\{(X^{M},\Breve{X}^{M})\}_{M \geq M_0}$ is tight. For this, it suffices to check that each sequence $\{X^{M}\}_{M \geq M_0}$ and $\{\Breve{X}^{M}\}_{M \geq M_0}$ is tight, which we do by means of Theorem 7.2 in Chapter 3 of Ethier \& Kurtz \cite{EthierKurtz}. Condition $(a)$ there (compact containment) is satisfied, because of \eqref{eq:TM_grows_to_infty} and because for $\tilde{M} \geq M \geq M_0$ we have that $X^{\tilde{M}}(\cdot \wedge T^M_{X^{\tilde{M}}})$ under $\PP^{\tilde{M}}$ has the same law as $X^M(\cdot \wedge T^M_{X^M})$ under $\PP^M$, where $T^M_{X^{\tilde{M}}}:= \inf\{t \geq 0 \:|\: X^{\tilde{M}}(t) \notin \X_M\}$. To verify condition $(b)$ in Theorem 7.2 of \cite{EthierKurtz}, for $t_0 > 0$ fixed and $\eta > 0$, let $M_{\eta} \geq M_0$ be such that $\PP^M[T^M_{X^M} < t_0] \leq \frac{\eta}{2}$ for all $M \geq M_{\eta}$. Then,
	\begin{align*}
		\PP^M[w'(X^M,\delta,t_0) \geq \eta] &\leq  \PP^M[w'(X^M,\delta,t_0) \geq \eta \:;\: T^M_{X^M} \geq t_0] + \PP^M[T^M_{X^M} < t_0] \\
		&\leq \tilde{\PP}[w'(\tilde{X},\delta,t_0) \geq \eta \:;\: \tau_{\tilde{X}}^M \geq t_0] + \frac{\eta}{2} \\
		&\leq \tilde{\PP}[w'(\tilde{X},\delta,t_0) \geq \eta] + \frac{\eta}{2},
	\end{align*}
	where  $w'(\cdot,\cdot,\cdot)$ is the modulus of continuity, as defined in Equation (6.2), Chapter 3 of \cite{EthierKurtz},  $\tilde{X}$ under $\tilde{\PP}$ is a realization of the Markov chain associated with the infinitesimal generator $Q$ that starts with $x^{\circ}$, and  $\tau_{\tilde{X}}^M := \inf\{ t \geq 0 \:|\: \tilde{X}(t) \notin \X_M \}$. Since $\tilde{X}$ under $\tilde{\PP}$ is a single process with right-continuous paths having finite left-limits, the tightness applies to it and so the term $\tilde{\PP}[w'(\tilde{X},\delta,t_0) \geq \eta]$ can be made less than $\frac{\eta}{2}$ by choosing $\delta$ sufficiently small and so condition $(b)$ of Theorem 7.2 of \cite{EthierKurtz} is satisfied. It follows that $\{X^{M}\}_{M \geq M_0}$ is tight.  Similar reasoning yields tightness for $\{\Breve{X}^{M}\}_{M \geq M_0}$. 
	
	It follows that there exists a probability space $(\Omega,\F,\PP)$ with two processes $X$ and $\Breve{X}$ defined there, having paths that are right-continuous with finite left-limits, and a subsequence $\{M_k\}_{k \geq 1}$ such that $M_k \to \infty$ as $k \to \infty$, and the sequence $\{(X^{M_k},\Breve{X}^{M_k})\}_{k \geq 1}$ converges in distribution to the pair of processes $(X,\Breve{X})$. To identify the law of the limit, note that since $\{Q^{M_k}\}_{k \geq 1}$ converges pointwise to $Q$, for any function $f$ with bounded support in $\X$, $f(X(t)) - \int_0^t Qf(X(s))ds$ will inherit the martingale property of $f(X^{M_k}(t)) - \int_0^t Q^{M_k}f(X^{M_k}(s))ds$. It follows from the martingale characterization that $X$ is a continuous-time Markov chain with infinitesimal generator $Q$ (see Chapter 4 in Ethier \& Kurtz \cite{EthierKurtz}). Similarly, $\Breve{X}$ will be a continuous-time Markov chain with infinitesimal generator $\Breve{Q}$. In addition, the processes have inherited initial conditions $X(0)=\initialx$ and $\Breve{X}(0)=\initialxbreve$.
	
	Finally, to show that \eqref{eq:CouplingCondition} holds, consider the set
	\begin{equation}
		F = \{ (f,g) \in \D([0,\infty),\X^2) \:|\: f(t) \preccurlyeq_A g(t) \text{ for all } t \geq 0 \},
	\end{equation}
	which is closed in the Skorokhod topology. From \eqref{eq:StoppedXaredominated} we know that the stopped processes satisfy $\PP^{M_k}[(X^{M_k}(\cdot \wedge T^{M_k}),\Breve{X}^{M_k}(\cdot \wedge T^{M_k})) \in F] =1$ for every $k \geq 1$. Furthermore, from \eqref{eq:TM_grows_to_infty} we know that $T^{M_k} \longrightarrow \infty$ in probability as $k \to \infty$. The reader may verify that this last fact, along with the convergence of $(X^{M_k},\Breve{X}^{M_k})$ to $(X,\Breve{X})$, implies that $(X^{M_k}(\cdot \wedge T^{M_k}),\Breve{X}^{M_k}(\cdot \wedge T^{M_k}))$ converges in distribution to $(X,\Breve{X})$ as $k \to \infty$. By the Portmanteau Theorem (see Theorem 2.1 in Billingsley \cite{Billingsley2nd}),
	\begin{equation}
		1=\limsup_{k \to \infty} \PP^{M_k}[(X^{M_k}(\cdot \wedge T^{M_k}),\Breve{X}^{M_k}(\cdot \wedge T^{M_k})) \in F] \leq \PP[(X,\Breve{X}) \in F]  
	\end{equation}
	and we obtain \eqref{eq:CouplingCondition}.
	
	\begin{remark}
		The proof of Theorem \ref{thm:MainResult} provides a method to simulate the sample paths for the continuous-time Markov chains $X$ and $\Breve{X}$ in a coupled manner for the case where \eqref{eq:uniformization_condition_intensity} holds. Roughly speaking, the procedure consists of determining $\lambda > 0$ as in \eqref{eq:LambdaCondition}, $\Phi_{\lambda},\Breve{\Phi}_{\lambda}$ as in \eqref{eq:DefinitionPhiLambda}, $Y,\Breve{Y}$ as in \eqref{eq:Y_and_YbreveMAIN} and $X,\Breve{X}$ as in \eqref{def:X_and_XbreveMAIN}. For the benefit of the reader, this method is described as an algorithm in SI - Section S.4, which yields coupled sample paths under the assumptions of Theorem  \ref{thm:InnerProductTheorem}, \ref{thm:SumTheorem2} and S.2.
	\end{remark}

	\subsection{Proof of Theorem \ref{thm:InnerProductTheorem}}    
	\label{sec:ProofInnerProductTheorem}

	By Theorem \ref{thm:MainResult}, it suffices to prove that for every $x,y \in \X$ such that $x \preccurlyeq_A y$, conditions \eqref{eq:CouplingCondition1} and  \eqref{eq:CouplingCondition2} hold. For this, we make some observations first. Consider $x,y \in \X$ such that $x \preccurlyeq_A y$ and let $1 \leq j \leq n$. Observe that $x \preccurlyeq_A y + v_j$ will hold if and only if $A(y+v_j -x) \geq 0$ which is equivalent to:
	\begin{equation}
		\label{eq:RelationAdotted_yplusvj}
		\inn{A_{i\bullet},y -x} + \inn{A_{i\bullet},v_j} \geq 0, \qquad \text{for every } 1 \leq i \leq m.
	\end{equation}
	Similarly, $x + v_j \preccurlyeq_A y$ will hold if and only if
	\begin{equation}
		\label{eq:RelationAdotted_xplusvj}
		\inn{A_{i\bullet},y -x} - \inn{A_{i\bullet},v_j} \geq 0, \qquad \text{for every } 1 \leq i \leq m.
	\end{equation}
	Since $x \preccurlyeq_A y$, then $\inn{A_{i\bullet},y-x} \geq 0$ for every $1 \leq i \leq m$. Now, consider $i \in \{1,...,m\}$ such that $\inn{A_{i\bullet},y-x} >0$. Since $A \in \Z^{m \times d}$ and $y-x \in \Z^d$, then $\inn{A_{i\bullet},y-x} \geq 1$. This yields that
	\begin{equation}
		\label{positivecond}
		\inn{A_{i\bullet},y-x} + \inn{A_{i\bullet},v_j} \geq 1 + \inn{A_{i\bullet},v_j}\ge 0,
	\end{equation}
	since $\inn{A_{i\bullet},v_j} \in \{-1,0,1\}$. Similarly,  $\inn{A_{i\bullet},y-x} - \inn{A_{i\bullet},v_j} \geq 1 - \inn{A_{i\bullet},v_j} \geq 0$. By observing that the interior of $K_A +x$ is of the form $\interior(K_A + x) = \{ y \in \R^d \:|\: Ax < Ay\}$, the latter argument shows that for every $x \in \X$ and $y \in \interior(K_A + x) \cap \X$, we have
	\begin{equation}
		\label{eq:JumpsInterior}
		x \preccurlyeq_A y + v_j \text{ and } x + v_j \preccurlyeq_A y, \qquad \text{ for every } 1 \leq j \leq n.
	\end{equation}
	
	Now, lets check condition \eqref{eq:CouplingCondition1}. For this, let $x,y \in \X$ be such that $x \preccurlyeq_A y$ and let $1 \leq j \leq n$ be such that $y + v_j \in \X \backslash (K_A + x)$. By \eqref{eq:JumpsInterior}, $y \notin \interior(K_A + x)$ and since $y \in K_A + x$, we must have $y \in \partial(K_A + x)= \{ z \in K_A +x \:|\: \inn{A_{i\bullet},z} = \inn{A_{i\bullet},x}\; \mathrm{for\; some\;} 1\le i \le m\}$, the boundary of $K_A + x$. Consider the set of indices $\Kb_{y} := \{ i\:|\: \inn{A_{i\bullet},y} = \inn{A_{i\bullet},x}, 1\le i\le m \}$, which is non-empty. Observe that for every $i \notin \Kb_y$, $\inn{A_{i\bullet},y-x} > 0$ and from \eqref{positivecond}, $\inn{A_{i\bullet},(y+v_j)-x} \geq 0$, while for $i \in \Kb_y$, $\inn{A_{i\bullet},(y+v_j)-x} = \inn{A_{i\bullet},v_j}$. From this, we can infer that there exists an $i_k \in \Kb_y$ such that $\inn{A_{i_k\bullet},v_j} < 0$. Indeed, if this was not the case, then $\inn{A_{i\bullet},(y+v_j)-x} \geq 0$ for every $i \in \Kb_y$ and consequently \eqref{eq:RelationAdotted_yplusvj} would hold. This contradicts the fact that $y + v_j \notin K_A + x$. By \eqref{eq:CouplingConditions_InnerProduct_I}, we know that $\inn{A_{i_k\bullet},v_j} < 0$ implies $\Breve{\rate}_j(y) \leq \rate_j(x)$ and we conclude that \eqref{eq:CouplingCondition1} holds.
	
	To check condition \eqref{eq:CouplingCondition2}, let $x,y \in \X$ be such that $x \preccurlyeq_A y$ and let $1 \leq j \leq n$ be such that
	$x + v_j \in \X$ and $y \notin K_A + x + v_j$. Again, by \eqref{eq:JumpsInterior}, we obtain that $y \in \partial(K_A + x)$ and $\Kb_{y} \neq \emptyset$. For every $i \notin \Kb_y$, $\inn{A_{i\bullet},y-(x+v_j)} \geq 0$, while for $i \in \Kb_y$, $\inn{A_{i\bullet},y-(x+v_j)} = -\inn{A_{i\bullet},v_j}$. From this, we can infer that there exists an $i_k \in \Kb_y$ such that $\inn{A_{i_k\bullet},v_j} > 0$. By \eqref{eq:CouplingConditions_InnerProduct_II}, we know that $\inn{A_{i_k\bullet},v_j} > 0$ implies $\Breve{\rate}_j(y) \geq \rate_j(x)$ and we conclude that \eqref{eq:CouplingCondition2} holds.
	
	\subsection{Proof of Theorem \ref{thm:SumTheorem2}}     
	\label{sec:ProofSumTheorem2}
	
	The proof of this result uses similar general ideas to the ones used in the proof of Theorem \ref{thm:MainResult}. However, since the conditions involve sums, the construction is somewhat different and more complex and we provide the details below. Let us consider again a non-empty set $\X \subseteq \Z_+^d$, a collection of distinct vectors $v_1,\ldots,v_n$ in $\Z^d \setminus \{0\}$ and two collections of non-negative functions on $\X$, $\rate=(\rate_1, \dots,\rate_n)$ and  $\Breve{\rate}= (\Breve{\rate}_1, \dots,\Breve{\rate}_n)$ such that \eqref{eq:LambdaNotOutofSpace} holds. In the following, let $A \in \Z^{m \times d}$ be a matrix with non-zero rows such that condition (i) of Theorem \ref{thm:SumTheorem2} holds. 
	
	We initially assume that $\sup_{x \in \X} \rate_j(x) < \infty$ and $\sup_{x \in \X} \Breve{\rate}_j(x) < \infty$ for every $1 \leq j \leq n$, and let $\lambda > 0$ such that \eqref{eq:LambdaCondition} holds. We shall relax these assumptions later. We start by defining functions analogous to $\Phi_{\lambda}$ and $\Breve{\Phi}_{\lambda}$ as defined in \eqref{eq:DefinitionPhiLambda}, although this time, the construction is more involved. 
	
	Recall that $s$ denotes the size of the set $\{Av_j \:|\: 1 \leq j \leq n\}$ and that the index sets $G^k \ne \emptyset$, $1 \le k \le s$, defined in \eqref{eq:sets_of_dot_product+-13}, are such that $Av_j =\eta^k$ for all $j \in G^k$, $1 \le k \le s$. Consider a bijection $\sigma:\{1,\ldots,n\} \longrightarrow \{1,\ldots,n\}$ such that the vectors $v_{\sigma(1)},\ldots, v_{\sigma(n)}$ have the property that the first $|G^1|$ vectors have indices in $G^1$, the next $|G^2|$ vectors have indices in $G^2$, and so on. More precisely, the bijection $\sigma$ is such that for $1 \leq k \leq s$, $Av_{\sigma(q)}=\eta^k$, whenever $\sum_{\ell=1}^{k-1}|G^{\ell}|+1 \leq q \leq \sum_{\ell=1}^{k}|G^{\ell}|$. Recall for this that a sum over an empty set is taken to equal zero.

	For $x \in \X$, we define a family of intervals $\{I^k(x)\:|\: 1 \leq k \leq s\}$ as follows.
	Let $p_0 := 0$, and for $1 \leq k \leq s$, inductively define $p_k :=\sum_{\ell=1}^{k}|G^{\ell}|$, and
	\begin{equation}
		\label{Int1}
		I^k(x) := \bigcup_{q=p_{k-1}+1}^{p_k} I^k_{q}(x),
	\end{equation}
	where for $p_{k-1}+1\le q \le p_k$,
	\begin{equation}
		\label{Int2}
		I^k_{q}(x):=\left[\frac{p_{k-1}}{n}+\sum_{\ell=p_{k-1}+1}^{q-1}\frac{\rate_{\sigma(\ell)}(x)}{\lambda},\frac{p_{k-1}}{n} +\sum_{\ell=p_{k-1}+1}^{q}\frac{\rate_{\sigma(\ell)}(x)}{\lambda}\right).
	\end{equation}
	The sets $I^k_q(x)$, with $1\le k\le s$ and $p_{k-1}+1 \leq q \leq p_k$, are mutually disjoint, and by \eqref{eq:LambdaCondition}, the length of $I^k(x)$ is less than $\frac{p_k-p_{k-1}}{n}=\frac{|G^k|}{n}$, and so the sum of the lengths of $\{I^k(x)\:|\: 1 \leq k \leq s\}$ is less than $\frac{1}{n} \sum_{k=1}^s |G^k|=1$. Now, let us define $\PPhi_{\lambda}(\cdot,\cdot): \X \times [0,1] \longrightarrow \X$ by 
	\begin{equation}\label{eq:DefinitionPhiLambdaBIS2A}
		\PPhi_{\lambda}(x,u) := x + \sum_{k=1}^{s} \sum_{q=p_{k-1}+1}^{p_k} v_{\sigma(q)}\one_{I^k_q(x)}(u), \qquad x \in \X, \; u \in [0,1].
	\end{equation}        
	Note that $Av_{\sigma(q)}=\eta^k$ for $p_{k-1}+1 \leq q \leq p_k$, $1 \le k \le s$. From the above properties of $I^k_q(x)$, we have that for any $u \in [0,1]$,
	either $u \notin \bigcup_{k=1}^{s} I^k(x)$ or $u \in I^k_{q}(x)$ for exactly one $k$ and $q$ such that $I^k_{q}(x) \ne \emptyset$. The latter condition implies, by \eqref{Int2}, that $\rate_{\sigma(q)}(x) > 0$ and then, by \eqref{eq:LambdaNotOutofSpace}, $x + v_{\sigma(q)} \in \X$. This shows that $\PPhi_{\lambda}(\cdot,\cdot)$ is well-defined as an $\X$-valued function. 
	
	In an analogous manner to that above, we can define intervals $\Breve{I}^k(x)$, $\Breve{I}^k_q(x), \: 1 \leq k \leq s, \: p_{k-1}+1 \leq q \leq p_k, \: x \in \X$ and a function $\Breve{\PPhi}_{\lambda}:\X \times [0,1] \longrightarrow \X$, as in \eqref{Int1} -- \eqref{eq:DefinitionPhiLambdaBIS2A}, but with $\Breve{\rate}_j(x)$, $\Breve{I}^k(x)$, $\Breve{I}^k_q(x)$,  $\Breve{\PPhi}_{\lambda}$ in place of $\rate_j(x)$, $I^k(x)$, $I^k_q(x)$, $\PPhi_{\lambda}$.

	\begin{lemma}
		\label{lem:PhiLambda_ineqBIS2}
		Suppose that $x,y \in \X$ are such that $x \preccurlyeq_A y$ and the following hold: whenever $y \in \partial_i(K_A+x) \cap \X$ for some $1 \leq i \leq m$, we have
		\begin{equation}
			\label{eq:CouplingConditionBISPROOF22}
			\sum_{j \in G^{k}} \Breve{\rate}_j(y) \leq \sum_{j \in G^{k}} \rate_j(x), \quad \text{for every } k \text{ such that } \eta^k_i <0,
		\end{equation}
		and  
		\begin{equation}
			\label{eq:CouplingConditionBISPROOF21}
			\sum_{j \in G^{k} } \Breve{\rate}_j(y) \geq \sum_{j \in G^{k}} \rate_j(x), \quad \text{for every } k \text{ such that } \eta^k_i >0.
		\end{equation}
		Then, for each $u \in [0,1]$,
		\begin{equation}
			\label{eq:ClaimPhi_Lambda_II2}
			\PPhi_{\lambda}(x,u) \preccurlyeq_A  \Breve{\PPhi}_{\lambda}(y,u).
		\end{equation}
	\end{lemma}
	\begin{proof}
		First, we note that $\PPhi_{\lambda}, \Breve{\PPhi}_{\lambda}$ have the following properties: for every $u \in [0,1]$, $1\le k \le s$, $j \in G^k$, 
		\begin{equation}
			\text{ if } \PPhi_{\lambda}(x,u) = x + v_j \text{, then } \Breve{\PPhi}_{\lambda}(y,u) \in \{y + v_\ell: \ell \in G^k \} \cup \{y\},
		\end{equation}
		since $I_{\sigma^{-1}(j)}^k(x), \Breve{I}_{\sigma^{-1}(\ell)}^k(y) \subseteq [\frac{p_{k-1}}{n},\frac{p_{k}}{n})$ for $\ell \in G^k$. Similarly,
		\begin{equation}
			\label{eq:PhiLambdaCondition1BIS00}
			\text{ if } \Breve{\PPhi}_{\lambda}(y,u) = y +v_j \text{, then } \PPhi_{\lambda}(x,u) \in \{x + v_\ell: \ell \in G^k\} \cup \{x\}.
		\end{equation}
		Furthermore, for $1\le k \le s$, $j \in G^k$, if $\sum_{\ell \in G^{k}} \Breve{\rate}_\ell(y) \geq \sum_{\ell \in G^{k}} \rate_\ell(x)$, then 
		\begin{equation}
			\PPhi_{\lambda}(x,u) = x +v_j \text{ implies that } \Breve{\PPhi}_{\lambda}(y,u) =y +v_\ell \: \text{ for some } \: \ell \in G^k,
		\end{equation}
		since under the condition, $I^k(x) \subseteq \Breve{I}^k(y)$. Similarly, if $ \sum_{\ell \in G^{k}} \Breve{\rate}_\ell(y) \leq \sum_{\ell \in G^{k}} \rate_\ell(x)$, then
		\begin{equation}
			\Breve{\PPhi}_{\lambda}(y,u) = y +v_j \text{ implies that } \PPhi_{\lambda}(x,u) =x +v_\ell \: \text{ for some } \: \ell \in G^k.
		\end{equation}
		We also have that, for $1 \le k \le s$ and $j \in G^k$, $x \preccurlyeq_A y + v_j$ if and only if
		\begin{equation}
			\label{eq:RelationAdotted_yplusvjBIS}
			\inn{A_{i\bullet},y -x} + \inn{A_{i\bullet},v_j} \geq 0, \qquad \text{for every } 1 \leq i \leq m.
		\end{equation}
		Similarly, $x + v_j \preccurlyeq_A y$ if and only if
		\begin{equation}
			\label{eq:RelationAdotted_xplusvjBIS}
			\inn{A_{i\bullet},y -x} - \inn{A_{i\bullet},v_j} \geq 0, \qquad \text{for every } 1 \leq i \leq m.
		\end{equation}
		Furthermore, for $1 \le k \le s$ and $j,\ell \in G^k$, since $Av_j=Av_{\ell}$ and $x \preccurlyeq_A y$, then
		\begin{equation}
			\label{eq:RelationAdotted_xplusvjyplusvjBIS}
			\inn{A_{i\bullet},y -x} + \inn{A_{i\bullet},(v_j-v_\ell)}= \inn{A_{i\bullet},y-x} \geq 0, \qquad \text{for every } 1 \leq i \leq m.
		\end{equation}
		
		To prove \eqref{eq:ClaimPhi_Lambda_II2}, we first consider the situation where $y \in \interior(K_A + x) = \{ w \in \R^d \:|\: Ax < Aw\}$. Then, for each $1 \leq i \leq m$,  $\inn{A_{i\bullet},y-x} >0$ and since $A \in \Z^{m \times d}$ and $y-x \in \Z^d$, we have $\inn{A_{i\bullet},y-x} \geq 1$. This implies that for $1 \le k \le s$ and $j \in G^k$,
		\begin{equation}
			\label{positivecondBIS1}
			\inn{A_{i\bullet},y-x} + \inn{A_{i\bullet},v_j} \geq 1 + \inn{A_{i\bullet},v_j}\ge 0, \qquad \text{for every } 1 \leq i \leq m,
		\end{equation}
		since $\inn{A_{i\bullet},v_j} \in \{-1,0,1\}$ by condition $(i)$ of Theorem \ref{thm:SumTheorem2}. Similarly, for $1 \le k \le s$ and $j \in G^k$,
		\begin{equation}
			\label{positivecondBIS2}
			\inn{A_{i\bullet},y-x} - \inn{A_{i\bullet},v_j} \geq 1 - \inn{A_{i\bullet},v_j} \geq 0, \qquad \text{for every } 1 \leq i \leq m.
		\end{equation}
		It follows from \eqref{eq:RelationAdotted_xplusvjyplusvjBIS} -- \eqref{positivecondBIS2} that
		if $y \in \interior(K_A + x) \cap \X$, then for any $1 \leq k \leq s$ and $j, \ell \in G^k$:
		\begin{equation}
			\label{eq:JumpsInteriorBIS}
			x \preccurlyeq_A y + v_j, \: x + v_j \preccurlyeq_A y \: \text{ and } \: x + v_\ell \preccurlyeq_A y + v_j.
		\end{equation}
		We also have, by assumption, that $x  \preccurlyeq_A y$. It follows that if $y \in \interior(K_A + x) \cap \X$, then $\{x,x +v_\ell \:|\: \ell \in G^k\}\preccurlyeq_A \{y,y+v_j \:|\: j \in G^k\}$ for $1 \le k \le s$ and consequently \eqref{eq:ClaimPhi_Lambda_II2} holds for all $u \in [0,1]$.
		
		Now, we turn to the other situation where $y \in \partial_i(K_A+x) \cap \X$ for some $1 \leq i \leq m$. Then $\Kb_{y} := \{ i\:|\: \inn{A_{i\bullet},y} = \inn{A_{i\bullet},x}, 1\le i\le m \}$ is non-empty. Let $u\in[0,1]$. We consider two cases.
		
		\medskip
		\textbf{Case 1:} $\Breve{\PPhi}_{\lambda}(y,u) = y +v_j$ for some $1 \leq j \leq n$.
		
		\medskip
		Fix such an index $j$. Consider the unique $1 \leq k \leq s$ such that $j \in G^k$. Then, by  \eqref{eq:PhiLambdaCondition1BIS00}, either $\PPhi_{\lambda}(x,u) = x + v_\ell$ for some $\ell \in G^k$, or $\PPhi_{\lambda}(x,u) = x$. 
		\begin{enumerate}
			\item[a)]
			If $\PPhi_{\lambda}(x,u) = x + v_\ell$ for some $\ell \in G^k$, then, since $x \preccurlyeq_A y$ and $Av_j = Av_\ell$, we have $x + v_\ell \preccurlyeq_A y +v_j$. Hence, $\PPhi_{\lambda}(x,u) \preccurlyeq_A  \Breve{\PPhi}_{\lambda}(y,u)$ and \eqref{eq:ClaimPhi_Lambda_II2} holds. 
			\item[b)]
			If $\PPhi_{\lambda}(x,u) = x$, we claim that $y + v_j \in K_A + x$. To see this, observe that for every $i \notin \Kb_y$, $\inn{A_{i\bullet},y-x} > 0$ and as for \eqref{positivecondBIS1}, $\inn{A_{i\bullet},(y+v_j)-x} \geq 0$, while for $i \in \Kb_y$, $\inn{A_{i\bullet},(y+v_j)-x} = \inn{A_{i\bullet},v_j}\in \{-1,0,1\}$. For each $i \in \Kb_y$, if $\inn{A_{i\bullet},v_j}=-1$, then by \eqref{eq:CouplingConditionBISPROOF22}, we would have $\sum_{\ell \in G^{k}} \Breve{\rate}_\ell(y) \leq \sum_{\ell \in G^{k}} \rate_\ell(x)$, which would imply that $\PPhi_{\lambda}(x,u) = x +v_\ell$ for some $\ell \in G^k$, but this contradicts the assumption that $\PPhi_{\lambda}(x,u) = x$. So we must have $\inn{A_{i\bullet},v_j} \geq 0$ and hence $\inn{A_{i\bullet},(y+v_j)-x} \geq 0$ for all $i \in \Kb_y$. Thus, $y + v_j \in K_A + x$ and so $\PPhi_{\lambda}(x,u) = x \preccurlyeq_A  y + v_j = \Breve{\PPhi}_{\lambda}(y,u)$ holds.
		\end{enumerate}
		
		\medskip
		\textbf{Case 2:} $\Breve{\PPhi}_{\lambda}(y,u) = y$. Again, we consider two subcases.
		
		\begin{enumerate}
			\item[a)]
			If $\PPhi_{\lambda}(x,u)=x$, then \eqref{eq:ClaimPhi_Lambda_II2} holds, because $x \preccurlyeq_A  y$.
			\item[b)]
			If $\PPhi_{\lambda}(x,u) = x +v_j$ for some $1 \leq j \leq n$, we claim that $y \in K_A + x + v_j$ for the corresponding value of $j$. To see this, fix the value of $j$ for which $\PPhi_{\lambda}(x,u) = x +v_j$, let $1 \leq k \leq s$ be such that $j \in G^k$, and observe that for every $i \notin \Kb_y$, $\inn{A_{i\bullet},y-x} > 0$ and as for \eqref{positivecondBIS2}, $\inn{A_{i\bullet},y-(x+v_j)} \geq 0$, while for $i \in \Kb_y$, $\inn{A_{i\bullet},y-(x+v_j)} = - \inn{A_{i\bullet},v_j}\in \{-1,0,1\}$. For each $i \in \Kb_y$, if $\inn{A_{i\bullet},v_j}=1$, then by \eqref{eq:CouplingConditionBISPROOF21}, we would have $\sum_{\ell \in G^{k} } \Breve{\rate}_\ell(y) \geq \sum_{\ell \in G^{k}} \rate_\ell(x)$, which would imply that $\Breve{\PPhi}_{\lambda}(y,u) = y +v_\ell$ for some $\ell \in G^k$. This would contradict the assumption that $\Breve{\PPhi}_{\lambda}(y,u) = y$. So we must have $\inn{A_{i\bullet},v_j} \leq 0$ and hence $\inn{A_{i\bullet},y-(x+v_j)}=\inn{A_{i\bullet},y-x}-\inn{A_{i\bullet},v_j} \geq 0$ for all $i \in \Kb_y$. Thus, we have $y \in K_A + x + v_j$ and then $\PPhi_{\lambda}(x,u) = x + v_j \preccurlyeq_A y = \Breve{\PPhi}_{\lambda}(y,u)$.
			
		\end{enumerate}  
		\vskip -0.2in    
	\end{proof}
	
	In order to prove Theorem \ref{thm:SumTheorem2}, from here on we can follow a similar procedure to the one used in the proof of Theorem \ref{thm:MainResult} after Lemma \ref{lem:PhiLambda_ineq} was proved there. For the case where \eqref{eq:uniformization_condition_intensity} holds, we define two discrete-time processes, $Y=(Y_k)_{k \geq 0}$ and $\Breve{Y}=(\Breve{Y}_k)_{k \geq 0}$, by defining $Y_0 := \initialx$, $\Breve{Y}_0 := \initialxbreve$, and for $k \geq 0$,
	\begin{equation}
		Y_{k+1} := \PPhi_{\lambda}(Y_k,U_{k+1}), \qquad  \Breve{Y}_{k+1} := \Breve{\PPhi}_{\lambda}(\Breve{Y}_k,U_{k+1}),
	\end{equation}
	and define $X$ and $\Breve{X}$ using these and an independent Poisson process $N$ as in \eqref{def:X_and_XbreveMAIN}. For the case where \eqref{eq:uniformization_condition_intensity} does not hold, we can use a truncation procedure similar to that for Theorem \ref{thm:MainResult}. In both cases, we use Lemma \ref{lem:PhiLambda_ineqBIS2} instead of Lemma \ref{lem:PhiLambda_ineq}.

	\section{Conclusion}
	\label{sec:conclusion}
	
	In this work, we first reviewed the concept of Stochastic Chemical Reaction Networks (SCRNs), a class of continuous-time Markov chain models frequently used to describe the stochastic behavior of chemical reaction systems. We also gave the definitions of preorder and increasing set considered in this paper. In Section \ref{sec:CouplingOfCRNs}, we presented the main theoretical results of this paper. We first derived, by exploiting uniformization and then coupling of stochastic processes (see Grassmann \cite{uniformization} and Keilson \cite{KeilsonRarityExponentiality}), three theorems which give practical sufficient conditions for stochastic dominance of one continuous-time Markov chain over another. More precisely, these theorems provide conditions under which, when one or more parameters is changed monotonically, the system is almost surely ``higher'' with respect to a certain preorder. While the first theorem (Theorem \ref{thm:MainResult}) can be used for any SCRN, it has extensive conditions to check. The second set of theorems (Theorems \ref{thm:InnerProductTheorem}, \ref{thm:SumTheorem2}) can be used for more specific SCRN classes, but they have assumptions that only need to be checked at the boundary of certain translated convex cones. All these theorems can be applied to SCRNs with either finite or countably many states. In Section \ref{sec:MonotonicPropMFPTandSD}, we exploited these tools to develop two theorems to specifically study the monotonicity properties of stationary distributions and mean first passage times depending on system parameters.
	
	Subsequently, in Section \ref{sec:examples}, we presented some illustrative examples to highlight the advantages of using our theoretical tools in order to study the stochastic behavior of SCRNs. Specifically, we focused on two common models for enzymatic kinetics (see Michaelis \& Menten \cite{MM1913}, Kang et al. \cite{KangKhudaBukhshKoepplRempala}, Del Vecchio \& Murray \cite{BFS} and Anderson et al. \cite{2010Anderson}), on a model inspired by Braess's paradox (see Calvert et al. \cite{braessQueuing}) and on a recently developed model describing the main interactions among histone modifications alone, and together with an expressed protein (see Bruno et al. \cite{BrunoDelVecchio}). In these illustrative examples we see that our sufficient conditions can be easy to check and our results can be also used to study networks with a countably infinite number of states. Furthermore, the conclusions obtained by using our theorems are true for trajectories of the Markov chains, yielding results for both transient and steady state behavior.
	
	Overall, in this paper we derived and presented theorems that can be used for the theoretical study of monotonicity of SCRNs associated to a variety of chemical reaction systems. Future work will include the adaptation of our theoretical tools to other forms of monotonicity for SCRNs (see Definition 5.1.1 in Muller \& Stoyan \cite{MullerStoyan} as an example), the investigation of possible correlations between the network graph properties and the monotonicity properties of the SCRN (extension of the work of Angeli et al. \cite{2006Sontag} to SCRNs), and the application of our results to deterministic chemical reaction network through appropriate limits.\\

	\noindent
	{\bf Supplementary information (SI) file:} file containing detailed mathematical derivations for some of our examples, a generalization of Theorem \ref{thm:SumTheorem2}, and an algorithm for coupled stochastic simulation.\\

	\noindent
	{\bf Acknowledgements:} S.B. was supported by NSF Collaborative Research grant MCB-2027949 (PI: D.D.V.). R.J.W., F.C and Y.F. were supported in part by NSF Collaborative Research grant MCB-2027947 (PI: R.J.W.) and by the Charles Lee Powell Foundation (PI: R.J.W.). We are grateful to the anonymous referees for several very helpful comments.\\ 
	
	\noindent
	{\bf Ethics declarations - conflict of interest:} The authors declare that they have no conflicts of interest.\\
	
	\noindent
	{\bf Data availability:} Data sharing not applicable to this article as no datasets were generated or analysed during the current study.\\

		\newpage	
	
	\appendixpageoff
	\appendixtitleoff
	\renewcommand{\appendixtocname}{Supplementary Information}
	\begin{appendices}
		\setcounter{section}{19}
		\crefalias{section}{supp}
		\setcounter{figure}{0}
		\renewcommand{\thefigure}{S.\arabic{figure}}

		\setcounter{equation}{0}
		\renewcommand{\theequation}{S.\arabic{equation}}
		
		\pagenumbering{arabic}
		\renewcommand*{\thepage}{\arabic{page}}

		\begin{center}
			{\LARGE Comparison Theorems for Stochastic Chemical Reaction Networks}\\
		\end{center}
		
		\begin{center}
			{\large Felipe A. Campos$^{1,*}$, Simone Bruno$^{2,*}$, Yi Fu$^{1}$, Domitilla Del Vecchio$^{2}$, and Ruth J. Williams$^{1}$}\\
		\end{center}

		\begin{center}
			$^{1}$\textit{\small Department of Mathematics, University of California, San Diego, 9500 Gilman Drive, La Jolla CA 92093-0112. Email: {\tt\small (fcamposv,yif064,rjwilliams)@ucsd.edu}}

			$^2$\textit{\small Department of Mechanical Engineering, Massachusetts Institute of Technology, 77 Massachusetts Avenue, Cambridge, MA 02139. Emails: {\tt\small (sbruno,ddv)@mit.edu}}
			
			$^{*}$\textit{\small These authors contributed equally: F. A. Campos and S. Bruno}
		\end{center}
		
		\section*{Supplementary Information (SI) file}

		\subsection{Criteria for Positive Recurrence and Exponential Ergodicity with Application to Examples \ref{exEK2} and \ref{ex:HistoneModificationAndProtein}}
		
		\subsubsection{Foster-Lyapunov Conditions for Positive Recurrence and Exponential Ergodicity}
		\label{sec:appendixA}
		
		Here we recall fairly general conditions for positive recurrence and exponential ergodicity of a continuous-time Markov chain. Such conditions are well known and are usually referred to as Foster-Lyapunov-type conditions. We also apply these to the Markov chains in Examples \ref{exEK2} and \ref{ex:HistoneModificationAndProtein}, in subsections \ref{appendix:example4.2} and \ref{appendix:example4.5}.
		
		\begin{theorem}
			\label{thm:FL}
			Let $X$ be an irreducible continuous-time Markov chain with state space $\X$ and infinitesimal generator $Q$. Suppose $V:\X \rightarrow \R_+$ is norm-like, that is $\{x \in \X: V(x) \leq a\}$ is compact \footnote{Since $X$ is a Markov chain, $\X$ is finite or countable, and we endow it with the usual discrete topology consisting of all subsets of $\X$.} for each $a \in \R_+$. Further assume that for some $c > 0$, $d > 0$ and a compact set $C$,
			\begin{equation}
				\label{eq:FLcriterion}
				Q V(x) \leq -c + d \one_C (x), \quad \text{ for all } x \in \X.
			\end{equation} 
			Then, $X$ is non-explosive and positive recurrent, and has a unique stationary distribution $\pi$. If instead of \eqref{eq:FLcriterion}, we have that for some $c' > 0$ and $d' > 0$, 
			\begin{equation}
				\label{eq:expFLcriterion}
				Q V(x) \leq -c' V(x) + d', \quad \text{ for all } x \in \X,
			\end{equation} 
			then \eqref{eq:FLcriterion} automatically holds and the consequences stated above hold, and in addition, the stationary distribution satisfies
			$$\pi (V) = \sum_{x \in \X} \pi_x V(x) < \infty,$$
			and there is $0 < B < \infty$ and $0 < \beta < 1$ such that for all $t \geq 0$ and $x \in \X$,
			\begin{equation}
				\label{eq:expErgodicity}
				\sum_{y \in \X} |P_{xy} (t) - \pi_y| \leq \norm{P_{x \bullet} (t) - \pi}_{V+1} \leq B (V(x) + 1) \beta^t,
			\end{equation}
			where $P_{xy} (t) = \PP[X (t)=y | X(0) = x]$ and
			\begin{equation*}
				\norm{P_{x\bullet} (t) - \pi}_{V+1} = \sup_{|g| \leq V+1} \left| \sum_{y \in \X} (P_{xy} (t) - \pi_y) g(y) \right|.
			\end{equation*}
		\end{theorem}
		
		\begin{remark}
			When the second inequality in \eqref{eq:expErgodicity} holds, we say that the Markov chain is exponentially ergodic in the $(V+1)$-norm.
		\end{remark}
		
		\begin{proof}
			We will verify the sufficient conditions for each of non-explosion, positive recurrence and exponential ergodicity given in Meyn \& Tweedie \cite{bib:MT1993}. Note that $X$ is a Borel right-process under the definition in Sharpe \cite{Sharpe}. In addition, since the state space is discrete, each compact set is finite and therefore petite.
			
			For $m \in \Z_+$, if $Q_m$ is the infinitesimal generator for the Markov chain $X$ killed upon exit from $O_m=\{x \in \X: V(x) \leq m\}$  \footnote{Upon exit from $O_m$, the killed process goes to a cemetery state $\Delta_m$ in $\X \setminus O_m$ where $V(\Delta_m)=\min\{V(x): x \in \X \setminus O_m\}$.}, then $Q_m V(x) \leq Q V(x)$ for $x \in O_m$. It then follows from \eqref{eq:FLcriterion} that conditions (CD0) and (CD2) (with $f=1$) in \cite{bib:MT1993} hold with $Q_m$ in place of $\mathscr{A}_m$ there. By Theorem 2.1 and Theorem 4.2 in \cite{bib:MT1993}, the Markov chain is non-explosive and positive recurrent, and it has a unique stationary distribution $\pi$. 
			
			On the other hand, if \eqref{eq:expFLcriterion} holds, then \eqref{eq:FLcriterion} holds, using the norm like property of $V$. Furthermore, \eqref{eq:expFLcriterion} implies that conditions (CD0), (CD2) (with $f = V+1$) and (CD3) in \cite{bib:MT1993} hold with $Q_m$ in place of $\mathscr{A}_m$ there. By Theorem 2.1, Theorem 4.2 and Theorem 6.1 in \cite{bib:MT1993}, the Markov chain is non-explosive and positive recurrent, with a unique stationary distribution $\pi$ such that $\pi(V) < \infty$, and it is exponentially ergodic in the $(V+1)$-norm, that is the second inequality in \eqref{eq:expErgodicity} holds for all $t \geq 0$ and $x \in \X$. For fixed $t \geq 0$ and $x \in \X$, setting $g (y)=\sgn(P_{xy} (t) - \pi_y)$, for $y \in \X$, we have that $|g| \leq 1 \leq V+1$, and
			\begin{equation*}
				\sum_{y \in \X} |P_{xy} (t) - \pi_y| = \left| \sum_{y \in \X} (P_{xy} (t) - \pi_y) g (y) \right| \leq \norm{P_{x \bullet} (t) - \pi}_{V+1},
			\end{equation*}
			yielding the first inequality in \eqref{eq:expErgodicity}.
		\end{proof}
		
		\subsubsection{Application to Example \ref{exEK2}}
		\label{appendix:example4.2}
		
		For Example \ref{exEK2}, we first show that the Markov chain is irreducible. For this, consider $x^\circ=(0,0,\mathrm{E_{tot}},0)$ and any fixed state $x= (x_1,x_2,x_3,\mathrm{E_{tot}} - x_3): 0 \leq x_3 \leq \mathrm{E_{tot}}$. Starting at $x^\circ$, by having reaction ${\large \textcircled{\small 5}}$ fire $x_1+\mathrm{E_{tot}}-x_3+x_2$ times in succession, then having reaction ${\large \textcircled{\small 1}}$, immediately followed by reaction ${\large \textcircled{\small 3}}$, fire $x_2$ times in succession and then reaction ${\large \textcircled{\small 1}}$ fire $\mathrm{E_{tot}}-x_3$ times, without any other reactions firing, we see that the Markov chain can transition with positive probability from $x^\circ$ to $x$. Since each reaction is reversible, it also follows that the Markov chain can transition from $x$ to $x^\circ$ with positive probability. Thus, the Markov chain is irreducible.

		Next we will introduce a norm-like function $V$ and show that \eqref{eq:FLcriterion} holds. For each $x \in \X$, let 
		$$V(x)=x_1^2 + ((2 \mathrm{E_{tot}} - 1) b + 1) x_2 + b x_4^2,$$ where
		\begin{equation}
			\label{eq:Ex4.2coef}
			b = \frac{1 + (\kappa_5 +  \kappa_2 \mathrm{E_{tot}} + \kappa_3 \mathrm{E_{tot}}) + \frac{(2 \kappa_5 + \kappa_6 + 2 \kappa_2 \mathrm{E_{tot}})^2}{8\kappa_6}}{\kappa_2 \mathrm{E_{tot}} (2 \mathrm{E_{tot}}-1)}.
		\end{equation}
		Notice that $b>0$ since $\mathrm{E_{tot}} \geq 1$. Then, for each $a \in \R_+$, $\{x \in \X: V(x) \leq a\}$ consists of finitely many states, and for any $x \in \X$,
		\footnotesize
		\begin{eqnarray*}
			Q V(x) &=& \sum_{j=1}^6 \rate_j (x) \cdot (V(x+v_j)-V(x)) \\
			&=& \kappa_1 x_1 x_3 \cdot (((x_1-1)^2 + b (x_4+1)^2) - (x_1^2 + b x_4^2)) \\
			&& \quad + \kappa_2 x_4 \cdot (((x_1+1)^2 + b (x_4-1)^2) - (x_1^2 + b x_4^2)) \\
			&& \quad + \kappa_3 x_4 \cdot ((((2 \mathrm{E_{tot}} - 1) b + 1) (x_2+1) + b (x_4-1)^2) - (((2 \mathrm{E_{tot}} - 1) b + 1) x_2 + b x_4^2)) \\
			&& \quad + \kappa_4 x_2 x_3 \cdot ((((2 \mathrm{E_{tot}} - 1) b + 1) (x_2-1) + b (x_4+1)^2) - (((2 \mathrm{E_{tot}} - 1) b + 1) x_2 + b x_4^2))\\
			&& \quad + \kappa_5 \cdot ((x_1+1)^2 - x_1^2) + \kappa_6 x_1 \cdot ((x_1-1)^2 - x_1^2)\\
			&=& \kappa_1 x_1 x_3 \cdot (-2x_1+1 + b (2x_4+1)) + \kappa_2 x_4 \cdot (2x_1+1 + b (-2x_4+1)) \\
			&& \quad + \kappa_3 x_4 \cdot (((2 \mathrm{E_{tot}} - 1) b + 1) + b (-2x_4+1)) \\
			&& \quad + \kappa_4 x_2 x_3 \cdot ( -((2 \mathrm{E_{tot}} - 1) b + 1) + b (2x_4+1))\\
			&& \quad + \kappa_5 \cdot (2x_1+1) + \kappa_6 x_1 \cdot (-2x_1+1)\\
			&=& - (2 \kappa_1 x_3 + 2 \kappa_6)  \cdot x_1^2 + (2 \kappa_5 + \kappa_6 + \kappa_1 (1 + b) x_3 + 2 b \kappa_1 x_3 x_4 + 2 \kappa_2 x_4) \cdot x_1\\
			&& \quad + \kappa_4 (-2 b \mathrm{E_{tot}} x_3 + 2 b x_3 x_4 + (2b - 1) x_3) \cdot x_2 \\
			&& \quad + \kappa_5 + (\kappa_2 (1+b) + \kappa_3 (2 b \mathrm{E_{tot}} +1)) x_4  - 2 b (\kappa_2 + \kappa_3) x_4^2\\
			&=& - (2 \kappa_1 x_3 + 2 \kappa_6)  \cdot x_1^2 + (2 \kappa_5 + \kappa_6 + \kappa_1 (1 + b) x_3 + 2 b \kappa_1 x_3 x_4 + 2 \kappa_2 x_4) \cdot x_1\\
			&& \quad + \kappa_4 (-2 b x_3^2 + (2b - 1) x_3) \cdot x_2 + \kappa_5 + (\kappa_2 (1+b) + \kappa_3 (2 b \mathrm{E_{tot}} +1)) x_4  - 2 b (\kappa_2 + \kappa_3) x_4^2.
		\end{eqnarray*}
		\normalsize
		The last equality uses the fact that $x_3=\mathrm{E_{tot}}-x_4$. For the following, we note that $-2 b x_3^2 + (2b - 1) x_3 = (-2b x_3 + (2b - 1)) x_3 \leq -1$ when $x_3 \geq 1$.
		
		We now consider two cases for $QV(x)$: when $x_3=0$ and $x_3>0$. For the first case, when $x_3=0$, we have $x_4=\mathrm{E_{tot}}$ and 
		\small
		\begin{equation*}
			Q V(x) = - 2 \kappa_6 \cdot x_1^2 + (2 \kappa_5 + \kappa_6 + 2 \kappa_2 \mathrm{E_{tot}}) \cdot x_1+ \kappa_5 + (\kappa_2 (1 + b) + \kappa_3) \mathrm{E_{tot}} - 2 b \kappa_2 \mathrm{E_{tot}^2}.
		\end{equation*}
		\normalsize
		As a quadratic function, the last expression 
		is bounded above by
		\begin{equation*}
			b \kappa_2 \mathrm{E_{tot}} (1-2 \mathrm{E_{tot}}) + (\kappa_5 +  \kappa_2 \mathrm{E_{tot}} + \kappa_3 \mathrm{E_{tot}}) + \frac{(2 \kappa_5 + \kappa_6 + 2 \kappa_2 \mathrm{E_{tot}})^2}{8\kappa_6} = -1,
		\end{equation*}
		since $b$ is chosen as in \eqref{eq:Ex4.2coef}. For the second case when $x_3 \in \{ 1, \dots, \mathrm{E_{tot}} \}$, we have 
		\begin{eqnarray*}
			Q V(x) &\leq& - 2 \kappa_6 \cdot x_1^2 + (2 \kappa_5 + \kappa_6 + \kappa_1 (1+b) x_3 + 2 b \kappa_1 x_3 x_4 + 2 \kappa_2 x_4) \cdot x_1\\
			&& \quad - \kappa_4 \cdot x_2 + \kappa_5 + (\kappa_2 (1+b) + \kappa_3 (2 b \mathrm{E_{tot}} +1)) x_4  - 2 b (\kappa_2 + \kappa_3) x_4^2\\
			&\leq& - 2 \kappa_6 \cdot x_1^2 + \kappa_6 \cdot x_1^2 + \left( \frac{2 \kappa_5 + \kappa_6 + \kappa_1 (1+b) x_3 + 2 b \kappa_1 x_3 x_4 + 2 \kappa_2 x_4}{2 \sqrt{\kappa_6}} \right)^2\\
			&& \quad - \kappa_4 \cdot x_2 + \kappa_5 + (\kappa_2 (1+b) + \kappa_3 (2 b \mathrm{E_{tot}} +1)) x_4  - 2 b (\kappa_2 + \kappa_3) x_4^2\\
			&\leq& - \kappa_6 \cdot x_1^2 - \kappa_4 \cdot x_2 + \frac{(2 \kappa_5 + \kappa_6 + \kappa_1 (1+b) \mathrm{E_{tot}} + 2 b \kappa_1 \mathrm{E_{tot}^2} + 2 \kappa_2 \mathrm{E_{tot}})^2}{4 \kappa_6}\\
			&& \quad + \kappa_5 + (\kappa_2 (1+b) + \kappa_3 (2 b \mathrm{E_{tot}} +1)) \mathrm{E_{tot}},
		\end{eqnarray*}
		where we have used the fact that $2ab \leq a^2 + b^2$ for the second inequality and the last expression will be less than or equal to $-1$ whenever $x_1$ or $x_2$ is sufficiently large. 
		
		Let $C=\{x \in \X: QV(x)>-1\}$. Then, by the above, $C$ consists of finitely many points, which implies that $C$ is a compact set. Then, \eqref{eq:FLcriterion} holds with $c=1$ and $d = 1+\max_{x \in C} QV(x) \vee 0$. It follows by Theorem \ref{thm:FL} that the Markov chain in Example \ref{exEK2} is non-explosive and positive recurrent, and has a unique stationary distribution $\pi$.
		
		\subsubsection{Application to Example \ref{ex:HistoneModificationAndProtein}}
		\label{appendix:example4.5}
		
		For Example \ref{ex:HistoneModificationAndProtein}, we first check that the Markov chain is irreducible. For this, consider $x^\circ=(0,\Dtot,0)$ and any fixed state $x= (x_1,x_2,x_3): 0 \leq x_1 + x_2 \leq \mathrm{D_{tot}}$. Starting at $x^\circ$, by having reaction ${\large \textcircled{\small 5}}$ fire $x_3$ times, then reaction ${\large \textcircled{\small 3}}$ fire $\Dtot-x_2$ times, and finally, reaction ${\large \textcircled{\small 2}}$ fire $x_1$ times, without any other reactions firing, we see that the Markov chain can transition from $x^\circ$ to $x$ with positive probability. For the reverse transition, by having reaction ${\large \textcircled{\small 6}}$ fire $x_3$ times, then reaction ${\large \textcircled{\small 4}}$ fire $x_1$ times, and finally, reaction ${\large \textcircled{\small 1}}$ fire $\Dtot-x_2$ times, we see that the Markov chain can transition from $x$ to $x^\circ$ with positive probability. Thus, the Markov chain is irreducible.

		Next we will introduce a norm-like function $V$ and show that \eqref{eq:expFLcriterion} holds. For each $x \in \X$, let $$V(x)=x_3.$$
		Then $\{x \in \X: V(x) \leq a\} = \{x \in \Z_+^3: x_1+x_2 \leq \Dtot, x_3 \leq a\}$ consists of finitely many states for each $a \in \R_+$, and for any $x \in \X$,
		\begin{eqnarray*}
			Q V(x) &=& \sum_{j=1}^6 \rate_j (x) \cdot (V(x+v_j)-V(x)) = \sum_{j=1}^6 \rate_j (x) \cdot V(v_j) \\
			&=& \rate_1 (x) \cdot 0 + \rate_2 (x) \cdot 0 + \rate_3 (x) \cdot 0 + \rate_4 (x) \cdot 0+ \kappa_{5a} x_2 \cdot 1 + \kappa_{6a} x_3 \cdot (-1)\\
			&\leq& - \kappa_{6a} x_3 + \kappa_{5a} \Dtot = -c' V(x) + d',
		\end{eqnarray*}
		where $c' = \kappa_{6a}$ and $d' = \kappa_{5a}$. Therefore, we conclude by Theorem \ref{thm:FL} that the Markov chain in Example \ref{ex:HistoneModificationAndProtein} is non-explosive and positive recurrent with a unique stationary distribution $\pi$ such that $\pi(V)<\infty$, and it is exponentially ergodic in the $(V+1)$-norm.
		
		\subsection{Derivation of Markov chain transition directions, $v_j$}
		
		\subsubsection{Example \ref{exEK2}}\label{ex2SI}
		The set of reactions associated to the chemical reaction system in Fig. \ref{fig:EK2}(a) is given by $\Rs=\{(v^-_1,v^+_1),(v^-_2,v^+_2),(v^-_3,v^+_3),(v^-_4,v^+_4),$ $(v^-_5,v^+_5),(v^-_6,v^+_6)\}$, where $(v^-_1,v^+_1),$ $(v^-_2,v^+_2),$ $(v^-_3,v^+_3)$ are defined as in Example \ref{exEK}, and $v^-_4=(0,1,1,0)^T$, $v^+_4= (0,0,0,1)^T$, $v^-_5=v^+_6=(0,0,0,0)^T$, $v^+_5=v^-_6=(1,0,0,0)^T$. Then, the potential transitions of the Markov chain are in six possible directions, $v_j=v^+_j-v^-_j$ for $j=1,...,6$, where $v_1=-v_2=(-1,0,-1,1)^T$, $v_3= - v_4 =(0,1,1,-1)^T$, and $v_5=-v_6=(1,0,0,0)^T$.

		\subsubsection{Example \ref{ex:Braess}}\label{ex3SI}
		The set of five reactions associated to the chemical reaction system in Fig. \ref{fig:braess}(a) is given by $\Rs=\{(v^-_1,v^+_1),(v^-_2,v^+_2),(v^-_3,v^+_3),(v^-_4,v^+_4),(v^-_5,v^+_5)\}$, where $v^-_1=v^-_3=(1,0,0,0)^T$, $v^+_1=v^-_2=v^-_5=(0,1,0,0)^T$, $v^+_3=v^-_4=v^+_5=(0,0,1,0)^T$, $v^+_2=v^+_4=(0,0,0,1)^T$. Then, the potential transitions of the Markov chain are in five possible directions, $v_j=v^+_j-v^-_j$, $j=1,...,5$, where $v_1=(-1,1,0,0)^T$, $v_2=(0,-1,0,1)^T$, $v_3=(-1,0,1,0)^T$, $v_4=(0,0,-1,1)^T$ and $v_5=(0,-1,1,0)^T$.

		\subsection{A generalization of Theorem \ref{thm:SumTheorem2}}
		
		Here, we provide a more general version of Theorem \ref{thm:SumTheorem2}. The simpler form given as Theorem \ref{thm:SumTheorem2} in the main text was used there because it is more straightforward to state and the conditions are easier to verify. However, the more general version of the theorem provided in this section can be useful in some cases, such as Example
		\ref{ex:Braess} (see Section \ref{appendix:A22}). 
		
		The generalization relies on the idea that grouping of vectors can be more general than what is described in \eqref{eq:sets_of_dot_product+-13}, and so we introduce the following assumption. 
		
		\begin{assumption}
			\label{as:SumTheorem3}
			Consider a collection of distinct vectors $v_1,\ldots,v_n$ in $\Z^d \setminus \{0\}$ and a matrix $A \in \Z^{m \times d}$ with non-zero rows. Suppose that there exists a partition  \footnote{In particular, $\{G^1, \dots, G^s\}$ is a finite collection of non-empty disjoint sets of distinct numbers, the union of which is $\{1,\dots,n\}$.} $\{G^1, \dots, G^s\}$ of $\{1,\dots,n\}$ and an associated bijection $\sigma:\{1,\ldots,n\} \longrightarrow \{1,\ldots,n\}$ such that for $p_0:=0$ and $1 \leq k \leq s$, with $p_k:=\sum_{\ell=1}^k |G^\ell|$, we have $\sigma(q) \in G^k$ for $p_{k-1}+1 \leq q \leq p_k$, 
			and whenever $p_{k-1}+2 \leq q \leq p_k$, we also have for each $1 \leq i \leq m$ that either $\inn{A_{i\bullet},v_{\sigma(q)}}$ is equal to $\inn{A_{i\bullet},v_{\sigma(q-1)}}$ or it is 0.
		\end{assumption}
		
		\begin{theorem}
			\label{thm:SumTheorem3}
			Consider a non-empty set $\X \subseteq \Z_+^d$, suppose that Assumption \ref{as:SumTheorem3} holds and consider two collections of non-negative functions on $\X$, $\rate=(\rate_1, \dots,\rate_n)$ and  $\Breve{\rate}= (\Breve{\rate}_1, \dots,\Breve{\rate}_n)$, such that \eqref{eq:LambdaNotOutofSpace} holds and the associated continuous-time Markov chains do not explode in finite time. Further suppose that both of the following conditions hold:
			\begin{enumerate}
				\item[(i)]
				For each $1 \leq j \leq n$, the vector $A v_j$ has entries in $\{-1,0,1\}$ only.
				\item[(ii)]
				For each $x \in \X$, $1 \leq i \leq m$ and $y \in \partial_i(K_A+x) \cap \X$ we have that for each $1 \le k \le s$,
				\begin{equation*}
					\sum_{j \in G^{k,-}_i} \Breve{\rate}_j (y) \leq \sum_{j \in G^{k,-}_i} \rate_j(x), \quad \text{ where } G^{k,-}_i = \{j \in G^k \:|\: \inn{A_{i\bullet},v_j} = -1 \},
				\end{equation*}
				and
				\begin{equation*}
					\sum_{j \in G^{k,+}_i} \Breve{\rate}_j(y) \geq \sum_{j \in G^{k,+}_i} \rate_j(x), \quad \text{ where } G^{k,+}_i = \{j \in G^k \:|\: \inn{A_{i\bullet},v_j} = 1 \}.
				\end{equation*}
			\end{enumerate}
			Then, for each pair $\initialx,\initialxbreve \in \X$ such that $\initialx \preccurlyeq_A \initialxbreve$, there exists a probability space $(\Omega,\F,\PP)$ with two continuous-time Markov chains $X = \{X(t), \: t \geq 0\}$ and $\Breve{X}=\{\Breve{X}(t), \: t \geq 0\}$ defined there, each having state space $\X \subseteq \Z^d_+$, with infinitesimal generators $Q$ and $\Breve{Q}$, associated with $\rate$ and $\Breve{\rate}$ respectively, with initial conditions $X(0)=\initialx$ and $\Breve{X}(0)=\initialxbreve$ and such that:
			\begin{equation}
				\PP\left[X(t) \preccurlyeq_A \Breve{X}(t) \text{ for every } t \geq 0 \right]=1.   
			\end{equation}
		\end{theorem}
		
		\begin{remark}
			Under Assumption \ref{as:SumTheorem3}, for a given $1 \leq k \leq s$ and $1 \leq i \leq m$, at most one of $G^{k,-}_i$ and $G^{k,+}_i$ is non-empty. If $G^{k,-}_i \neq \emptyset$, then $G^{k,-}_i = \{\sigma(q) \:|\: p_{k-1}+1 \leq q \leq q^*\}$ where $q^* = \max \{q \:|\: p_{k-1}+1 \leq q \leq p_k \text{ and } \inn{A_{i\bullet},v_{\sigma(q)}} = -1 \}$. On the other hand, if $G^{k,+}_i \neq \emptyset$, then $G^{k,+}_i = \{\sigma(q) \:|\: p_{k-1}+1 \leq q \leq q^*\}$ where $q^* = \max \{q \:|\: p_{k-1}+1 \leq q \leq p_k \text{ and } \inn{A_{i\bullet},v_{\sigma(q)}} = 1 \}$. Furthermore, in either case, for $q^* < q \leq p_k$, we have $\inn{A_{i\bullet},v_{\sigma(q)}} = 0$.
		\end{remark}
		
		The proof for Theorem \ref{thm:SumTheorem3} can be found in Section \ref{appendix:A21}. With Theorem \ref{thm:SumTheorem3} in place, we can extend Theorems \ref{thm:comparison_MFPT} and \ref{comparStationaryDistributionMAIN} by adding an additional alternative condition $(iv)$: Assumption \ref{as:SumTheorem3} holds, and conditions $(i)$ and $(ii)$ in Theorem \ref{thm:SumTheorem3} are satisfied. If this is satisfied instead of one of $(i)$ -- $(iii)$ in Theorems \ref{thm:comparison_MFPT} and \ref{comparStationaryDistributionMAIN}, then the conclusions of these theorems about (mean) first passage times and stationary distributions will still hold.
		
		\subsubsection{Proof of Theorem \ref{thm:SumTheorem3}}
		\label{appendix:A21}
		
		We initially assume that $\sup_{x \in \X} \rate_j(x) < \infty$ and $\sup_{x \in \X} \Breve{\rate}_j(x) < \infty$ for every $1 \leq j \leq n$, and let $\lambda > 0$ such that \eqref{eq:LambdaCondition} holds. We shall relax these assumptions later. Further suppose that Assumption \ref{as:SumTheorem3} and condition (i) of Theorem \ref{thm:SumTheorem3} both hold. For $x \in \X$, define $I^k(x)$, $I^k_q(x)$, $\PPhi_{\lambda}$, $\Breve{I}^k(x)$, $\Breve{I}^k_q(x)$ and $\Breve{\PPhi}_{\lambda}$ in the same manner as in the proof of Theorem \ref{thm:SumTheorem2} (see \eqref{Int1} -- \eqref{eq:DefinitionPhiLambdaBIS2A}), with $\{G^k \:|\: 1 \leq k \leq s\}$, $\{p_k \:|\: 0 \leq k \leq s\}$ and $\sigma$ as in Assumption \ref{as:SumTheorem3}. Our proof of Theorem \ref{thm:SumTheorem3} has some elements that are the same as those for the proof of Theorem \ref{thm:SumTheorem2}. However, some additional elements are needed. We give the details for completeness. As for Theorem \ref{thm:SumTheorem2}, $\PPhi_{\lambda}(\cdot,\cdot)$ and $\Breve{\PPhi}_{\lambda}(\cdot,\cdot)$ are well-defined as $\X$-valued functions.
		
		\begin{lemma}
			\label{lem:PhiLambda_ineqSumTheorem} 
			Suppose that $x,y \in \X$ are such that $x \preccurlyeq_A y$ and the following hold: whenever $y \in \partial_i(K_A+x) \cap \X$ for some $1 \leq i \leq m$, we have that for each $1 \le k \le s$,
			\begin{equation}
				\label{eq:CouplingConditionSumTheoremProof1}
				\sum_{j \in G^{k,-}_i} \Breve{\rate}_j (y) \leq \sum_{j \in G^{k,-}_i} \rate_j(x), \quad \text{ where } G^{k,-}_i = \{j \in G^k \:|\: \inn{A_{i\bullet},v_j} = -1 \},
			\end{equation}
			and
			\begin{equation}
				\label{eq:CouplingConditionSumTheoremProof2}
				\sum_{j \in G^{k,+}_i} \Breve{\rate}_j(y) \geq \sum_{j \in G^{k,+}_i} \rate_j(x), \quad \text{ where } G^{k,+}_i = \{j \in G^k \:|\: \inn{A_{i\bullet},v_j} = 1 \}.
			\end{equation}
			Then, for each $u \in [0,1]$:
			\begin{equation}
				\label{eq:ClaimPhi_Lambda_SumTheorem}
				\PPhi_{\lambda}(x,u) \preccurlyeq_A  \Breve{\PPhi}_{\lambda}(y,u).
			\end{equation}
		\end{lemma}
		
		\begin{proof}
			First, we note that $\PPhi_{\lambda}, \Breve{\PPhi}_{\lambda}$ have the following properties: for every $u \in [0,1]$, $1\le k \le s$, $j \in G^k$, 
			\begin{equation}
				\text{ if } \PPhi_{\lambda}(x,u) = x + v_j \text{, then } \Breve{\PPhi}_{\lambda}(y,u) \in \{y + v_\ell: \ell \in G^k \} \cup \{y\},
			\end{equation}
			since $I_{\sigma^{-1}(j)}^k(x), \Breve{I}_{\sigma^{-1}(\ell)}^k(y) \subseteq [\frac{p_{k-1}}{n},\frac{p_{k}}{n})$ for $\ell \in G^k$. Similarly,
			\begin{equation}
				\label{eq:PhiLambdaCondition1SumTheorem}
				\text{ if } \Breve{\PPhi}_{\lambda}(y,u) = y +v_j \text{, then } \PPhi_{\lambda}(x,u) \in \{x + v_\ell: \ell \in G^k\} \cup \{x\}.
			\end{equation}
			
			Furthermore, for $1 \leq k \leq s$, $1 \leq i \leq m$, $j \in G^{k,+}_i$, if $\sum_{\ell \in G^{k,+}_i} \Breve{\rate}_\ell (y) \geq \sum_{\ell \in G^{k,+}_i} \rate_\ell (x)$, then 
			\begin{equation}
				\label{eq:PhiLambdaCondition3SumTheorem}
				\PPhi_{\lambda}(x,u) = x +v_j \text{ implies that } \Breve{\PPhi}_{\lambda}(y,u) =y +v_\ell \: \text{ for some } \: \ell \in G^{k,+}_i,
			\end{equation}
			since under the condition, we have $\cup_{q=p_{k-1}+1}^{q^*} I^k_q (x) \subseteq \cup_{q=p_{k-1}+1}^{q^*} \Breve{I}^k_q (y)$ where $q^* = \max \{q \:|\: p_{k-1}+1 \leq q \leq p_k \text{ and } \inn{A_{i\bullet},v_{\sigma(q)}} = 1 \}$ and, by Assumption \ref{as:SumTheorem3}, $G^{k,+}_i = \{\sigma(q) \:|\: p_{k-1}+1 \leq q \leq q^*\}$.
			Similarly, for $1 \leq k \leq s$, $1 \leq i \leq m$, $j \in G^{k,-}_i$, if $\sum_{\ell \in G^{k,-}_i} \Breve{\rate}_\ell (y) \leq \sum_{\ell \in G^{k,-}_i} \rate_\ell (x)$, then
			\begin{equation}
				\label{eq:PhiLambdaCondition4SumTheorem}
				\Breve{\PPhi}_{\lambda}(y,u) = y +v_j \text{ implies that } \PPhi_{\lambda}(x,u) =x +v_\ell \: \text{ for some } \: \ell \in G^{k,-}_i.
			\end{equation}
			
			We also have that, for $1 \le k \le s$ and $j \in G^k$, $x \preccurlyeq_A y + v_j$ if and only if
			\begin{equation}
				\label{eq:RelationAdotted_yplusvjSumTheorem}
				\inn{A_{i\bullet},y -x} + \inn{A_{i\bullet},v_j} \geq 0, \qquad \text{for every } 1 \leq i \leq m.
			\end{equation}
			Similarly, $x + v_j \preccurlyeq_A y$ if and only if
			\begin{equation}
				\label{eq:RelationAdotted_xplusvjSumTheorem}
				\inn{A_{i\bullet},y -x} - \inn{A_{i\bullet},v_j} \geq 0, \qquad \text{for every } 1 \leq i \leq m.
			\end{equation}
			
			Furthermore, for $1 \le k \le s$ and $j,\ell \in G^k$, $x + v_\ell \preccurlyeq_A y  + v_j$ if and only if
			\begin{equation}
				\label{eq:RelationAdotted_xplusvjyplusvjSumTheorem}
				\inn{A_{i\bullet},y -x} + \inn{A_{i\bullet},v_j-v_\ell} \geq 0, \qquad \text{for every } 1 \leq i \leq m.
			\end{equation}
			
			To prove \eqref{eq:ClaimPhi_Lambda_SumTheorem}, we first consider the situation where $y \in \interior(K_A + x) = \{ w \in \R^d \:|\: Ax < Aw\}$. Then, for each $1 \leq i \leq m$, $\inn{A_{i\bullet},y-x} >0$ and since $A \in \Z^{m \times d}$ and $y-x \in \Z^d$, we have $\inn{A_{i\bullet},y-x} \geq 1$. This implies that for each $1 \le k \le s$ and $j \in G^k$,
			\begin{equation}
				\label{positivecondSumTheorem1}
				\inn{A_{i\bullet},y-x} + \inn{A_{i\bullet},v_j} \geq 1 + \inn{A_{i\bullet},v_j}\ge 0, \qquad \text{for every } 1 \leq i \leq m,
			\end{equation}
			since $\inn{A_{i\bullet},v_j} \in \{-1,0,1\}$ by condition $(i)$ of Theorem \ref{thm:SumTheorem3}. Similarly, for each $1 \le k \le s$ and $j \in G^k$,
			\begin{equation}
				\label{positivecondSumTheorem2}
				\inn{A_{i\bullet},y-x} - \inn{A_{i\bullet},v_j} \geq 1 - \inn{A_{i\bullet},v_j} \geq 0, \qquad \text{for every } 1 \leq i \leq m.
			\end{equation}
			
			In addition, for $1 \leq k \leq s$ and $j,\ell \in G^k$,
			\begin{equation}
				\label{positivecondSumTheorem3}
				\inn{A_{i\bullet},y -x} + \inn{A_{i\bullet},v_j-v_\ell} \geq 1 + \inn{A_{i\bullet},v_j-v_\ell} \geq 0, \quad \text{ for every } 1 \leq i \leq m,
			\end{equation}
			since, by Assumption \ref{as:SumTheorem3}, if $\inn{A_{i\bullet},v_j} \neq 0$, then either $\inn{A_{i\bullet},v_\ell}=\inn{A_{i\bullet},v_j}$ or $\inn{A_{i\bullet},v_\ell}=0$.
			It follows from \eqref{positivecondSumTheorem1} -- \eqref{positivecondSumTheorem3} that
			if $y \in \interior(K_A + x) \cap \X$, then for any $1 \leq k \leq s$ and $j, \ell \in G^k$:
			\begin{equation}
				\label{eq:JumpsInteriorSumTheorem}
				x \preccurlyeq_A y + v_j, \: x + v_j \preccurlyeq_A y \: \text{ and } \: x + v_\ell \preccurlyeq_A y + v_j.
			\end{equation}
			
			We also have, by assumption, that $x  \preccurlyeq_A y$. It follows that if $y \in \interior(K_A + x) \cap \X$, then $\{x,x +v_\ell \:|\: \ell \in G^k\}\preccurlyeq_A \{y,y+v_j \:|\: j \in G^k\}$ for $1 \le k \le s$ and consequently \eqref{eq:ClaimPhi_Lambda_SumTheorem} holds for all $u \in [0,1]$.
			
			Now, we turn to the other situation where $y \in \partial_i(K_A+x) \cap \X$ for some $1 \leq i \leq m$. Then $\Kb_{y} := \{ i\:|\: \inn{A_{i\bullet},y} = \inn{A_{i\bullet},x}, 1\le i\le m \}$ is non-empty. Let $u\in[0,1]$. We consider two cases.
			
			\medskip
			\textbf{Case 1:} $\Breve{\PPhi}_{\lambda}(y,u) = y +v_j$ for some $1 \leq j \leq n$.
			
			\medskip
			Fix such an index $j$. Consider the unique $1 \leq k \leq s$ such that $j \in G^k$. Then, by \eqref{eq:PhiLambdaCondition1SumTheorem}, either $\PPhi_{\lambda}(x,u) = x + v_\ell$ for some $\ell \in G^k$, or $\PPhi_{\lambda}(x,u) = x$.
			\begin{enumerate}
				\item[a)]
				Suppose $\PPhi_{\lambda}(x,u) = x + v_\ell$ for some $\ell \in G^k$. Observe that for every $i \notin \Kb_y$, $\inn{A_{i\bullet},y-x} > 0$ and as for \eqref{positivecondSumTheorem3}, $\inn{A_{i\bullet},(y+v_j)-(x+v_\ell)} \geq 0$, while for $i \in \Kb_y$, $\inn{A_{i\bullet},(y+v_j)-(x+v_\ell)} = \inn{A_{i\bullet},v_j}- \inn{A_{i\bullet},v_\ell}$. For each $i \in \Kb_y$,
				
				\begin{enumerate}
					\item[i)]
					if $\inn{A_{i\bullet},v_j}=-1$ and $\inn{A_{i\bullet},v_\ell}=0$, then $j \in G^{k,-}_i$ and $\ell \notin G^{k,-}_i$. By \eqref{eq:CouplingConditionSumTheoremProof1}, we would then have $\sum_{r \in G^{k,-}_i} \Breve{\rate}_r (y) \leq \sum_{r \in G^{k,-}_i} \rate_r (x)$, which would imply by \eqref{eq:PhiLambdaCondition4SumTheorem} that $\PPhi_{\lambda}(x,u) =x +v_r$ for some $r \in G^{k,-}_i$. Since we are assuming that $\PPhi_{\lambda}(x,u) = x + v_{\l}$ and we know the vectors $v_1,\ldots,v_n$ are distinct, we obtain that $\l=r \in G_i^{k,-}$. This contradicts $\ell \notin G^{k,-}_i$.
					\item[ii)]
					if $\inn{A_{i\bullet},v_j}=0$ and $\inn{A_{i\bullet},v_\ell}=1$, then $j \notin G^{k,+}_i$ and $\ell \in G^{k,+}_i$. By \eqref{eq:CouplingConditionSumTheoremProof2}, we would then have $\sum_{r \in G^{k,+}_i} \Breve{\rate}_r (y) \geq \sum_{r \in G^{k,+}_i} \rate_r (x)$, which would imply by \eqref{eq:PhiLambdaCondition3SumTheorem} that $j \in G^{k,+}_i$. This contradicts $j \notin G^{k,+}_i$.
					\item[iii)]
					in all the other cases, that is when $$(\inn{A_{i\bullet},v_j},\inn{A_{i\bullet},v_\ell}) \in \{(1,0),(0,-1),(1,1),(0,0),(-1,-1)\},$$ we have $\inn{A_{i\bullet},(y+v_j)-(x+v_\ell)} \geq 0.$
				\end{enumerate}

				Combining the above for case a), we see that $\inn{A_{i\bullet},(y+v_j)-(x+v_\ell)} \geq 0$ for each $1 \le i \le m$, which implies that $\PPhi_{\lambda}(x,u) = x+v_\ell \preccurlyeq_A y+v_j= \Breve{\PPhi}_{\lambda}(y,u)$.
				
				\item[b)]
				Suppose $\PPhi_{\lambda}(x,u) = x$. We claim that $y + v_j \in K_A + x$. To see this, observe that for every $i \notin \Kb_y$, $\inn{A_{i\bullet},y-x} > 0$ and as for \eqref{positivecondSumTheorem1}, $\inn{A_{i\bullet},(y+v_j)-x} \geq 0$, while for $i \in \Kb_y$, $\inn{A_{i\bullet},(y+v_j)-x} = \inn{A_{i\bullet},v_j}\in \{-1,0,1\}$. For each $i \in \Kb_y$, if $\inn{A_{i\bullet},v_j}=-1$, which means $j \in G^{k,-}_i$, then by \eqref{eq:CouplingConditionSumTheoremProof1} we would have $\sum_{\ell \in G^{k,-}_i} \Breve{\rate}_\ell (y) \leq \sum_{\ell \in G^{k,-}_i} \rate_\ell (x)$, which would imply by \eqref{eq:PhiLambdaCondition4SumTheorem} that $\PPhi_{\lambda}(x,u) = x +v_\ell$ for some $\ell \in G^{k,-}_i$, but this contradicts the assumption that $\PPhi_{\lambda}(x,u) = x$. So we must have $\inn{A_{i\bullet},v_j} \geq 0$ and hence $\inn{A_{i\bullet},(y+v_j)-x} \geq 0$ for all $i \in \Kb_y$. Thus, $y + v_j \in K_A + x$ and so $\PPhi_{\lambda}(x,u) = x \preccurlyeq_A  y + v_j = \Breve{\PPhi}_{\lambda}(y,u)$.
			\end{enumerate}
			
			\medskip
			\textbf{Case 2:} $\Breve{\PPhi}_{\lambda}(y,u) = y$. Again, we consider two subcases.
			\begin{enumerate}
				\item[a)]
				If $\PPhi_{\lambda}(x,u)=x$, then \eqref{eq:ClaimPhi_Lambda_SumTheorem} holds, because $x \preccurlyeq_A y$.
				\item[b)]
				If $\PPhi_{\lambda}(x,u) = x +v_j$ for some $1 \leq j \leq n$, we claim that $y \in K_A + x + v_j$ for the corresponding value of $j$. To see this, fix the value of $j$ for which $\PPhi_{\lambda}(x,u) = x +v_j$, let $1 \leq k \leq s$ be such that $j \in G^k$, and observe that for every $i \notin \Kb_y$, $\inn{A_{i\bullet},y-x} > 0$ and as for \eqref{positivecondSumTheorem2}, $\inn{A_{i\bullet},y-(x+v_j)} \geq 0$, while for $i \in \Kb_y$, $\inn{A_{i\bullet},y-(x+v_j)} = - \inn{A_{i\bullet},v_j}\in \{-1,0,1\}$. For each $i \in \Kb_y$, if $\inn{A_{i\bullet},v_j}=1$, which means $j \in G^{k,+}_i$, then by \eqref{eq:CouplingConditionSumTheoremProof2}, we would have $\sum_{\ell \in G^{k,+}_i} \Breve{\rate}_\ell (y) \geq \sum_{\ell \in G^{k,+}_i} \rate_\ell (x)$, which would imply by \eqref{eq:PhiLambdaCondition3SumTheorem} that $\Breve{\PPhi}_{\lambda}(y,u) = y +v_\ell$ for some $\ell \in G^{k,+}_i$, but this contradicts the assumption that $\Breve{\PPhi}_{\lambda}(y,u) = y$. So we must have $\inn{A_{i\bullet},v_j} \leq 0$ and hence $\inn{A_{i\bullet},y-(x+v_j)} \geq 0$ for all $i \in \Kb_y$. Thus, we have $y \in K_A + x  + v_j$ and then $\PPhi_{\lambda}(x,u) = x  + v_j \preccurlyeq_A  y = \Breve{\PPhi}_{\lambda}(y,u)$.
			\end{enumerate}  
			\vskip -0.1in    
		\end{proof}
		
		In order to prove Theorem \ref{thm:SumTheorem3}, from here on we can follow a similar procedure to the one used in the proof of Theorem \ref{thm:MainResult} after Lemma \ref{lem:PhiLambda_ineq} was proved there. For the case where \eqref{eq:uniformization_condition_intensity} holds, we define two discrete-time processes, $Y=(Y_k)_{k \geq 0}$ and $\Breve{Y}=(\Breve{Y}_k)_{k \geq 0}$, by defining $Y_0 := \initialx$, $\Breve{Y}_0 := \initialxbreve$, and for $k \geq 0$,
		\begin{equation}
			Y_{k+1} := \PPhi_{\lambda}(Y_k,U_{k+1}), \qquad  \Breve{Y}_{k+1} := \Breve{\PPhi}_{\lambda}(\Breve{Y}_k,U_{k+1}),
		\end{equation}
		and define $X$ and $\Breve{X}$ using these and an independent Poisson process $N$ as in \eqref{def:X_and_XbreveMAIN}. For the case where \eqref{eq:uniformization_condition_intensity} does not hold, we can use a truncation procedure similar to that for Theorem \ref{thm:MainResult}. In both cases, we use Lemma \ref{lem:PhiLambda_ineqSumTheorem} instead of Lemma \ref{lem:PhiLambda_ineq}.
		
		\subsubsection{Two other $A$ matrices for Example \ref{ex:Braess}}
		\label{appendix:A22}
		
		Let
		\begin{equation}\label{matrixAbraess2}
			A= \begin{bmatrix}
				-1 & 0 & 0 & 0\\
				0 & 0 & -1 & 0\\
				0 & -1 & -1 & 0
			\end{bmatrix}.
		\end{equation}
		For $x \in \X$, consider infinitesimal transition rates $\Breve{\rate}_1(x),\Breve{\rate}_2(x), \Breve{\rate}_3(x), \Breve{\rate}_4(x)$ and $\Breve{\rate}_5(x)$ defined as for $\rate_1(x),\rate_2(x),\rate_3(x),\rate_4(x)$ and $\rate_5(x)$ in (\ref{ratesBraess}), but with $\Breve{\kappa}_i$ in place of $\kappa_i$ where $\Breve{\kappa}_i=\kappa_i$, for $i = 1,2,3,4$, and $\Breve{\kappa}_5 \leq \kappa_5$. Suppose that $\kappa_2 > \kappa_4$. Now, let us verify that the assumptions of Theorem \ref{thm:SumTheorem3} hold. Condition $(i)$ holds since $Av_1=(1,0,-1)^T$, $Av_2=(0,0,1)^T$, $Av_3=(1,-1,-1)^T$, $Av_4=(0,1,1)^T$ and $Av_5=(0,-1,0)^T$. Assumption \ref{as:SumTheorem3} holds with $G^1 = \{3,1\}$, $G^2 = \{4,2\}$, $G^3 = \{5\}$ and $\sigma(1)=3$, $\sigma(2)=1$, $\sigma(3)=4$, $\sigma(4)=2$, $\sigma(5)=5$. To verify that condition $(ii)$ of Theorem \ref{thm:SumTheorem3} holds, fix $x\in \X$ and first consider $y\in \partial_1(K_A+x) \cap \X $, where $\partial_1(K_A+x) \cap \X = \{ w \in \X \:|\: x_1 = w_1, x_3 \geq w_3, x_2 + x_3 \geq w_2 + w_3, x_4 \le w_4\}$. Given that $\inn{A_{1\bullet},v_1}=\inn{A_{1\bullet},v_3}=1$, we need to check that $\rate_1(x) + \rate_3(x) \le \Breve{\rate}_1(y) + \Breve{\rate}_3(y)$. Since $y \in \partial_1(K_A+x) \cap \X$, then $\rate_1(x) = \kappa_1 x_1 =\kappa_1 y_1 = \Breve{\kappa}_1 y_1 = \Breve{\rate}_1(y)$ and $\rate_3(x) = \kappa_3 x_1 =\kappa_3 y_1 = \Breve{\kappa}_3 y_1 = \Breve{\rate}_3(y)$, and so the desired inequality holds with equality. Secondly, consider $y \in \partial_2(K_A+x) \cap \X = \{ w \in \X \:|\: x_1 \geq w_1, x_3 = w_3, x_2 \geq w_2, x_4\le w_4\}$. Given that $\inn{A_{2\bullet},v_3}=\inn{A_{2\bullet},v_5}=-1$ and $\inn{A_{2\bullet},v_4}=1$, we need to check that $\rate_3(x) \ge \Breve{\rate}_3(y)$, $\rate_4(x) \le \Breve{\rate}_4(y)$ and $\rate_5(x) \ge \Breve{\rate}_5(y)$. Since $y \in \partial_2(K_A+x) \cap \X$, then $\rate_3(x) = \kappa_3 x_1 \geq \kappa_3 y_1 = \Breve{\kappa}_3 y_1 = \Breve{\rate}_3(y)$, $\rate_4(x) = \kappa_4 x_3 = \kappa_4 y_3 = \Breve{\kappa}_4 y_3 = \Breve{\rate}_4(y)$ and $\rate_5(x) = \kappa_5 x_2 \geq \kappa_5 y_2 \geq \Breve{\kappa}_5 y_2 = \Breve{\rate}_5(y)$, and so the desired inequality holds. Lastly, consider $y \in \partial_3(K_A+x) \cap \X = \{ w \in \X \:|\: x_1 \geq w_1, x_3 \geq w_3, x_2 + x_3 = w_2 + w_3, x_4\le w_4\}$. Given that $\inn{A_{1\bullet},v_1}=\inn{A_{1\bullet},v_3}=-1$ and $\inn{A_{1\bullet},v_2}=\inn{A_{1\bullet},v_4}=1$, we need to check that $\rate_1(x) +\rate_3(x) \ge \Breve{\rate}_1(y) + \Breve{\rate}_3(y)$ and $\rate_2(x) + \rate_4(x) \le \Breve{\rate}_2(y) + \Breve{\rate}_4(y)$. For $y\in \partial_3(K_A+x) \cap \X $, since $\kappa_2 > \kappa_4$ was assumed, we have that $\rate_2(x) + \rate_4(x) =  \kappa_2 x_2 + \kappa_4 x_3 = (\kappa_2-\kappa_4) x_2 + \kappa_4 (x_2 + x_3) \le (\kappa_2-\kappa_4) y_2 + \kappa_4 (y_2 + y_3) = \Breve{\kappa}_2 y_2 + \Breve{\kappa}_4 y_3 = \Breve{\rate}_2(y) + \Breve{\rate}_4(y)$ and $\rate_1(x) = \kappa_1 x_1 \ge \kappa_1 y_1 = \Breve{\kappa}_1 y_1  = \Breve{\rate}_1(y)$, $\rate_3(x) = \kappa_3 x_1 \ge \kappa_3 y_1 = \Breve{\kappa}_3 y_1 = \Breve{\rate}_3(y)$. Thus, the conditions of Theorem \ref{thm:SumTheorem3} are satisfied and so the conclusion of that theorem holds.
		
\begin{figure}[t!]
			\centering
			\includegraphics[scale=0.66]{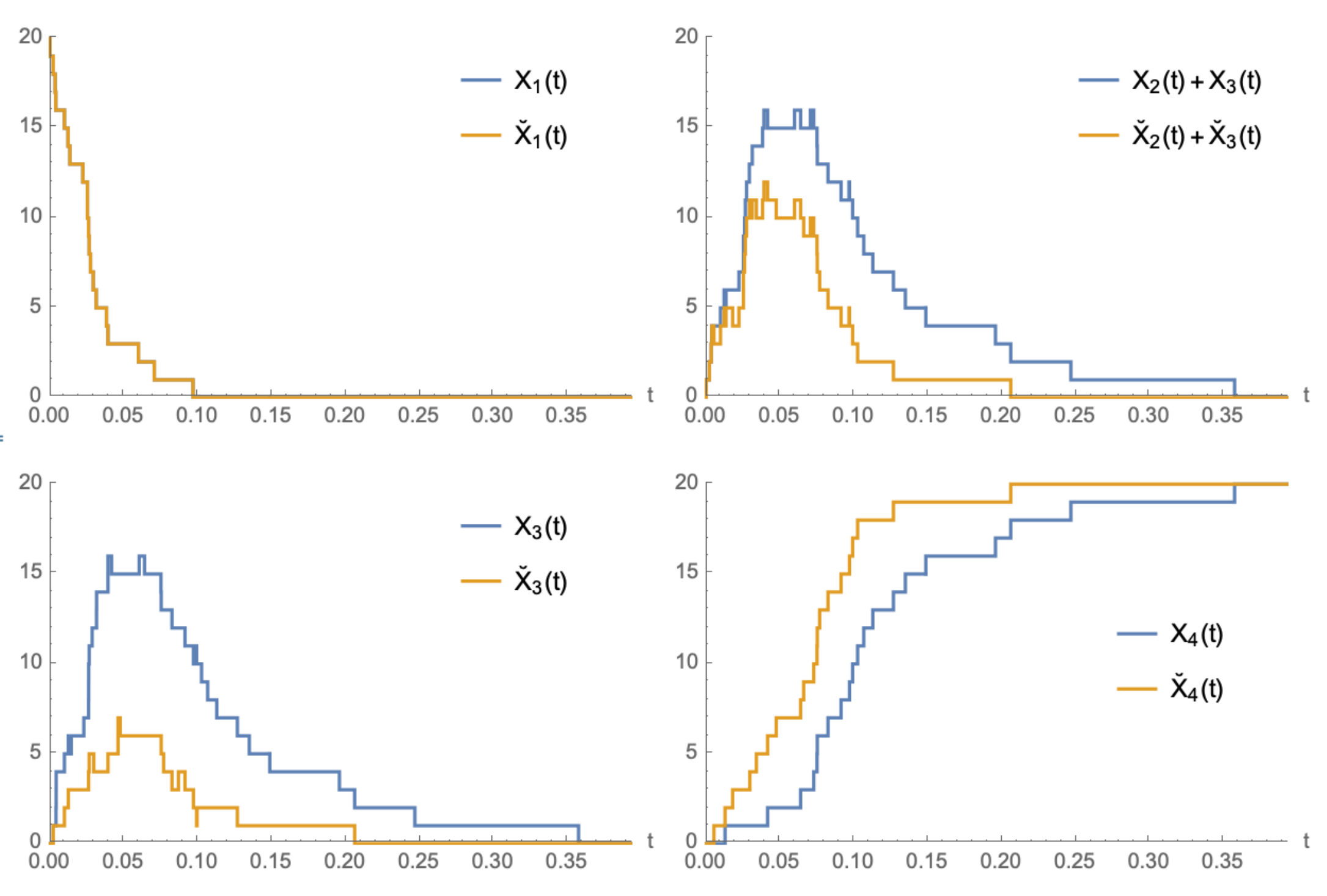}
			\caption{ \small { \bf A typical coupled realization of sample paths for Example \ref{ex:Braess}.} Here $\mathrm{S_{tot}}=20$, $\kappa_1=\Breve{\kappa}_1=30$, $\kappa_2=\Breve{\kappa}_2=50$,
				$\kappa_3=\Breve{\kappa}_3=10$,
				$\kappa_4=\Breve{\kappa}_4=10$,
				$\kappa_5=1000$ and $\Breve{\kappa}_5=10$. Both processes $X$ and $\Breve{X}$ start at $(\mathrm{S_{tot}},0,0,0)$. As shown in Section \ref{appendix:A22}, for these parameters, we have, almost surely, $X(t) \preccurlyeq_A \Breve{X}(t) \text{ for every } t \geq 0$ where the matrix $A$ is given in \eqref{matrixAbraess2}. A coupled realization using the algorithm described in Section \ref{sec:AlgorithmStochasticSimulation} is plotted to illustrate this result. In particular, we see in this sample that $X_1 (t) \geq \Breve{X}_1 (t)$, $X_2 (t) + X_3 (t) \geq \Breve{X}_2 (t) + \Breve{X}_3 (t)$ and $X_3 (t) \geq \Breve{X}_3 (t)$ for all times $t$. Moreover, the first passage time to $(0,0,0,\mathrm{S_{tot}})$ for $X$ (which is equivalent to the first time to get to the state where $X_4 = \mathrm{S_{tot}}$) is larger than for $\Breve{X}$. Since this first passage time result is true for all coupled samples, we can conclude that the mean first passage time from $(\mathrm{S_{tot}},0,0,0)$ to $(0,0,0,\mathrm{S_{tot}})$ is larger for $X$ than for $\Breve{X}$.
			}
			\label{fig:SPR}
		\end{figure}

		Let $\Gamma = \{(0,0,0,\mathrm{S_{tot}})\}$. This is an increasing set in $\X$ with respect to the relation $\preccurlyeq_A$. Let $T_{(0,0,0,\mathrm{S_{tot}})}$, respectively $\Breve{T}_{(0,0,0,\mathrm{S_{tot}})}$ be the first time that the Markov chain $X$, respectively $\Breve{X}$, reaches the set $\Gamma$. Then, by the generalization of Theorem \ref{thm:comparison_MFPT}, if $X(0)=\Breve{X}(0)=(\mathrm{S_{tot}},0,0,0)$, we have that $\Breve{T}_{(0,0,0,\mathrm{S_{tot}})} \preccurlyeq_{st} T_{(0,0,0,\mathrm{S_{tot}})}$. It follows that increasing $\kappa_5$ will increase the mean first passage time from $(\mathrm{S_{tot}},0,0,0)$ to $(0,0,0,\mathrm{S_{tot}})$ when $\kappa_2 > \kappa_4$ (See Figure \ref{fig:SPR}). Indeed, when $\kappa_2 > \kappa_4$, it takes a longer time to get to $(0,0,0,\mathrm{S_{tot}})$ from $(\mathrm{S_{tot}},0,0,0)$ if reaction ${\large \textcircled{\small 5}}$ is added to the system without that reaction.
		
		On the other hand, suppose
		\begin{equation*}\label{matrixAbraess3}
			A= \begin{bmatrix}
				-1 & 0 & 0 & 0\\
				0 & -1 & 0 & 0\\
				0 & -1 & -1 & 0
			\end{bmatrix} 
		\end{equation*}
		and infinitesimal transition rates $\Breve{\rate}_1(x),\Breve{\rate}_2(x), \Breve{\rate}_3(x), \Breve{\rate}_4(x)$ and $\Breve{\rate}_5(x)$ are defined as for $\rate_1(x)$,$\rate_2(x)$, $\rate_3(x)$,$\rate_4(x)$ and $\rate_5(x)$ in (\ref{ratesBraess}), but with $\Breve{\kappa}_i$ in place of $\kappa_i$ where $\Breve{\kappa}_i=\kappa_i$, for $i = 1,2,3,4$, $\Breve{\kappa}_5 \geq \kappa_5$, and $\kappa_2 < \kappa_4$. We can verify that the assumptions of Theorem \ref{thm:SumTheorem3} hold, as follows. Condition $(i)$ holds since $Av_1=(1,-1,-1)^T$, $Av_2=(0,1,1)^T$, $Av_3=(1,0,-1)^T$, $Av_4=(0,0,1)^T$ and $Av_5=(0,1,0)^T$. Assumption \ref{as:SumTheorem3} holds with $G^1 = \{1,3\}$, $G^2 = \{2,4\}$, $G^3 = \{5\}$ and $\sigma(1)=1$, $\sigma(2)=3$, $\sigma(3)=2$, $\sigma(4)=4$, $\sigma(5)=5$. 
		To verify that condition $(ii)$ of Theorem \ref{thm:SumTheorem3} holds, let $x \in \X$, and first consider $y\in \partial_1(K_A+x) \cap \X $, where $\partial_1(K_A+x) \cap \X = \{ w \in \X \:|\: x_1 = w_1, x_2 \geq w_2, x_2 + x_3 \geq w_2 + w_3, x_4 \le w_4\}$. Given that $\inn{A_{1\bullet},v_1}=\inn{A_{1\bullet},v_3}=1$, we need to check that $\rate_1(x) + \rate_3(x) \le \Breve{\rate}_1(y) + \Breve{\rate}_3(y)$. Since $y \in \partial_1(K_A+x) \cap \X$, then $\rate_1(x) = \kappa_1 x_1 =\kappa_1 y_1 = \Breve{\kappa}_1 y_1 = \Breve{\rate}_1(y)$ and $\rate_3(x) = \kappa_3 x_1 =\kappa_3 y_1 = \Breve{\kappa}_3 y_1 = \Breve{\rate}_3(y)$, and so the desired inequality holds with equality. Secondly, consider $y \in \partial_2(K_A+x) \cap \X = \{ w \in \X \:|\: x_1 \geq w_1, x_2 = w_2, x_2 + x_3 \geq w_2 + w_3, x_4\le w_4\}$. Given that $\inn{A_{2\bullet},v_2}=\inn{A_{2\bullet},v_5}=1$ and $\inn{A_{2\bullet},v_1}=-1$, we need to check that $\rate_1(x) \ge \Breve{\rate}_1(y)$, $\rate_2(x) \le \Breve{\rate}_2(y)$ and $\rate_5(x) \le \Breve{\rate}_5(y)$. Since $y \in \partial_2(K_A+x) \cap \X$, then $\rate_1(x) = \kappa_1 x_1 \geq \kappa_1 y_1 = \Breve{\kappa}_1 y_1 = \Breve{\rate}_1(y)$, $\rate_2(x) = \kappa_2 x_2 = \kappa_2 y_2 = \Breve{\kappa}_2 y_2 = \Breve{\rate}_2(y)$ and $\rate_5(x) = \kappa_5 x_2 = \kappa_5 y_2 \leq \Breve{\kappa}_5 y_2 = \Breve{\rate}_5(y)$, and so the desired inequality holds. Lastly, consider $y \in \partial_3(K_A+x) \cap \X = \{ w \in \X \:|\: x_1 \geq w_1, x_2 \geq w_2, x_2 + x_3 = w_2 + w_3, x_4\le w_4\}$. Given that $\inn{A_{1\bullet},v_1}=\inn{A_{1\bullet},v_3}=-1$ and $\inn{A_{1\bullet},v_2}=\inn{A_{1\bullet},v_4}=1$, we need to check that $\rate_1(x) +\rate_3(x) \ge \Breve{\rate}_1(y) + \Breve{\rate}_3(y)$ and $\rate_2(x) + \rate_4(x) \le \Breve{\rate}_2(y) + \Breve{\rate}_4(y)$. For $y\in \partial_3(K_A+x) \cap \X$, since $\kappa_2 < \kappa_4$, we have that $\rate_1(x) + \rate_3(x) = (\kappa_1 + \kappa_3) x_1 \ge (\kappa_1 + \kappa_3) y_1 = (\Breve{\kappa}_1 + \Breve{\kappa}_3) y_1  = \Breve{\rate}_1(y) + \Breve{\rate}_3(y)$ and $\rate_2(x) + \rate_4(x) = \kappa_2 x_2 + \kappa_4 x_3 = (\kappa_4-\kappa_2) x_3 + \kappa_2 (x_2 + x_3) \le (\kappa_4-\kappa_2) y_3 + \kappa_2 (y_2 + y_3) = \Breve{\kappa}_2 y_2 + \Breve{\kappa}_4 y_3 = \Breve{\rate}_2(y) + \Breve{\rate}_4(y)$. Thus, the conditions of Theorem \ref{thm:SumTheorem3} are satisfied. In particular, we can conclude when $\kappa_2 < \kappa_4$, that it takes less time to get to $(0,0,0,\mathrm{S_{tot}})$ from $(\mathrm{S_{tot}},0,0,0)$ if reaction ${\large \textcircled{\small 5}}$ is added to the system without that reaction\footnote{The system without reaction ${\large \textcircled{\small 5}}$ can be obtained from the system with reaction ${\large \textcircled{\small 5}}$ by setting $\kappa_5 =0$. While strictly speaking a zero rate constant is not within our definition of mass action kinetics, our theory does cover propensity functions with such a zero rate constant.}.
		\subsection{An algorithm for coupled stochastic simulation}
		\label{sec:AlgorithmStochasticSimulation}
		
		We now provide an algorithm for stochastic simulation of the coupled continuous-time Markov chains $X$ and $\Breve{X}$ under the conditions of Theorems \ref{thm:MainResult}, \ref{thm:InnerProductTheorem}, \ref{thm:SumTheorem2} or \ref{thm:SumTheorem3}, when the transitions rates are bounded on the state space, i.e., when \eqref{eq:uniformization_condition_intensity} holds. 
		
		\SetKwComment{Comment}{// }{ //}
		
		\begin{algorithm}[H]
			
			\SetKwFunction{TransitionDTMC}{TransitionDTMC}
			\SetKwFunction{Exponential}{Exponential}
			\SetKwFunction{Uniform}{Uniform}
			
			\TitleOfAlgo{Stochastic simulation for coupled continuous-time Markov chains $X$ and $\Breve{X}$.}

			\KwData{Integer $n \geq 1$, real $T >0$, set $\X \subseteq \Z_+^d$, vectors $v_1,\ldots,v_n$ in $\Z^d \setminus \{0\}$, $\initialx,\initialxbreve$ in $\X$, functions $\rate=(\rate_1, \dots,\rate_n)$ and  $\Breve{\rate}= (\Breve{\rate}_1, \dots,\Breve{\rate}_n)$ and integer $a \in \{1,2,3,4\}$ to indicate which theorem is invoked (\ref{thm:MainResult}, \ref{thm:InnerProductTheorem}, \ref{thm:SumTheorem2} or \ref{thm:SumTheorem3}).}
			\KwResult{Sample of initial time and subsequent potential jump times $T_0,T_1,\ldots,T_{N}$ and associated states $X(T_0), X(T_1),\ldots,X(T_{N})$ and $\Breve{X}(T_0),\Breve{X}(T_1),$ $\ldots,\Breve{X}(T_{N})$ for the continuous-time Markov chains $X$ and $\Breve{X}$ in the time interval $[0,T]$.}
			$\lambda \gets 1 +  n\max\left\{\sup_{x \in \X} \sum_{j=1}^n\rate_j(x), \sup_{x \in \X} \sum_{j=1}^n\Breve{\rate}_{j}(x) \right\}$ \;
			$K \gets 0$\;
			$T_0 \gets 0$\;
			\While{$T_K \leq T$}{
				$T_{K+1} \gets T_K +$ \Exponential{$\lambda$}\;
				$K \gets K+1$\;
			}
			$Y_0,\Breve{Y}_0 \gets \initialx,\initialxbreve$\;
			$X(T_0),\Breve{X}(T_0) \gets \initialx,\initialxbreve$\;
			$N \gets K-1$\;
			\If{$N \geq 1$}{
				\For{$k\leftarrow 0$ \KwTo $N-1$}{
					$U \gets$ \Uniform{$[0,1]$}\;
					$Y_{k+1}, \Breve{Y}_{k+1} \gets$ \TransitionDTMC{$\rate,\lambda,Y_k,U,a$}, \TransitionDTMC{$\Breve{\rate},\lambda,\Breve{Y}_k,U,a$}\;
					$X(T_{k+1}),\Breve{X}(T_{k+1}) \gets Y_{k+1}, \Breve{Y}_{k+1}$\;
				}  
			}
		\end{algorithm}
		
		The random variables $T_1,\ldots,T_{N}$ are called potential jump times because it could be that $X(T_k) = X(T_{k+1})$ or $\Breve{X}(T_k) = \Breve{X}(T_{k+1})$ for some $0 \leq k \leq N-1$. Letting $T_{N+1}:= T$, the trajectories of $X$ are given by $X(t) = X(T_k)$ for $T_k \leq t < T_{k+1}$, $0 \leq k \leq N$, and similarly for the trajectories of $\Breve{X}$.
		
		The function \texttt{TransitionDTMC} can be found below. This function is meant to replicate $\Phi_{\lambda}(x,u)$ in \eqref{eq:DefinitionPhiLambda} for the case of Theorems \ref{thm:MainResult} and \ref{thm:InnerProductTheorem}, and $\PPhi_{\lambda}(x,u)$ in \eqref{eq:DefinitionPhiLambdaBIS2A} for the case of Theorems \ref{thm:SumTheorem2} and \ref{thm:SumTheorem3}.
		
		\SetKwProg{Fn}{Function}{}{}
		
		\begin{algorithm}[H]
			\Fn{\TransitionDTMC{$\rate,\lambda,x,u,a$}}{
				\KwData{Integer $n \geq 1$, set $\X \subseteq \Z_+^d$, vectors $v_1,\ldots,v_n$ in $\Z^d \setminus \{0\}$.}
				\KwIn{Function $\rate=(\rate_1, \dots,\rate_n)$, $\lambda > 0, x \in \X, u \in [0,1]$. Integer $a \in \{1,2,3,4\}$ to indicate which theorem is invoked (\ref{thm:MainResult}, \ref{thm:InnerProductTheorem}, \ref{thm:SumTheorem2} or \ref{thm:SumTheorem3}). For the case of Theorem \ref{thm:SumTheorem2} or \ref{thm:SumTheorem3}, include partition $\{G^1, \dots, G^s\}$ and bijection $\sigma$.}
				\KwOut{State $x+v \in \X$.}
				$v \gets 0$\;
				\If(\tcp*[h]{The case of Theorem \ref{thm:MainResult} or \ref{thm:InnerProductTheorem}.}){$a \in \{1,2\}$}{
					\For{$j \gets 1$ \KwTo $n$}{
						\If{$\frac{j-1}{n} \leq u < \frac{j-1}{n} + \frac{\rate_j(x)}{\lambda}$}{
							$v  \gets v_j$\;
						}
					}
				}
				\Else(\tcp*[h]{The case of Theorem \ref{thm:SumTheorem2} or \ref{thm:SumTheorem3}.}){
					$p_0 \gets 0$\;
					\For{$k \gets 1$ \KwTo $s$}{
						$p_k \gets p_{k-1} + |G^k|$\;
						\For{$q \gets p_{k-1}+1$ \KwTo $p_k$}{
							\If{$\frac{p_{k-1}}{n}+\sum_{\ell=p_{k-1}+1}^{q-1}\frac{\rate_{\sigma(\ell)}(x)}{\lambda} \leq u < \frac{p_{k-1}}{n} +\sum_{\ell=p_{k-1}+1}^{q}\frac{\rate_{\sigma(\ell)}(x)}{\lambda}$}{
								$v  \gets v_{\sigma(q)}$\;
							}
						}
					}
				}
				\KwRet{x+v}
			}
		\end{algorithm}
		
		\begin{remark}
			The above algorithm can be adapted to provide simultaneous stochastic simulation for $X$ and $\Breve{X}$ when transition rates are not bounded on the state space, by applying the algorithm on a sequence of bounded sets, which expand to the whole state space. This employs a sequence of successively defined stopping times $\tau_0=0$, $\tau_{\l}=\inf\{ t \geq \tau_{\l-1} \:|\: X(t) \notin C_{\l} \text{ or } \Breve{X}(t) \notin C_{\l} \}$, $\l=1,2,\dots$, where the $C_{\l}$ are compact, $C_{\l} \subseteq C_{\l+1}$ for $\l=1,2,\dots$ and $\cup_{\l=1}^{\infty}C_{\l}=\X$. The simulation of the pair $(X,\Breve{X})(t)$ for $\tau_{\l} \leq t < \tau_{\l+1}$ uses the above algorithm on $C_{\l}$ for $\l=1,2,\dots$.
		\end{remark}

		\vfill\break

	\end{appendices}

\end{document}